%% file: lefschetzsym.tex
\documentclass[11pt]{amsart}
\usepackage{eucal}
\usepackage[all]{xy}
\usepackage{amsfonts,amssymb}
\usepackage{epsfig}
\usepackage[usenames]{color}
\usepackage{float} 
\xyoption{arc}

\newtheorem{theorem}{Theorem}[section]
\newtheorem{lemma}[theorem]{Lemma}

\newtheorem{example}[theorem]{Example}
\newtheorem{corollary}[theorem]{Corollary}
\newtheorem{defn}[theorem]{Definition}

\newtheorem{fig}[theorem]{Figure}

\newcommand{\C}{\mathbb C}
\newcommand{\R}{\mathbb R}
\renewcommand{\P}{\mathbb P}
\renewcommand{\O}{\EuScript O}

\newcommand{\D}{\mathbb D}
\newcommand{\F}{\mathbb F}

\newcommand{\Z}{\mathbb Z}
\newcommand{\N}{\mathbb N}

\newcommand{\nhd}{\mbox{nhd}}
\newcommand{\nsmbox}[1]{\mbox{\rm #1}}

\title{Lefschetz fibrations and symplectic homology}
\author{Mark McLean}

\begin{document}

\begin{abstract}
We show that for each $k > 3$ there are infinitely
many finite type Stein manifolds diffeomorphic to 
Euclidean space $\R^{2k}$ which are pairwise distinct as 
symplectic manifolds. 
\end{abstract}

\maketitle

\bibliographystyle{hamsplain}


\tableofcontents

\section{Introduction}

This paper is about the symplectic topology of Stein manifolds.
If we have a symplectic manifold $(V,\omega)$, then we say
it carries a Stein structure if there exists a complex 
structure $J$ and
an exhausting (i.e. proper and bounded from below)
plurisubharmonic function $\phi : V \rightarrow \R$ such that $\omega = -dd^c \phi$,
where $d^c$ is defined by $d^c(a)(X):= da(JX)$.
The triple $(V,J,\phi)$ is called a Stein manifold.
We say that a Stein manifold is of finite type if
$\phi$ only has finitely many critical points, each of which is non-degenerate.

We define an equivalence relation $\sim$ on Stein manifolds by:
$A \sim B$ if there exists a sequence of Stein manifolds
$F_0,F_1,\cdots, F_n$ such that
\begin{enumerate}
\item $F_0 = A$ and $F_n = B$
\item $F_i$ is either symplectomorphic or Stein deformation equivalent 
to $F_{i+1}$.
(Stein deformation is defined later in Definition \ref{defn:steindeformation}.)
\end{enumerate}

The aim of this paper is the following theorem:
\begin{theorem} \label{thm:infinitelymanysteinmanifolds}
Let $k \geq 4$. 
There exists a family of finite type
Stein manifolds $X_i$ diffeomorphic to $\R^{2k}$
indexed by $i \in \N$ such that
\[i \neq j \Rightarrow X_i \nsim X_j\]
\end{theorem}

We also have the following corollary:
\begin{corollary} \label{corollary:cotangent}
Let $M$ be a compact manifold of dimension $4$ or higher.
There exists a family of finite type
Stein manifolds $X^M_i$ diffeomorphic to $T^*M$
indexed by $i \in \N$ such that
\[i \neq j \Rightarrow X^M_i \nsim X^M_j\]
\end{corollary}
We will prove this at the end of this introduction.
The results of
Paul Seidel and Ivan Smith in \cite{SS:rama}
show that there exists a finite type Stein manifold
diffeomorphic to $\R^{4k}$ but not symplectomorphic
to $\R^{4k},k \geq 2$. They show this by constructing
an affine variety which has a Lagrangian
torus which cannot be moved off itself
by a Hamiltonian isotopy. In fact they
show that none of these Stein manifolds can
be embedded in a subcritical Stein manifold.
Also, they have an argument (explained in
\cite{Seidel:biasedview}) that shows that the
symplectic homology groups of these
varieties are non-trivial.
In this paper we strengthen this result in two ways:
\begin{enumerate}
\item \label{item:allevendimensions}
we give examples in all even dimensions $\geq 8$ (not just in dimension
$4k$ where $k \geq 2$); 
\item
we also show there are countably many pairwise distinct
examples in each of these dimensions. \label{item:countablymany}
\end{enumerate}
(\ref{item:allevendimensions}) is straightforward but (\ref{item:countablymany})
is much harder and involves various new ideas.
We show in Corollary \ref{corollary:embeddinginsubcritical}
that these manifolds cannot be embedded in a subcritical Stein manifold.


There is no analogue of our result in dimension 4 because
any finite type Stein manifold diffeomorphic to $\R^4$ is
symplectomorphic to $\R^4$ (see the introduction to \cite{SS:rama}).
Having said that,
Gompf in \cite{Gompf:handlebody} constructs uncountably many non-finite type
Stein manifolds
which are homeomorphic to $\R^4$, but are pairwise not 
diffeomorphic to each other.
We hope to address the question of whether Theorem 
\ref{thm:infinitelymanysteinmanifolds}
holds in dimension 6, and whether we can distinguish the
contact boundaries of these manifolds up to contactomorphism,
in future work.

\bigskip

We will now construct an example of a family of Stein manifolds 
$(X_n)_{n \in \N}$ as in Theorem 
\ref{thm:infinitelymanysteinmanifolds} in dimension 8.
Let $V := \{x^7+y^2+z^2+w^2 = 0\} \subset \C^4$ and
consider a smooth point, say $p := (0,0,1,i) \in V$.
Let $H$ be the blowup of $\C^4$ at $p$.
Then $X := H \setminus \widetilde{V}$ is a Stein manifold where
$\widetilde{V}$ is the proper transform of $V$.
The variety $X$ is called the Kaliman modification of $(\C^4,V,p)$
(see \cite{Kaliman:eisenman}).
We will think of this modification in two stages:
\begin{enumerate}
\item Cut out the hypersurface $V$ in $\C^4$ to get $Z := \C^4 \setminus V$.
\item \label{item:kalimanblowup}
Blow up $Z$ at infinity to get $X$.
\end{enumerate}
Operation (\ref{item:kalimanblowup})
attaches a $2$-handle along a
knot which is transverse to the contact structure.
All our Stein manifolds in this family will be constructed from $X$.
If we have two Stein manifolds $A$ and $B$, then
it is possible to construct their {\it end connected sum}
$A \#_e B$ (see Theorem \ref{thm:endconnectsum}).
Roughly what we do here is join $A$ and $B$ with a $1$-handle,
and then extend the Stein structure over this handle.
Finally, $X_n := \#_{i=1}^n X$ is our family of Stein manifolds.

\subsection{Sketch of the proof of the main theorem} \label{subsection:sketchofproof}

We will now give an outline of the proof. We 
will only consider the examples $X_n$
in dimension 8 constructed above, as the higher dimensional
examples are similar.
For each Stein manifold $Y$ with 
a trivialisation of the canonical bundle, we have an integer graded
commutative $\Z/2\Z$ algebra $SH_{n+*}(Y)$ where
$SH_*(Y)$ is called symplectic homology 
\footnote{With our convention, the pair-of-pants product makes
$SH_*$ (and not $SH^*$) a unital ring.}.
The reason why we write $SH_{n+*}(Y)$ instead of $SH_*(Y)$
is because we want the unit to be in degree $0$ and not in degree $n$.
If $Y_1$ and $Y_2$ are Stein manifolds with
$H_1(Y_1) = H_1(Y_2) = 0$ and $Y_1 \sim Y_2$, then
$SH_*(Y_1) = SH_*(Y_2)$ (see \cite[Section 7]{Seidel:biasedview}).
The reason why we need $H_1(Y_1) = H_1(Y_2) = 0$ is because symplectic
homology is only known to be invariant up to exact symplectomorphism.
For each Stein manifold $Y$, we can define another invariant
$i(Y)$ which is the number of idempotents of $SH_*(Y)$
(this invariant might be infinite).
Hence all we need to do is show that for $i \neq j$,
$i(X_i) \neq i(X_j)$.
The next fact we need is that for any two Stein manifolds $Y_1$ and $Y_2$,
$SH_*(Y_1 \#_e Y_2) = SH_*(Y_1) \times SH_*(Y_2)$.
This means that $SH_*(X_n) = \prod_{i=1}^n SH_*(X)$,
and hence $i(X_n) = i(X)^n$. So all we need to do is
show that $1 < i(X) < \infty$.
If $SH_*(X) \neq 0$, $i(X) > 1$ since we have $0$ and $1$;
but since $SH_*(X)$ can a priori be infinite dimensional in each
degree, finiteness of $i(X)$ is much harder.
Most of the work in this paper involves proving $i(X) < \infty$.

For any Stein manifold $Y$, $SH_*(Y)$ is $\Z$ graded by
the Conley-Zehnder index taken with negative sign.
The group $SH_*(Y)$ has a ring structure making
$SH_{n+*}(Y)$ into a $\Z / 2\Z$ graded algebra.
This ring is also a $H_1(Y)$ graded algebra. This means
that as a vector space, it is of the form
$\bigoplus_{i \in H_1(Y)}R_i$ and if $x \in R_a$ and $y \in R_b$
then their product $xy$ is in $R_{a+b}$.
Hence idempotents in $SH_{n+*}(Y)$ are contained in
$SH_n(Y)$ and are a linear combination of elements with grading in the torsion part of $H_1(Y)$
(see Lemma \ref{lemma:finitelymanyidempotents}). The problem
is that $H_1(X)=0$.
In order to find out which elements of $SH_{n+*}(X)$ are idempotents,
we will show that $SH_*(X)$ is isomorphic as a ring to $SH_*(Z)$ where
$Z = \C^4 \setminus V$ was defined above. 
Because $Z$ is so much simpler than $X$ and $H_1(Z) \neq 0$, 
it is possible by a direct
calculation to show that $SH_{n+*}(Z)$ has finitely many idempotents.

Proving that $SH_*(X) \cong SH_*(Z)$ relies on the following theorem.
This theorem is the heart of the proof. We let
$E' \rightarrow \C$, $E'' \rightarrow \C$ be
Lefschetz fibrations, and  $F'$ (resp. $F''$) be smooth fibres of 
$E'$ (resp. $E''$). Let $F'$ and $F''$ be Stein domains with
$F''$ a holomorphic and symplectic submanifold of $F'$.

\begin{theorem} \label{theorem:fibrationaddition}
Suppose $E'$ and $E''$ satisfy the following properties:
 \begin{enumerate}
   \item \label{item:subfibration}
$E''$ is a subfibration of $E'$.
   \item \label{item:monodromysupport}
The support of all the monodromy maps of $E'$ are
contained in the interior of $E''$.
   \item \label{item:complexcurvecriterion}
Any holomorphic curve in $F'$
with boundary inside $F''$ must be contained in $F''$.
 \end{enumerate}
Then $SH_*(E') \cong SH_*(E'')$.
\end{theorem}

Remark 1: There exist Lefschetz fibrations $E'$, $E''$ 
with the above properties such that as convex symplectic manifolds,
$E'$ (resp. $E''$) is convex deformation equivalent to $X$ (resp. $Z$).
This is because we can choose an
algebraic Lefschetz fibration on $Z$ where the closures of
all the fibres pass through $p$. Then blowing up $Z$ at infinity
(operation (\ref{item:kalimanblowup}) of the Kaliman modification)
is the same as blowing up each fibre at infinity
and keeping the same monodromy.
Hence $SH_*(X) \cong SH_*(Z)$. 

Remark 2: Given varieties $X$ and $Z$ in dimension $4$ such that
$X$ is obtained from $Z$ by blowing up at infinity,
there are Lefschetz fibrations $E'$, $E''$
satisfying properties (\ref{item:subfibration}) and
(\ref{item:monodromysupport}) such that 
as convex symplectic manifolds,
$E'$ (resp. $E''$) is convex deformation equivalent to $X$ (resp. $Z$).
These do not satisfy property
(\ref{item:complexcurvecriterion}) because $E'$ is obtained from
$E''$ by filling in a boundary component of the fibres with a disc.

We will prove Theorem \ref{theorem:fibrationaddition}
in two stages. In stage (i), we
construct a $\Z / 2\Z$ graded algebra $SH_*^{\nsmbox{lef}}(E)$
related to a Lefschetz fibration $E$ and show it is
equal to symplectic homology. This is
covered in sections \ref{section:lefschetzfibrationproofs}
and \ref{section:bettercofinal}. In stage (ii),
we prove that
$SH_*^{\nsmbox{lef}}(E') \cong SH_*^{\nsmbox{lef}}(E'')$.
This is covered in section \ref{section:transferisomorphism}.
In a little more detail:

\smallskip

(i) Let $F$ be a smooth fibre of $E$ and $\D$ a disc in $\C$.
In section \ref{section:lefschetzcofinal}, we show (roughly) that
the chain complex $C$ for $SH_*(E)$ is generated by:
\begin{enumerate}
\item \label{item:criticalmorsepoints}
critical points of some Morse function on $E$;
\item \label{item:fixedmonodromypoints}
two copies of fixed points of iterates of the
monodromy map around a large circle;
\item \label{item:productorbits}
pairs $(\Gamma, \gamma)$ where $\Gamma$ is a Reeb orbit
on the boundary of $F$ and $\gamma$ is either a Reeb
orbit of $\partial \D$ or a fixed point in the interior of $\D$.
\end{enumerate}

\begin{figure}[b]
\end{figure}
\begin{figure}[b]
\end{figure}

This is done in almost exactly the same way as the proof of the
K\"unneth formula for symplectic homology \cite{Oancea:kunneth}.
The differential as usual
involves counting cylinders connecting the orbits and
satisfying the perturbed Cauchy -Riemann equations.
The orbits in (\ref{item:criticalmorsepoints}) and 
(\ref{item:fixedmonodromypoints}) actually form a subcomplex 
$C^{\nsmbox{lef}}$, and we define
Lefschetz symplectic homology $SH_*^{\nsmbox{lef}}(E)$
to be the homology of this subcomplex.
Let $\Phi : SH_k^{\nsmbox{lef}}(E) \rightarrow SH_k(E)$ be the map
corresponding to the inclusion $C^{\nsmbox{lef}} \hookrightarrow C$.
We need to show that $\Phi$ is an isomorphism.
In order to do this we will show that the elements of $C$ in (\ref{item:productorbits})
have very high index.
The orbits $\gamma$ in (\ref{item:productorbits}) actually correspond to orbits
of some Hamiltonian on the complex plane $\C$. We can ensure that
this Hamiltonian has orbits of index as large as we like.
In particular, we can assume that the index of $\gamma$ is so large that
the index of the pair $(\gamma, \Gamma)$ is greater than any number
we like (because the index of $(\gamma, \Gamma)$ is the sum
$\nsmbox{index}(\gamma) + \nsmbox{index}(\Gamma)$).
This means for any $k \in \Z$, the map $SH_k^{\nsmbox{lef}}(E) \rightarrow SH_k(E)$
is an isomorphism if we ensure the indices of the orbits $(\gamma, \Gamma)$
are all greater than $k+1$.
Hence $SH_*(E) \cong SH_*^{\nsmbox{lef}}(E)$.

(ii)
Let $C'$ (resp. $C''$) be the chain complex for the group
$SH_*^{\nsmbox{lef}}(E')$ (resp. $SH_*^{\nsmbox{lef}}(E'')$).
The fibration $E' \setminus E''$ is a trivial fibration
$\D \times W$.
We have a short exact sequence
$0 \rightarrow B \rightarrow C' \rightarrow C'' \rightarrow 0$
where $B$ is generated by orbits of the form $(\gamma,\Gamma)$
in $\D \times W$.
The orbit $\Gamma$ is a critical point of some Morse function on $W$
and $\gamma$  is either a Reeb
orbit of $\partial \D$ or a fixed point in the interior of $\D$.
In this case the homology of the chain complex $B$ is actually a product
$SH_*(\D) \otimes X$ where $X$ is the relative cohomology 
group $H^{n-*}(F', F'')$.
Because $SH_*(\D)=0$, we have that $H_*(B)=0$ which
gives us our isomorpism between $SH_*^{\nsmbox{lef}}(E')$ 
and $SH_*^{\nsmbox{lef}}(E'')$.
Property
(\ref{item:complexcurvecriterion}) in Theorem
\ref{theorem:fibrationaddition} is needed here to ensure
that the above exact sequence exists. If we didn't have
this property, then there would be some spectral
sequence from $SH_*^{\nsmbox{lef}}(E'')$ (with an extra grading
coming from the $H_1(E'')$ classes of these orbits)
to $SH_*^{\nsmbox{lef}}(E')$.

Lefschetz symplectic homology was partially inspired
by Paul Seidel's Hochshild homology conjectures 
\cite{Seidel:hochschildhomology},
which also relate symplectic homology to
Lefschetz fibrations. His conjectures would in particular
prove Theorem \ref{theorem:fibrationaddition}.

\proof of Corollary \ref{corollary:cotangent}.
There is a standard Stein structure on $T^*M$
such that $SH_0(M)$ is a non-trivial
finite dimensional $\Z / 2\Z$ vector space.
By Lemma \ref{lemma:finitelymanyidempotents} this
means that $i(T^*M) < \infty$. Also $0 \in SH_*(T^*M)$
is an idempotent which means that $0<i(T^*M)$.
We let $X_i$ be defined as in the proof of
the main theorem \ref{thm:infinitelymanysteinmanifolds}.
We define 
\[X^M_i := T^*M \#_e X_i\]
Then $i(X^M_i) = i(T^*M)i(X_i)$.
These numbers are all different as $0 < i(T^*M) < \infty$
and $i(X_i) \neq i(X_j)$ for $i \neq j$.
\qed

\bigskip

{\bf Acknowledgements:}
I would like to thank my supervisor Ivan Smith for checking this paper
and for giving many useful suggestions. I would also like to thank
Paul Seidel, Kai Cieliebak, Alexandru Oancea,
Fr\'ed\'eric Bourgeois, Domanic Joyce, Burt Totaro,
Jonny Evans and Jack Waldron
for giving useful comments.

\subsection{Notation}

Throughout this paper we use the following notation:
\begin{enumerate}
\item $M,M',M'',\dots$ are manifolds (with or without boundary).
\item $\partial M$ is the boundary of $M$.
\item $(E,\pi),(E',\pi'),(E'',\pi'')$ are exact Lefschetz fibrations
(See Definition \ref{defn:lefschetzfibration}).
\item If we have some data $X$ associated to $M$ (resp. $E$), then \\
 $X,X',X'',\dots$ are data
 associated to $M,M',M'',\dots$ \\ (resp.  $E,E',E'',\dots$).
 For instance $\partial M''$ is the boundary of $M''$.

\item $\omega$ is a symplectic form on $M$ or $E$.
\item $\theta$ is a $1$-form such that
$d\theta=\omega$.
\item $(M,\theta)$ is an exact symplectic manifold.
\item $J$ is an almost complex structure compatible with $\omega$.
\item If $(M,\theta)$ is a compact convex symplectic manifold
(See Definition \ref{defn:compactconvexsymplecticmanifold}), 
then $(\widehat{M},\theta)$ is the completion of $(M,\theta)$
(Lemma \ref{lemma:completion}).
Similarly by Definition \ref{defn:halfconvexlefschetzcompletion},
$(E,\pi)$ can be completed to
$(\widehat{E},\pi)$ (we leave $\pi$ and $\theta$ as they are by
abuse of notation).
\item $(M,\theta_t)$ is a convex symplectic or Stein deformation.
\item $F$ will denote a smooth fibre of $(E,\pi)$.
\item If we have some subset $A$ of a topological space, then we will let
$\nhd(A)$ be some open neighbourhood of $A$.

\end{enumerate}

\bigskip

\section{Background}

\subsection{Stein manifolds}

We will define Stein manifolds as in \cite{SS:rama}.
We let $M$ be a manifold and $\theta$ a $1$-form where
$\omega := d\theta$ is a symplectic form.

\begin{defn} \label{defn:compactconvexsymplecticmanifold}
$(M,\theta)$ is called a {\bf compact convex symplectic manifold}
if $M$ is a compact manifold with boundary and the
$\omega$-dual of $\theta$ is transverse to $\partial M$ and pointing
outwards.
A {\bf compact convex symplectic deformation} is a family of
compact convex symplectic manifolds $(M,\theta_t)$ parameterized
by $t \in [0,1]$.
We will let $\lambda$ be the vector field which is $\omega$-dual
to $\theta$.
\end{defn}

Usually, a compact convex symplectic manifold is called
a convex symplectic domain. We have a natural contact form
$\theta |_{\partial M}$ on $\partial M$, and hence we call
this the contact boundary of $M$.

\begin{defn} \label{defn:convexsymplecticmanifold}
Let $M$ be a manifold without boundary.
We say that $(M,\theta)$ is a
{\bf convex symplectic manifold} if there exist
constants $c_1 < c_2 < \cdots$ tending to infinity
and an exhausting
function  $\phi : M \rightarrow \R$ such that
$(\{\phi \leq c_i\},\theta)$ is a
compact convex symplectic manifold for each $i$.
Exhausting here means proper and bounded from below.
If the flow of $\lambda$ exists for all
positive time, then $(M,\theta)$ is called {\bf complete}.
If there exists a constant $c>0$ such that for all $x \geq c$,
$(\{\phi \leq c\},\theta)$ is a
compact convex symplectic manifold, then we say that
$(M,\theta)$ is of {\bf finite type}.
\end{defn}

\begin{defn} \label{defn:convexsymplecticdeformation}

Let $(M,\theta_t)$ be a smooth family of convex symplectic manifolds
with exhausting functions $\phi_t$. Suppose that for each $t \in [0,1]$,
there are constants 
$c_1 < c_2 < \cdots$ tending to infinity and an $\epsilon>0$
such that for each $s$ in $(t-\epsilon, t+\epsilon)$ and
$i \in \N$,
$(\{\phi_s \leq c_i\},\theta_s)$ is a
compact convex symplectic manifold. 
Then $(M,\theta_t)$ is called a {\bf convex symplectic deformation}.
\end{defn}
The constants $c_1 < c_2 < \cdots$
mentioned in this definition
depend on $t$ but not necessarily in a continuous way.
The nice feature of convex symplectic manifolds is that we have some control
over how they behave near infinity. That is, the
level set $\phi^{-1}(c_k)$ is a contact manifold for all $k$. 
Also note that if we have a convex symplectic manifold
of finite type then it has a cylindrical end. That is,
there exists a manifold $N$ with a contact form $\alpha$ such that
at infinity, $M$ is symplectomorphic to $(N \times [1,\infty), d(r \alpha))$
($r$ is a coordinate on $[1,\infty)$).

\begin{lemma} \label{lemma:completion}
A compact convex symplectic manifold $M$ can
be completed to a finite type complete 
convex symplectic manifold 
$(\widehat{M},\theta)$.
\end{lemma}
This is explained for instance in 
\cite[section 1.1]{Viterbo:functorsandcomputations}. The proof
basically involves gluing a cylindrical end onto $\partial M$.
Let $(M,\theta)$ be a complete convex symplectic manifold.
Let $(M',\theta')$ be a compact convex symplectic manifold
which is a codimension $0$ exact submanifold
of $(M,\theta)$ (i.e. $\theta|_{M'} = \theta' + dR$ for
some smooth function $R$ on $M'$).
\begin{lemma} \label{lemma:completionsubmanifold}
We can extend the embedding $M' \hookrightarrow M$ to an
embedding $\widehat{M'} \hookrightarrow M$.
\end{lemma}
\proof
There exists a function $R : M' \rightarrow \R$ such that
$\theta' = \theta + dR$. We can extend $R$ over the whole
of $M$ such that $R=0$ outside some compact subset $K$
containing $M'$.
Let $\lambda'$ be the $\omega$-dual
of $\theta'$. Let $F_t : M \rightarrow M$  be the flow
of $\lambda'$. This exists for all time because $M$ is complete
and $\lambda'= \lambda$ outside $K$.
We have an embedding $\Phi : (\partial M') \times [1,\infty) \rightarrow M$
defined by $\Phi(a,t) = F_{\log{t}}(a)$. This attaches a cylindrical
end to $M'$ inside $M$, hence we have an exact embedding
$\widehat{M'} \rightarrow M$ extending the embedding of $M'$.
\qed

\begin{defn} \label{defn:steinmanifold}
A {\bf Stein manifold} $(M,J,\phi)$ is a complex manifold \\
$(M,J)$ with an exhausting plurisubharmonic function
$\phi : M \rightarrow \R$ (i.e. $\phi$ is proper and
bounded from
below and $-dd^c(\phi) > 0$ where $d^c = J^*d$).
A Stein manifold is called {\bf subcritical} if $\phi$
is a Morse function with critical points of index 
$< \frac{1}{2}\nsmbox{dim}_{\R}M$.
A manifold of the form
$\phi^{-1}((-\infty,c])$ is called a {\bf Stein domain}.
\end{defn}

We can perturb $\phi$ so that it becomes a Morse function.
From now on, if we are dealing with a Stein manifold, we will
always assume that $\phi$ is a Morse function.
The index of a critical point of $\phi$ is always less
than or equal to $\frac{1}{2} \nsmbox{dim}_{\R}M$.
This is because the unstable manifolds of these critical points
are isotropic submanifolds of $M$.
Note that the definition of a Stein manifold in  
\cite[Section 2]{Eliashberg:symplecticgeometryofplushfns}
is that it is a
closed holomorphic submanifold of $\C^N$ for some $N$. This
has an exhausting plurisubharmonic function $|z|^2$.
An important example of a subcritical Stein manifold is
$(\C^n,i,|z|^2)$.
The Stein manifold $(M,J,\phi)$ is 
a convex symplectic manifold $(M,\theta := -d^c \phi)$.
Note that
$\lambda := \nabla \phi$ where
$\nabla$ is taken with respect to the metric $\omega(\cdot,J(\cdot))$.
It is easy to see that $\nabla \phi$ is a Liouville vector field
transverse to a regular level set of $\phi$ and pointing outwards.
We call a Stein manifold complete or of finite type if the associated
convex symplectic structure is complete or of finite type respectively.

\begin{defn} \label{defn:steindeformation}
If $(J_t,\phi_t)$ is a smooth family of Stein
structures on $M$, then it is called a {\bf Stein deformation}
if the function $(t,x) \longrightarrow \phi_t(x)$ is
proper and for each $t \in [0,1]$, there exists $c_1 < c_2 < \dots$ 
tending to infinity and an $\epsilon>0$
such that for each $s$ in $(t-\epsilon, t+\epsilon)$ 
we have that
$c_k$ is a regular value of $\phi_s$.
This induces a corresponding convex symplectic deformation.
\end{defn}

\begin{example} \label{example:algebraicsteinmanifold}
An affine algebraic subvariety $M$ of $\C^N$ admits a Stein
structure. This is because it has a natural embedding
in $\C^N$, so we can restrict the plurisubharmonic function
$\|z\|^2$ to this variety to make it into a Stein manifold.
We can also use the following method to find a plurisubharmonic
function on $M$.
We first compactify $M$ by finding a projective variety $X$
with complex structure $i$
and an effective ample divisor $D$ such that $M = X \setminus D$
(for instance we can embed $M$ in $\C^N \subset \P^N$
and then let $X$ be the closure of $M$ in $\P^N$).
There exists an ample line bundle $E \longrightarrow X$
associated to the divisor $D$. 
Choose a holomorphic section $s$ of $E$ such that
$D = s^{-1}(0)$. Then ampleness means that we can
choose a metric $\|.\|$ such that its curvature form
$\omega := iF_{\nabla}$ is a positive $(1,1)$-form.
Hence we have a Stein structure 
\[(M := X \setminus D, J:=i, \phi := -\nsmbox{log}\|s\|) \]
Note that by \cite[Lemma 8]{SS:rama}, 
this is
of finite type.
\end{example}

For affine algebraic varieties we need a uniqueness theorem.
This basically comes from the text following \cite[Lemma 4.4]{Seidel:biasedview}.
The problem is that this text only deals with the case
when our compactification divisor $D$ is a normal crossing divisor.
But the methods used carry over to general divisors by using
Hironaka's resolution theorem. In fact, the paper
\cite{Seidel:biasedview} uses Hironaka's resolution theorem in the
same way as we will use it. We will repeat the argument here.
Let $M$ be an algebraic variety. 
\begin{lemma} \label{lemma:algebraicuniqueness}
If $h_1$ and $h_2$
are Stein functions on $M$ constructed as in the previous
example, then $(1-t)h_1 + th_2$ is a Stein deformation.
\end{lemma}
\proof
If we can prove that all the critical points of $(1-t)h_1 + th_2$
stay inside a compact set for all $t \in [0,1]$, then this will prove that it is a Stein deformation.
The way we do this is to look at this function in a neighbourhood of
some compactification divisor of $M$.
For $i=1,2$, let $X_i$ be projective varieties compactifying $M$. 
Let $E_i$ be ample line bundles on $X_i$ whose associated divisor
is effective and ample and has support equal to $X_i \setminus M$.
Suppose $h_i$ is equal to $-\nsmbox{log}\|s_i\|_i$ where $s_i$
is a section of the line bundle $E_i$ such that $s_i^{-1}(0)$ has
support equal to $X_i \setminus M$ and $\|\cdot\|_i$ is a metric
on $E_i$.

By Hironaka's resolution theorem, we have that there exists a compactification
$X$ by a normal crossing divisor $D$ such that it dominates $X_1$ and $X_2$ (i.e there
is a surjective morphism $X \rightarrow X_i$ such that away from $D$ and $D_i$ it is an isomorphism).
We pull back the line bundles, sections and metrics to $X$ and Stein functions via these morphisms.
We now just apply \cite[Lemma 4.4]{Seidel:biasedview}, and this gives us our result. 
\qed

\bigskip

The following operation constructs a new Stein manifold from two
old ones. This is used to construct our infinite family of Stein
manifolds. We will let $(M,J,\phi),(M',J',\phi')$ be complete
finite type Stein manifolds. Because these manifolds are complete
and of finite type, they are the completions of compact
convex symplectic manifolds $N$, $N'$ respectively.
In fact $N = \{\phi \leq R\}, N' = \{\phi' \leq R\}$ for some
arbitrarily large $R$. Let $p$ (resp. $p'$) be a point in $\partial N$
(resp. $\partial N'$). The following theorem is proved in greater
generality in \cite{Eliashberg:steintopology} and
\cite[Theorem 9.4]{CieliebakEliashberg:symplecticgeomofsteinmflds}.

\begin{theorem} \label{thm:endconnectsum}
There exists a connected finite type Stein manifold \\ $(M'',J'',\phi'')$
such that $N'' := \{\phi'' \leq R\}$ is biholomorphic to the disjoint
union of $N$ and $N'$ with $\phi''|_N = \phi$ on $N$ and 
$\phi''|_{N'} = \phi'$ on $N'$. Also,  the only critical point
of $\phi''$ outside $N''$ has index $1$.
\end{theorem}

In this theorem, what we are doing is joining $N$ and $N'$ with
a $1$-handle and then extending the Stein structure over this
handle, and then completing this manifold so that it becomes
a Stein manifold. The Stein manifold $M''$ is called the
{\bf end connect sum} of $M$ and $M'$, and we define
$M \#_e M'$ as this end connected sum.
If $M$ and $M'$ are Stein manifolds diffeomorphic to $\C^n$, then
$M \#_e M'$ is also diffeomorphic to $\C^n$.

Remark: There is an example due to the author of a
non-finite type Stein manifold which is not
equivalent to any finite type Stein manifold.
This is described by Seidel in \cite[Section 7]{Seidel:biasedview},
and is constructed as an infinite end-connect-sum.

\subsection{Lefschetz fibrations} \label{section:lefschetzfibrations}

Throughout this section we will let $E$ be a compact manifold with corners 
whose boundary is the union of two faces $\partial_h E$ and $\partial_v E$ meeting
in a codimension $2$ corner. We will also assume that $\Omega$ is
a $2$-form on $E$ and $\Theta$ a $1$-form satisfying $d\Theta = \Omega$.
We let
$S$ be a surface with boundary.
Let $\pi : E \rightarrow S$ be a smooth map with only finitely many
critical points  (i.e. points where $d\pi$ is not surjective).
Let $E^{\nsmbox{crit}} \subset E$ 
be the set of critical points of $\pi$ and
$S^{\nsmbox{crit}} \subset S$ the set of critical values of $\pi$.
For $s \in S$, let $E_s$ be the fibre $\pi^{-1}(s)$.

\begin{defn} \label{defn:symplecticfibration}
If for every $s \in S$ we have that $\Omega$ is a symplectic form
on $E_s \setminus E^{\nsmbox{crit}}$ then we say that
$\Omega$ {\bf is compatible with} $\pi$.
\end{defn}

Note that if $\Omega$ is compatible with $\pi$ then there is a
natural connection
(defined away from the critical points) for $\pi$
defined by the horizontal plane distribution which is
$\Omega$-orthogonal to each fibre.
If parallel transport along some path in the base
is well defined then it is an exact symplectomorphism
(an exact symplectomorphism is a diffeomorphism $\Phi$
between two symplectic manifolds $(M_1,d\theta_1)$
and $(M_2,d\theta_2)$ such that $\Phi^*\theta_2 = \theta_1 + dG$
where $G$ is a smooth function on $M_1$).
From now on we will assume that $\Omega$ is compatible with $\pi$.
We deal with Lefschetz fibrations as 
defined in \cite{Seidel:longexactsequence}.
We define $F$ to be some smooth fibre of $\pi$.

\begin{defn} \label{defn:lefschetzfibration}
$(E,\pi)$ is an {\bf exact Lefschetz fibration} if:

\begin{enumerate}

\item $\pi : E \longrightarrow S$
is a proper map with $\partial_v E = \pi^{-1}(\partial S)$ and
such that $\pi|_{\partial_v V} : \partial_v E \rightarrow \partial S$ is a smooth
fibre bundle. Also there is a neighbourhood  $N$ of $\partial_hE$ such that
$\pi|_N : N \rightarrow S$ is a product fibration $S \times \nhd(\partial F)$ where
$\Omega|_N$ and $\Theta|_N$ are pullbacks from the second factor of this product.

\item 
There is an integrable complex structure $J_0$ (resp. $j_0$) defined on some neighbourhood
of $E^{\nsmbox{crit}}$ (resp. $S^{\nsmbox{crit}}$) such that
$\pi$ is $(J_0,j_0)$
holomorphic near $E^{crit}$. 
At any critical point, the Hessian $D^2\pi$
is nondegenerate as a complex quadratic form. 
We also assume that there is at most one critical
point in each fibre.

\item $\Omega$ is a
K\"ahler form for $J_0$ near $E^{crit}$.

\end{enumerate}

\end{defn}

Sometimes we will need to define a Lefschetz fibration without
boundary. This is defined in the same way as an exact Lefschetz fibration
except that $E$, the fibre $F$ and the base $S$ are open manifolds without boundary.
We replace ``neighbourhood of $\partial_h E$'' in the above definition 
with an open set whose complement is relatively compact when
restricted to each fibre. We also replace ``$\partial_v E$''
with $\pi^{-1}(S \setminus K)$ where $K$ is a compact set in $S$.
Also $\pi$ is obviously no longer a proper map,
and we assume that the set of critical points is compact.
From now on we will let $(E,\pi)$ be an exact Lefschetz fibration.

\begin{lemma} \label{lemma:lefschetzsymplecticstructure} \cite[Lemma 1.5]{Seidel:longexactsequence}
If $\beta$ is a symplectic form on $S$ then
$\omega := \Omega + N\pi^*\beta$ is a symplectic form on $E$ for $N$
sufficiently large.
\end{lemma}

We really want our Lefschetz fibrations to be described as
finite type convex symplectic manifolds.

\begin{defn} \label{defn:halfconvexlefschetzfibration}
A {\bf compact convex Lefschetz fibration} is an exact Lefschetz fibration
$(E,\pi)$ such that $(F,\Theta|_F)$ is a compact convex symplectic manifold.
A {\bf compact convex Lefschetz deformation} is a smooth family
of compact convex Lefschetz fibrations parameterized by $[0,1]$.
\end{defn}

Note that by the triviality condition at infinity, all smooth
fibres of $\pi$ are
compact convex symplectic manifolds as long as the base $S$
is connected. From now on we will assume that $(E,\pi)$ is
a compact convex Lefschetz fibration.

\begin{theorem} \label{theorem:halfconvexstructure}
Let the base $S$ be a compact convex symplectic manifold $(S,\theta_S)$. 
There exists a constant $K>0$ such that
for all $k \geq K$ we have:
$\omega := \Omega + k\pi^{*}(\omega_S)$ is a symplectic form,
and the $\omega$-dual $\lambda$
of $\Theta + k\pi^{*}\theta_S$ is transverse
to $\partial E$ and pointing outwards.
\end{theorem}
(The proof is given in section \ref{section:lefschetzfibrationproofs}.)
Note that this theorem also implies that if we have a compact
convex Lefschetz deformation, then we have a corresponding
compact convex symplectic deformation because we can smooth the codimension $2$ corners
slightly. 

If we have a compact convex Lefschetz fibration, then we wish to
extend the Lefschetz fibration structure over the completion $\widehat E$ of $E$.
Here is how we naturally complete $(E,\pi)$: 
The horizontal boundary is a product
$\partial F \times S$. We can add a cylindrical
end $G:=(\partial F \times [1,\infty)) \times S$
to this in the usual way (i.e. as in Lemma \ref{lemma:completion}), extending
$\Theta$ over this cylindrical end by the 
$1$-form $r(\Theta|_{\partial F})$ where $r$ is the coordinate
for $[1,\infty)$. Let $E_1$ be the resulting
manifold. We also extend the map $\pi$ over $E_1$
by letting
$\pi|_G : G \rightarrow S$ be the natural projection.
This ensures that $\pi$ is compatible with
the natural symplectic form on $E_1$ defined as
in Lemma \ref{lemma:lefschetzsymplecticstructure}.
The fibres of $\pi$ are finite type complete
convex symplectic manifolds.
We now need to ``complete'' the vertical boundary of $E_1$ so
that we have a fibration over the completion $\widehat{S}$
of $S$.
Let $V := \partial_v E_1 := \pi^{-1}(\partial S)$.
We then attach $A := V \times [0,\infty)$
to $E_1$ by identifying $V \subset E_1$ with $V \times \{0\} \subset A$ to
create a new manifold $\widehat{E}$.
We can extend our map $\pi$ over $A$ by letting
$\pi|_A (v,r) = \pi|_V(v)$ where $v \in V \subset E_1$ and
$r \in [0,\infty)$.
Let $\pi_1 : A \rightarrow V$ 
be the natural projection onto $V$.
Let $\Theta_A := \pi_1^* \Theta|_V$.
We wish to join $\Theta$ with $\Theta_A$ inside $\widehat{E}$.
Let $\zeta : [0,\infty) \rightarrow [0,1]$ be a smooth function such that
$\zeta(x) = 0$ for $x < 0.5$  and $\zeta(x) = 1$ for $x > 1$.
The $1$-form $\Theta$ smoothly extends to a $1$-form $Q$ defined on $V \times [0,\epsilon] \subset A$
such that $dQ$ restricted to the fibres of $\pi$ is a symplectic form and
$Q = \Theta$ in the region $G' = (\partial F \times [1,\infty)) \times \widehat{S}$
(i.e. the region where the Lefschetz fibration is a product near infinity).
Let $\Theta_A'(v,r) := (1-\zeta(Kr/ \epsilon)) Q + \zeta(Kr/ \epsilon)\Theta_A$
where $K$ is a large constant.
We can now smoothly extend $\Theta$ over $A$ by the 1-form
$\Theta_A'$. Note that $d\Theta_A'$ is a symplectic form on all the
fibres for $K$ large enough because $\Theta|_V = \Theta_A|_{V \times \{0\}}$,
and hence $dQ$ restricted to the fibres is as close as we like to $d\Theta_A$
restricted to the fibres. This means that $(1-t)dQ + td\Theta_A$ is
arbitrarily close to a symplectic form when restricted to the fibres and hence it is a symplectic form
when restricted to the fibres for $t \in [0,1]$ and this implies that
$d\Theta_A'$ is a symplectic form on all the
fibres.

\begin{defn} \label{defn:halfconvexlefschetzcompletion}
$(\widehat{E},\pi)$ is called the {\bf completion} of $(E,\pi)$.
\end{defn}

Note that the base of our completed fibration is $(\widehat{S},\theta_S)$.

\begin{defn} \label{defn:convexlefschetzfibration}
Any fibration which is the completion
of a compact convex Lefschetz fibration is called a 
{\bf complete convex Lefschetz fibration}.
\end{defn}

Note that if we add  a large multiple of $\pi^*\theta_S$ to $\Theta$ then
$(\widehat{E},\theta)$ is a complete finite type convex symplectic manifold.
Lefschetz fibrations have well defined parallel transport maps due to the fact
that the fibration is trivial near the horizontal boundary of $E$.
Now we need to deal with almost complex structures on $\widehat{E}$,
as this will be useful when we later define $SH_*^{\nsmbox{lef}}$.
Let $J$ (resp. $j$) be an almost complex structure on $\widehat{E}$
(resp. $\widehat{S}$). 

\begin{defn} \label{defn:lefschetzalmostcomplexstructure}
We say that $(J,j)$ are {\bf compatible} with $(\widehat{E},\pi)$ if:

\begin{enumerate}

\item $\pi$ is
$(J,j)$-holomorphic, and that $J=J_0$ near $E^{\nsmbox{crit}}$ and
$j = j_0$ near $S^{\nsmbox{crit}}$.

\item $j$ is convex at infinity with respect to the convex symplectic
structure of $\widehat{S}$ 
(i.e. $\theta_S \circ j = dr$ for large $r$ where $\theta_S$ is the
contact form at infinity on $\widehat{S}$ and $r$ is the 
radial coordinate of the cylindrical end of $\widehat{S}$).

\item $J$ is a product $(J_F,j)$ on the region
$C \times S$ where $C$ is the cylindrical end $\partial F \times [1,\infty)$
of $\widehat{F}$, and $J_F$ is convex at infinity for $F$.
($\theta|_{\widehat{F}} \circ J_F = dr$ for large $r$ where $r$ is the radial coordinate
of the cylindrical end).


\item $\omega(\cdot,J \cdot)$ is symmetric and positive definite.

\end{enumerate}

\end{defn}

If $(\widehat{E},\pi)$ is a complete convex Lefschetz fibration then
the space of such almost complex structures is nonempty and contractible
(see \cite[section 2.2]{Seidel:longexactsequence}).
We wish to have a slightly larger class of almost complex structures.

\begin{defn} \label{defn:admissiblecomplexstructure}
We define
${\mathcal J}^h(\widehat{E})$ to be the space of almost complex structures
on $\widehat{E}$ such that for each $J$ in this space,
there exists a $(J_1,j_1)$ compatible with $(\widehat{E},\pi)$ and a compact set
$K \subset \widehat{E}$ with $J = J_1$ outside $K$ and
with $\omega(\cdot,J \cdot)$ symmetric and positive definite everywhere.
\end{defn}

The set of such complex structures is still contractible.

\subsection{Symplectic homology} \label{section:symplectichomology}

In this section we will discuss symplectic homology
as defined by Viterbo in \cite{Viterbo:functorsandcomputations}
for finite type Stein manifolds. For simplicity we will assume that
our homology theory has coefficients in $\Z /2$.

Let $(M,\theta)$ be a compact convex symplectic manifold.
Our manifold $\widehat{M}$
has a cylindrical
end symplectomorphic to $(N \times [1,\infty), d(r \alpha))$
where $r$ is a coordinate on $[1,\infty)$ and $\alpha$ is
a contact form on $N$.
We choose a smooth function 
$H : {\mathbb S}^1 \times \widehat{M} \longrightarrow \R$
such that each $H_t$ is of the form $H_t = ar+b$
at infinity where $a$ and $b$ are some constants independent of $t$.
We call such a Hamiltonian {\bf admissible}.
We also choose an ${\mathbb S}^1$ family of
almost complex structures $J_t$ compatible with
the symplectic form. 
We assume that
$J_t$ is convex with respect to this cylindrical end outside some large compact set
(i.e. $\theta \circ J_t = dr$). We call the constant
$a$ the {\bf slope at infinity}. We also say that $J_t$ is {\bf admissible}.

We define an ${\mathbb S}^1$ family
of vector fields $X_{H_t}$ by
$\omega(X_{H_t},\cdot) = dH_t(\cdot)$. The flow $\nsmbox{Flow}^t_{X_{H_t}}$
is called the {\bf Hamiltonian flow}.
A {\bf $1$-periodic orbit} of $H_t$ is a
path of the form $l : [0,1] \rightarrow \widehat{M}$,
$l(t) := \nsmbox{Flow}^t_{X_{H_t}}(p)$ where $p \in \widehat{M}$ and $l(0) = l(1)$.
We choose the cylindrical end and the slope of our
Hamiltonians so that the union of the
$1$-periodic orbits form a compact set.
Let $F := \nsmbox{Flow}^1_{X_{H_t}}$, then we have a correspondence
between 1-periodic orbits $x$ and fixed points $p$ of $F$. In particular
we say that $x$ is non-degenerate if 
$DF|_p : T_p\widehat{M} \rightarrow T_p\widehat{M}$
has no eigenvalue equal to $1$.
We can also assume that the $1$-periodic orbits of
our Hamiltonian flow $X_{H_t}$ are non-degenerate.
We call a Hamiltonian $H$ satisfying these conditions a
{\bf non-degenerate admissible Hamiltonian} or an
admissible Hamiltonian with non-degenerate orbits.

From now on we will assume that $c_1(M) = c_1(\widehat{M})=0$. If we are given a
trivialisation of the canonical bundle ${\mathcal K} \cong {\mathcal O}$,
then for each orbit $x$, we can define an index of $x$ called the
{\bf Robbin-Salamon} index (This is equal to the Conley-Zehnder
index taken with negative sign).
The choice of these indices depend on the choice of trivialisation
of ${\mathcal K}$ up to homotopy but the indices are canonical if
$H_1(M)=0$.
We denote this index by $\nsmbox{ind}(x)$.
Let
\[CF_k(M,H,J) := \displaystyle\bigoplus_{\nsmbox{Flow}^1_{X_{H_t}}(x)=x,ind(x)=k} \Z /2 \langle x \rangle\]
For a 1-periodic orbit $\gamma$ we define the {\bf action} $A_H(\gamma)$:
\[A_H(\gamma) := -\int_0^1 H(t,\gamma(t))dt -\int_\gamma \theta\]
This is the convention of \cite{Viterbo:functorsandcomputations} 
and \cite{Oancea:kunneth}. This differs in sign from
Seidel's convention in \cite{Seidel:biasedview}.
We will now describe the differential \[\partial : CF_k(M,H,J) \rightarrow CF_{k-1}(M,H,J)\]
We consider curves
$u : \R \times {\mathbb S}^1 \longrightarrow \widehat{M}$ satisfying
the perturbed Cauchy-Riemann equations:
\[ \partial_s u + J_t(u(s,t)) \partial_t u = \nabla^{g_t} H\]
where $\nabla^{g_t}$ is the gradient associated to the ${\mathbb S}^1$ family of metrics
$g_t := \omega(\cdot,J_t(\cdot))$.
For two periodic orbits $x_{-},x_{+}$ let
$\bar{U}(x_{-},x_{+})$ denote the set of all curves $u$ satisfying
the Cauchy-Riemann equations such that $u(s,\cdot)$ converges
to $x_{\pm}$ as $s \rightarrow \pm \infty$. This has a natural
$\R$ action given by replacing the coordinate $s$ with $s+v$
for $v \in \R$. Let $U(x_{-},x_{+})$ be equal to $\bar{U}(x_{-},x_{+}) / \R$.
For a $C^{\infty}$ generic admissible Hamiltonian and complex structure
we have that $U(x_{-},x_{+})$ is an $\nsmbox{ind}(x_{-}) - \nsmbox{ind}(x_{+}) -1$  dimensional
manifold (see \cite{FHS:transversalitysymplectic}). 
There is a maximum principle which ensures
that all elements of $U(x_{-},x_{+})$ stay inside
a compact set $K$ (see \cite[Lemma 1.5]{Oancea:survey}). 
Hence we can use a compactness theorem (see for instance
\cite{BEHWZ:compactnessfieldtheory}) which ensures that
if $\nsmbox{ind}(x_{-}) -1 = \nsmbox{ind}(x_{+})$, then
$U(x_{-},x_{+})$ is compact and hence a finite set.
Let $\# U(x_{-},x_{+})$ denote the number of
elements of $U(x_{-},x_{+})$ mod $2$. Then we have a differential:
\[\partial : CF_k(M,H,J) \longrightarrow CF_{k-1}(M,H,J) \]
\[\partial \langle x_{-} \rangle := \displaystyle \sum_{\nsmbox{ind}(x_{+}) = \nsmbox{ind}(x_{-}) -1 } \# U(x_{-},x_{+}) \langle x_{+} \rangle\]
By analysing the structure of 1-dimensional moduli spaces, one shows
 $\partial^2=0$ and defines
$SH_*(M,H,J)$ as the homology of the above chain complex.
As a $\Z / 2\Z$ module $CF_k(M,H,J)$ is independent of $J$,
but its boundary operator does depend on $J$. The homology
group $SH_*(M,H,J)$ depends on $M,H$ but is independent of $J$
up to canonical isomorphism.
Note that for each $f \in \R$ we have a subcomplex 
generated by orbits of action $\leq f$. 
The reason why this is a subcomplex is because if $\bar{U}(x_{-},x_{+})$
is a non-empty set, then $A_H(x_-) \geq A_H(x_+)$. This means that
the differential decreases action and hence if $x_-$ is an orbit of
action $\leq f$, then $\partial(x_-)$ is a linear combination of orbits
of action $\leq f$.
The homology
of such a complex is denoted by:
$SH_*^{\leq f}(M,H,J)$.

If we have two non-degenerate admissible
Hamiltonians $H_1 \leq H_2$ and two admissible
almost complex structures $J_1,J_2$, then there is a natural
map:
\[SH_*(M,H_1,J_1) \longrightarrow SH_*(M,H_2,J_2)\]
This map is called a {\bf continuation map}.
This map is defined from a map $C$ on the chain level as follows:
\[C : CF_k(M,H_1,J_1) \longrightarrow CF_{k}(M,H_2,J_2) \]\[\partial \langle x_{-} \rangle := \displaystyle \sum_{\nsmbox{ind}(x_{+}) = \nsmbox{ind}(x_{-}) } \# P(x_{-},x_{+}) \langle x_{+} \rangle\]
where $P(x_{-},x_{+})$ is the set of solutions of the following equations:
Let $K_s$ be a smooth increasing family of admissible Hamiltonians joining $H_1$ and $H_2$
and $J_s$ a smooth family of admissible almost complex structures joining $J_1$ and $J_2$.
Then for a $C^{\infty}$ generic family $(K_s,J_s)$, the set  $P(x_{-},x_{+})$ is the
set of solutions to the parameterized Floer equations
\[ \partial_s u + J_{s,t}(u(s,t)) \partial_t u = \nabla^{g_t} K_{s,t}\]
such that $u(s,\cdot)$ converges
to $x_{\pm}$ as $s \rightarrow \pm \infty$.
If we have another such increasing family admissible Hamiltonians joining $H_1$ and $H_2$
and another smooth family of admissible almost complex structures joining $J_1$ and $J_2$,
then the continuation map induced by this second family is the same as
the map induced by $(K_s,J_s)$. 
The composition of two continuation maps is a continuation map.

If we take the direct limit of all these maps with respect
to admissible Hamiltonians ordered by $\leq$, then
we get our symplectic homology groups $SH_*(M)$.
Supposing we have a family of Hamiltonians $(H_\lambda)_{\lambda \in \Lambda}$
ordered by $\leq$. We say that a family of Hamiltonians $(H_i)_{i \in I \subset \Lambda}$
is {\it cofinal} if for every $\lambda \in \Lambda$, there exists an $i \in I$
and a constant $a$ such that $H_\lambda \leq H_i + a$.
The fact that continuation maps are natural ensures that we
can also define $SH_*(M)$ as the direct limit of all these maps with respect to any cofinal family of Hamiltonians.
If we have a degenerate admissible Hamiltonian $H$, then
we can still define symplectic Homology $SH_*(H,J)$
as a direct limit $SH_*(H_k,J_k)$ where $H_k$ and $J_k$
are $C^{\infty}$ generic making $SH_*(H_k,J_k)$ well defined, and
such that $H_k$ and $J_k$ tend to $H$ and $J$.
The limit is taken with respect to continuation maps.
This is independent of choices of $H_k$ and $J_k$
due to the naturality of continuation maps
(see \cite[Remark 1.2]{Viterbo:functorsandcomputations}).

The symplectic homology groups also have a ring structure.
The product is called the pair of pants product and it makes $SH_{n+*}(M)$
into a $\Z / 2\Z$ graded algebra where $n$ is the complex dimension of $M$.
The product is defined in the following way:
Take an admissible Hamiltonian $H_1$ and another admissible Hamiltonian
$H_2$ such that $H_1 = H_2$ outside some large compact set, and such
that $H_1(1,\cdot) = H_2(0,\cdot)$ with all time derivatives.
Then $H_1 \# H_2$ is also admissible where $H_1 \# H_2$ is defined as:
\[ H_1 \# H_2 (t,x) := \left\{
              \begin{array}{ll}
                   2H_1(2t,x) & (t \in [0,\frac{1}{2}])\\
                   2H_2(2t - 1,x) & (t \in [\frac{1}{2},1])
              \end{array}
       \right. 
\]
 We define a chain map (maybe after perturbing the Hamiltonians $H_1,H_2,H_1 \# H_2$ slightly):
\[CF_k(H_1,J) \otimes CF_j(H_2,J) \rightarrow CF_{k + j - n}(H_1 \# H_2,J),\]
\[x_1 \otimes x_2 \rightarrow 
\sum_{\nsmbox{ind}(x_3) =  k+j-n} 
\# {\mathcal M}(H,J,x_1,x_2,x_3) \langle x_3 \rangle.\]
The set ${\mathcal M}(H,J,x_1,x_2,x_3)$ is the set of maps
$u : \P^1 \setminus \{0,1,\infty\} \rightarrow M$
which satisfy some Floer type equations and such that
each puncture converges to an orbit. The details
of this are explained in \cite{AbbondandoloSchwarz:noteonfloerloop},
\cite{abbondandoloschwarz:cotangentloopproduct}
and \cite[Section 8]{Seidel:biasedview}. 
This commutes with continuation maps and hence we can
take the direct limit of these maps with respect to the
ordering $\leq$ giving us a map:
\[SH_{n+i}(M) \otimes SH_{n+j}(M) \rightarrow SH_{n+i+j}(M).\]
This also has a unit in degree $n$ and is given by counting holomorphic planes
(\cite[Section 8]{Seidel:biasedview}).

Suppose that $(M',\theta')$ is a compact convex symplectic manifold which
is an exact submanifold of $M$,
then there exists a natural map
\[i : SH_*(M) \longrightarrow SH_*(M')\]
called the transfer map. 
The composition of two of these
transfer maps is another transfer map. 
These maps are introduced in \cite[Section 2]{Viterbo:functorsandcomputations}
and studied in \cite[Section 3.3]{Cieliebak:handleattach}.
Let $N$ be a general convex symplectic manifold.
Let $(W_j)_{j \in J}$ be the set of codimension $0$ compact convex
exact symplectic submanifolds. This is a directed system, where
the morphisms are just inclusion maps. If $W_{j_1} \subset W_{j_2}$,
then we have a transfer map $SH_*(W_{j_2}) \rightarrow SH_*(W_{j_1})$.
Because the transfer map points in the opposite direction
(i.e. from $W_{j_2}$ to $W_{j_1}$ instead of $W_{j_1}$ to $W_{j_2}$), we have
an inverse system $(SH_*(W_j))_{j \in J}$ where the morphisms
now are transfer maps.
Hence we can define the ring $SH_*(N)$ as the inverse limit of
this inverse system. Because the definition only involves
codimension $0$ exact symplectic manifolds, we have that
it is invariant under exact symplectomorphism.
We have that $SH_*(N)$ is invariant under convex symplectic deformations
(see \cite[Theorem 2.12]{McLean:thesis}).
We also have that for a compact convex symplectic manifold $M$,
$SH_*(M) = SH_*(\widehat{M})$ (see \cite[Theorem 2.12]{McLean:thesis}).

Symplectic homology $SH_*(M)$ also has a natural $H_1(M)$ grading. 
Each orbit $x$ of some Hamiltonian $H$ is contained in some
$H_1(M)$ class, and the definition of the differential ensures
that $\partial x$ is a linear combination of orbits in the
same $H_1(M)$ class. Also if we have two orbits $x_1$ and $x_2$ of some
Hamiltonian $H_1$ and $H_2$ respectively,
then they are in $H_1(M)$ classes $a$ and $b$. If we have
a pair of pants, where two of the boundaries are in classes
$a$ and $b$, then the third boundary is in the class $a+b$
(if we orient the boundaries correctly). This ensures that
the pants product of $x_1$ and $x_2$ is a linear combination
of orbits of $H_1 \# H_2$ in the homology class $a+b$. Hence the pair of
pants product for $SH_*(M)$ is additive in $H_1(M)$.

\begin{theorem} \label{thm:handleattaching}
\label{theorem:endconnsymhom}
Let $M,M'$ be finite type Stein manifolds of real dimension greater than $2$, then \\
$SH_*(M \#_e M') \cong SH_*(M) \times SH_*(M')$ as rings. 
Also the transfer map $SH_*(M \#_e M') \rightarrow SH_*(M)$
is just the natural projection
\[SH_*(M) \times SH_*(M') \twoheadrightarrow SH_*(M).\] 
\end{theorem}
Cieliebak in  \cite{Cieliebak:handleattach} showed that the above theorem
is true if we view $SH_*$ as a vector space.  
We will prove that we have a ring
isomorphism in section \ref{subsection:handleattaching}.

\subsection{Symplectic homology and Lefschetz fibrations} \label{subsection:symhomlefschetz}

We need three theorems which relate symplectic homology to
Lefschetz fibrations. These are the key ingredients in
proving that our exotic Stein manifolds are pairwise distinct.
The proofs of these theorems will be deferred to sections
\ref{section:lefschetzcofinal} and \ref{section:transferisomorphism}.
Theorem \ref{theorem:lefschetzcofinalfamily} is very close
to Oancea's K\"unneth formula \cite{Oancea:kunneth} but
theorems \ref{theorem:lefschetzcofinalsequence} and
\ref{theorem:lefschetztransferisomorphism} are new and the
main part of the story.
Throughout this section we will
let $\pi' : E' \rightarrow S'$ be a compact convex Lefschetz fibration
with fibre $F'$.
From now on we will assume that $c_1(E')=0$ and to make $SH_*(E')$
graded we will choose a trivialisation of the canonical bundle of $E'$.
Note that when we talk about symplectic homology of a
compact convex Lefschetz fibration, we mean the symplectic
homology of its completion with respect to the
convex symplectic structure.
The fibration $\widehat{E'}$ can be partitioned into three sets as follows:
\begin{enumerate}
\item $E' \subset \widehat{E'}$
\item $A := F'_e \times \widehat{S'}$, where 
$F'_e := \partial F' \times \R_{\geq 1}$ is the cylindrical
end of $\widehat{F'}$.
\item $B := \widehat{E'} \setminus (A \cup E')$
\end{enumerate}
\begin{figure}[b]
\end{figure}
The set $B$ is of the form 
$(A_1 \times \R_{\geq 1}) \bigsqcup (A_2 \times \R_{\geq 1}) \cdots \bigsqcup (A_n \times \R_{\geq 1})$,
where $A_i$ is a mapping torus of the monodromy symplectomorphism
around one of the boundary components of $S'$.
Here is a picture of the regions $E'$, $A$ and $B$. 
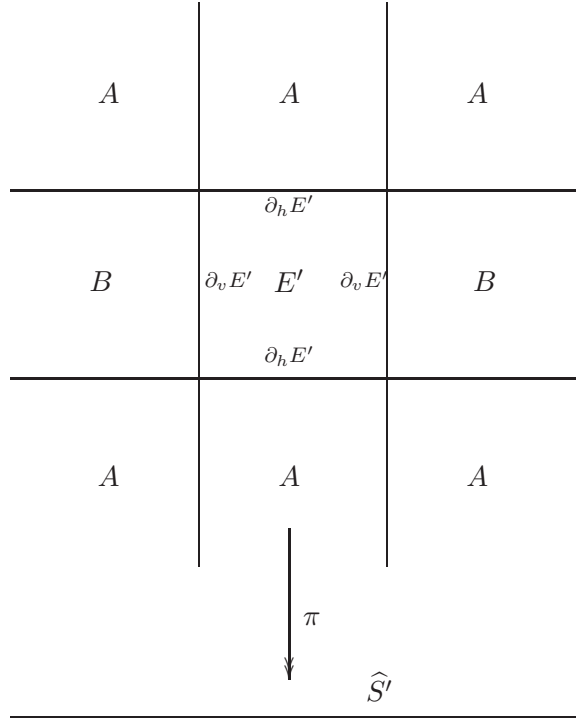
\begin{figure}[b]
\end{figure}
\begin{figure}[b]
\end{figure}
 \[
 \xy
 (0,-25)*{}="A"; (75,-25)*{}="B";
 (0,-50)*{}="C"; (75,-50)*{}="D";
 (25,0)*{} ="E"; (25,-75)*{} ="F";
 (50,0)*{} ="G"; (50,-75)*{} ="H";
 (37,-37)*{}= "I";
 (37,-70)*{}="J"; (37,-90)*{}="K";
 (0,-95)*{}="L"; (75,-95)*{}="M";
 "A"; "B" **\dir{-};
 "C"; "D" **\dir{-};
 "L"; "M" **\dir{-};
 "E"; "F" **\dir{-};
 "G"; "H" **\dir{-};
{\ar@{->>} "J";"K"};
 "I" *{E'};
 "I" + (10,0) *{\scriptstyle \partial_v E'};
 "I" + (-8,0) *{\scriptstyle \partial_v E'};
 "I" + (0,10) *{\scriptstyle \partial_h E'};
 "I" + (0,-10) *{\scriptstyle \partial_h E'};
 "E" + (-12,-12) *{A};
 "E" + (12,-12) *{A};
 "G" + (12,-12) *{A};
 "F" + (-12,12) *{A};
 "F" + (12,12) *{A};
 "H" + (12,12) *{A};
 "A" + (12,-12) *{B};
 "B" + (-12,-12) *{B};
 "K" + (3,8) *{\pi};
 "K" + (12,-1) *{\widehat{S'}};
 \endxy 
 \]

\bigskip

Let $\pi_1 : A \twoheadrightarrow F'_e$ be the natural
projection onto $F'_e$.

\begin{defn} \label{defn:lefschetzadmissiblehamiltonian}
Let $H_{S'}$ be an admissible Hamiltonian for the base $\widehat{S'}$.
Let $H_{F'}$ be an admissible Hamiltonian for the fibre $\widehat{F'}$.
We assume that $H_{F'} = 0$ on $F'\subset \widehat{F'}$. 
The map $H : \widehat{E'} \rightarrow \R$ is called a
{\bf Lefschetz admissible Hamiltonian} if 
$H|_A = \pi^* H_{S'} + \pi_1^* H_{F'}$ and
$H|_B = \pi^* H_{S'}$ outside some large compact set.
We say that $H$ has {\bf slope} $(a,b)$ if
$H_{S'}$ has slope $a$ at infinity and $H_F$ has slope
$b$ at infinity.
\end{defn}
Let $H$ be a Lefschetz admissible Hamiltonian and let $J$ 
be an admissible almost complex structure for $E'$. We will
call the pair $(H,J)$ a Lefschetz admissible pair.
For generic
$(H,J)$ we can define
$SH_*(E',H,J)$ (see section \ref{section:lefschetzcofinal} for
more details). If $(H_1,J_1)$ is another generic
Lefschetz admissible pair such that $H \leq H_1$,
then there is a continuation map
$SH_*(H,J) \rightarrow SH_*(H_1,J_1)$ induced by an
increasing homotopy from $H$ to $H_1$ through Lefschetz admissible Hamiltonians.
Hence, we have a direct limit
$SH_*^l(E') := \underset{(H,J)}{\varinjlim} \nsmbox{ } SH_*(H,J)$
with respect to the ordering $\leq$ on Hamiltonians $H$.
This has the natural structure of a ring with respect
to the pair of pants product.
\begin{theorem} \label{theorem:lefschetzcofinalfamily}
There is a ring isomorphism $SH_*(E') \cong SH_*^l(E')$.
\end{theorem}
This will be proved in section \ref{section:lefschetzcofinal}.
Let $\epsilon$ be smaller than the length of the shortest Reeb orbit
of $\partial F'$.
A Hamiltonian $H$ is called a {\bf half admissible Hamiltonian}
if it is Lefschetz admissible and has slope $(a,\epsilon)$.
We let $J$
be an admissible almost complex structure for $E'$.
\begin{defn} \label{lefschetzsymplectichomology}
We define
\[SH^{\nsmbox{lef}}_*(E') := \underset{(H,J)}{\varinjlim} \nsmbox{ } SH_*(H,J)\]
as the direct limit with respect to the ordering $\leq$ on
half admissible Hamiltonians $H$. This has the structure of a ring as usual.
\end{defn}
The difference between $SH_*(E')$ and $SH_*^l(E')$ is that
$SH_*(E')$ is defined using Hamiltonians which are linear
with respect to some fixed cylindrical end.
The ring $SH_*^l(E')$ is defined using Hamiltonians which are linear
in the horizontal and vertical directions with respect to some Lefschetz fibration.
The difference between $SH_*^{\nsmbox{lef}}(E')$ and the
other homology theories is that the slopes of a
cofinal family of half admissible Hamiltonians
do not have to tend to infinity pointwise in the
vertical direction. This has to be true for $SH_*(E')$ and
$SH_*^l(E')$ where the Hamiltonians have to get steeper and steeper
at infinity in all directions.
Because a half admissible Hamiltonian is Lefschetz admissible,
we have a natural ring homomorphism:
\[\Phi : SH_*^{\nsmbox{lef}}(E') \rightarrow SH_*^l(E').\]
This comes from continuation maps $SH_*(H,J) \rightarrow SH_*(K,J)$
where $H$ is a half admissible Hamiltonian and $K$ is a Lefschetz admissible Hamiltonian.
This is because $H \leq K$ when $K$ is large enough.

\begin{theorem} \label{theorem:lefschetzcofinalsequence}
If $S'=\D$, the unit disk, then $\Phi$ is an isomorphism of rings.
Hence by Theorem \ref{theorem:lefschetzcofinalfamily},
\[SH_*(E') \cong SH^{\nsmbox{lef}}_*(E')\]
as rings.
\end{theorem}
This will be proved in section \ref{section:bettercofinal}.
Let  $F'$ (resp. $F''$) be a smooth fibre of
$E'$ (resp. $E''$). Let $F'$ and $F''$ be Stein domains with
$F''$ a holomorphic and symplectic submanifold of $F'$.
\begin{theorem} \label{theorem:lefschetztransferisomorphism}
Suppose $E'$ and $E''$ satisfy the following properties:
 \begin{enumerate}
   \item \label{item:subfibration1}
$E''$ is a subfibration of $E'$ over the same base.
   \item \label{item:monodromysupport1}
The support of all the monodromy maps of $E'$ are
contained in the interior of $E''$.
   \item \label{item:complexcurvecriterion1}
Any holomorphic curve in $F'$
with boundary inside $F''$ must be contained in $F''$.
 \end{enumerate}
Then $SH_*^{\nsmbox{lef}}(E') \cong SH_*^{\nsmbox{lef}}(E'')$ as rings.
\end{theorem}
This theorem will be proved in section \ref{section:transferisomorphism}.
Combining this theorem with Theorem \ref{theorem:lefschetzcofinalsequence}
proves the key theorem \ref{theorem:fibrationaddition} in the
introduction of this paper.

\subsection{The Kaliman modification} \label{section:kalimanmodification}

In order to produce examples of exotic symplectic manifolds,
we first need to construct exotic algebraic varieties. One
tool used for constructing these manifolds is called the
Kaliman modification.
Our treatment follows section 4 of \cite{Zaidenberg:1998exot}.

Consider a triple $\left( M, D, C \right)$ where $C \subseteq  D \subseteq M$
are complex varieties. Let $M$ and $C$ be smooth, $D$ be an
irreducible hypersurface in $M$, and $C$ be a closed
subvariety contained in the smooth part of $D$ such that
$\nsmbox{dim}(C) < \nsmbox{dim}(D)$.
\begin{defn} \label{defn:kalimanMod}
(see \cite{Kaliman:eisenman})
The {\bf Kaliman modification} $M^\prime$ 
of $\left( M, D, C \right)$ is defined by
$M^\prime := \nsmbox{Kalmod}\left( M, D, C \right) = \tilde{M} 
\setminus \tilde{D}$ where $\tilde{M}$ is the
blowup of $M$ along $C$ and $\tilde{D}$ is the proper
transform of $D$ in $\tilde{M}$.
\end{defn}

The Kaliman modification of an affine variety is again
an affine variety (see \cite{Kaliman:eisenman}).


\begin{lemma} \label{lm:kaliman}
\cite[Theorem 3.5]{Kaliman:eisenman}
Suppose that (i) $D$ is a topological manifold, and (ii)
$D$ and $C$ are acyclic. Then $M^\prime$ is contractible
iff $M$ is.
\end{lemma}


\begin{example} \label{ex:tDP}
({\bf tom Dieck-Petrie surfaces} see
\cite{tDP:1990conaff, tDP:1989moh})
For $k>l\geq 2$ with $(k,l)$ coprime, the triple
$A_{k,l}:=\left( \C^2, \{ x^k-y^l=0 \}, \{\left(1,1\right)\} \right)$
satisfies the conditions of Lemma \ref{lm:kaliman}.
Hence $X_{k,l} = \nsmbox{Kalmod}\left( A_{k,l} \right)$ is contractible. 
Note: $X_{k,l}$ is isomorphic to
\[ \left\{ \frac{\left(xz+1\right)^k-\left(yz+1\right)^l-z}{z}=0 \right\}.\]
Here $x,y,z$ are the standard coordinates of $\C^3$. Also the numerator
of this fraction is divisible by $z$, hence the above fraction
is a polynomial.
\end{example}

Here is another construction:


\begin{example} \label{ex:repeatedKal}
({\bf Kaliman} \cite{Kaliman:eisenman})
If we have a contractible affine variety $M$ of complex dimension
$n$, then we can construct a contractible affine variety
\[ M_k := \nsmbox{Kalmod}\left( M \times \C, M \times \{ p_1,\dots,p_k \}, \{ \left(a_1,p_1\right),\dots,\left(a_k,p_k\right) \} \right) \] 
where $p_i$ are distinct points in $\C$ and $a_i$ are points in $M$.
This variety is contractible by a repeated application of Lemma \ref{lm:kaliman}, 
because it is a repeated Kaliman modification with $D$ isomorphic to $M$
and $C$ a point. 
There are obvious variants: replace $\C$ and
$\{ p_1,\dots,p_k \}$ with some contractible variety
and a disjoint union of contractible irreducible hypersurfaces, etc.
\end{example}

At the moment we are only discussing contractibility of varieties.
We need to produce varieties diffeomorphic to some $\C^n$.
We will use the h-cobordism theorem to achieve this stronger condition.

\begin{corollary} \label{cor:contractibleexotic}
(See \cite[Page 174]{DimcaChoudary:complexhypersurfaces}, \cite{Ramanujam:affineplane} and
\cite[Proposition 3.2]{Zaidenberg:1998exot})
Let $M$ be a contractible Stein manifold of finite type.
If $n := \nsmbox{dim}_{\C}M \geq 3$ then $M$ is diffeomorphic to $\C^n$.

\proof

Let $(J,\phi)$ be the Stein structure associated with $M$.
We can also assume that $\phi$ is a Morse function.
For $R$ large enough, the domain $M_R := \{\phi < R\}$
is diffeomorphic to the whole of $M$ as $M$ is of finite type.
We want to show that the boundary of $\bar{M}_R := \{\phi \leq R\}$
is simply connected, then the result follows from
the h-cobordism theorem.

The function $\psi := R - \phi$ only has critical
points of index $\geq n \geq 3$ because the function
$\phi$ only has critical points of index $\leq n$
(see \cite[Corollary 2.9]{Eliashberg:symplecticgeometryofplushfns}).
Viewing $\psi$ as a Morse function,
$\bar{M}_R$ is obtained from
$\partial \bar{M}_R$
by attaching handles of index $\geq 3$.
This does not change
$\pi_1$, hence $\partial \bar{M}_R$ is simply connected
because $\bar{M}_R$ is simply connected.
\qed
\end{corollary}

We now need a theorem which relates the Kaliman modification with
symplectic homology. We do this via Lefschetz fibrations.
Let $X$,$D$,$M$ be as in Example \ref{example:algebraicsteinmanifold}.
Let $Z$ be an irreducible
divisor in $X$ and $q \in (Z \cap M)$ a point in
the smooth part of $Z$. We assume there is a rational function $m$
on $X$ which is holomorphic on $M$ such that
$\overline{m^{-1}(0)}$ is reduced and irreducible and
$Z = \overline{m^{-1}(0)}$.
Let $M' := \nsmbox{Kalmod}(M,(Z \cap M),\{q\})$,
and let $M'' := M \setminus Z$.
Suppose also that $\nsmbox{dim}_{\C}X \geq 3$.
We also assume
that $c_1(M') = c_1(M'') = 0$.
\begin{theorem} \label{theorem:KalimanTransfer}
$SH_*(M'') = SH_*(M')$.
\end{theorem}

This theorem follows easily from the key theorem
\ref{theorem:fibrationaddition} and the following theorem:
\begin{theorem} \label{theorem:KalimanLefschetz}
There exist compact convex Lefschetz fibrations $E'' \subset E'$
respectively satisfying the conditions of Theorem
\ref{theorem:fibrationaddition} such that
$E'$ (resp. $E''$) is convex deformation equivalent to $M'$ (resp. $M''$).
\end{theorem}
This will be proved in the appendix (\ref{section:appendix}). 
The basic idea of the proof
is to use Lefschetz fibrations defined in an algebraic way.

\section{Proof of the main theorem}

\subsection{Construction of our exotic Stein manifolds} \label{section:exoticconstruction}

First of all, we will construct a Stein manifold $K_4$
diffeomorphic to $\C^4$. We will then construct
Stein manifolds $K_n$ diffeomorphic to $\C^n$
for all $n>3$ from $K_4$. Finally using end connect
sums we will construct infinitely many Stein manifolds $(K_n^k)_{k \in \N}$
diffeomorphic to $\C^n$ for all $n>3$.

We define the polynomial 
$P(z_0,\dots,z_3) := z_0^7 + z_1^2 + z_2^2 + z_3^2$
and $V := \{P = 0\} \subset \C^4$.
Let ${\mathbb S}^7$ be the unit sphere in $\C^4$.
\begin{theorem} \label{theorem:brieskorntopology}
{\bf \cite{Brieskorn:sphere}} $V \cap {\mathbb S}^7$ is homeomorphic
to ${\mathbb S}^5$.
\end{theorem}
Since $V$ is topologically the cone on the link $V \cap {\mathbb S}^7$,
\begin{corollary} \label{corollary:brieskorntopology}
$V$ is homeomorphic to $\R^6$.
\end{corollary}
Let $p \in V \setminus \{0\}$. We let
$K_4 := \nsmbox{Kalmod}(\C^4,V,\{p\})$.
Now by Corollary \ref{corollary:brieskorntopology} and
Lemma \ref{lm:kaliman} we have that $K_4$ is
contractible. Hence by Theorem \ref{cor:contractibleexotic}
we have that $K_4$ is diffeomorphic to $\C^4$.
We will now construct the varieties $K_n$ by 
induction. Suppose we have constructed the varieties
$K_4,\dots,K_n$, we wish to construct the variety
$K_{n+1}$. We do this using example \ref{ex:repeatedKal}.
This means that we will define 
$K_{n+1} := \nsmbox{Kalmod}(K_n \times \C, K_n \times \{0\}, (q,0))$
where $q$ is a point in $K_n$.
All these are affine varieties and hence have Stein structures by
Example \ref{example:algebraicsteinmanifold}.
Finally, we define
\[K_n^k := \underset{i=1 \dots k}{\#_e} K_n\]
which is the $k$ fold end connect sum of $K_n$.
The aim of this paper is to show that if 
$K_n^k \sim K_n^m$ then $k=m$.

\subsection{Proof of the main theorem in dimension 8} \label{section:dim8proof}
Here we will prove Theorem \ref{thm:infinitelymanysteinmanifolds} in dimension $8$.
In this section we wish to show that if $K_4^k \sim K_4^m$ then $k=m$.
Let $M' := K_4^1$.
By Theorem \ref{thm:handleattaching}, $SH_*(K_4^k) = \prod_{i=0}^k SH_*(M')$.
Hence if $i(M')$ is finite,
$i(K_4^k) = i(M')^k$ where $i(M)$ denotes the number of idempotents
of $SH_*(M)$ for any Stein manifold $M$. So in order to
distinguish these manifolds, we need to show that $1<i(M')<\infty$.
Let $M'' := \C^4 \setminus V$ where $V$ is defined in 
section \ref{section:exoticconstruction}.
By Theorem \ref{theorem:KalimanTransfer}, we have that
$SH_*(M'') = SH_*(M')$. We have that $1<i(M'')<\infty$ by
Theorem \ref{theorem:finitelymanyidempotentsofMprimeprime},
hence $1<i(M')<\infty$.

\subsection{Proof of the theorem in dimensions greater than 8} \label{section:dimhigherproof}
Here we will prove Theorem \ref{thm:infinitelymanysteinmanifolds} in dimensions greater than $8$.
Let $K_n := K_n^1$.
For each $n > 4$ we need to show that $1<i(K_n)<\infty$ in order
to distinguish $K_n^k$. This is done by induction.
Suppose that $1 < i(K_n) < \infty$ for some $n$, then we wish
to show that $1 < i(K_{n+1}) < \infty$. We have by Theorem
\ref{theorem:KalimanTransfer}, that 
$SH_*(K_{n+1}) \cong SH_*(K_n \times \C^*)$.
Let $B := K_n \times \C^*$. Let $SH_*^{\nsmbox{\tiny contr}}(\C^*)$ be the subring
of $SH_*(\C^*)$ with $H_1$ grading $0$.

One can check that $SH_*^{\nsmbox{\tiny contr}}(\C^*)$ is a subring
isomorphic to $H^{1-*}(\C^*)$.
In particular $SH_1^{\nsmbox{\tiny contr}}(\C^*) \cong \Z / 2$.
By the K\"unneth formula (see \cite{Oancea:kunneth}),
we have that 
\[SH_{(n+1) + *}(B) \cong SH_{n+*}(K_n) \otimes SH_{1+*}(\C^*).\]
This ring is naturally graded by $H_1(\C^*)$. Hence
any idempotent must be an element of
\[SH_{n+*}(K_n) \otimes SH_{1+*}^{\nsmbox{\tiny contr}}(\C^*) \subset  
SH_{n+*}(K_n) \otimes SH_{1+*}(\C^*)\]
by Lemma \ref{lemma:finitelymanyidempotents}.
The ring $SH_{1+*}^{\nsmbox{\tiny contr}}(\C^*)$ is naturally graded by the
Conley-Zehnder index taken with negative sign because 
$c_1(\C^*)=0$. This means that
any idempotents must live in:
\[SH_{n+*}(K_n) \otimes SH_1^{\nsmbox{\tiny contr}}(C^*) \cong SH_{n+*}(K_n) \otimes \Z / 2 \cong SH_{n+*}(K_n).\]
Hence $i(K_{n+1}) = i(K_n)$.
This means that by induction we have $1 < i(K_n) < \infty$ for all $n >3$
as we proved $1 < i(K_3) < \infty$ in section \ref{section:dim8proof}.
This proves our theorem.

\section{Lefschetz fibration proofs} \label{section:lefschetzfibrationproofs}

Here is the statement and proof of Theorem \ref{theorem:halfconvexstructure}:
{\it Let $(E,\pi)$ be a compact convex
Lefschetz fibration. There exists a constant $K > 0$ such that
for all $k \geq K$ we have:
$\omega := \Omega + k\pi^{*}(\omega_S)$ is a symplectic form,
and the $\omega$-dual $\lambda$
of $\theta := \Theta + k\pi^{*}\theta_S$ is transverse
to $\partial E$ and pointing outwards.}

\proof
We let $K$ be a large constant so that
$\omega := \Omega + \pi^{*}(K\omega_S)$
is a symplectic form
(see \cite[Lemma 1.5]{Seidel:longexactsequence}).
Let $\theta_S' := K\theta_S$ and
$\omega_S' = d\theta_S'$ and $\lambda_S'$
be the $\omega_S'$-dual of $\theta_S$.
Let $U \times V$ be some trivialisation
of $\pi$ around some point $p \in \pi^{-1}(\partial S)$
where $U \subset F$ and $V \subset S$. We let $V$ be
some small half disk around $\pi(p)$ and $U$ is
some small open ball.
Let $\pi^1 : U \times V \twoheadrightarrow U$ be
the natural projection.
Let $\lambda_F$ be the $\Omega|_F$-dual of $\Theta|_F$,
and $\lambda_Q$ be the horizontal lift of $\lambda_S'$.
The $\omega$-dual of $\Theta$ is equal to:
\[\lambda_F + W\]
where $W$ is $\omega$-orthogonal to the vertical plane field tangent to the fibres and
is equal to $0$ near the horizontal boundary of $E$.
The $\omega$-dual of $K\pi^*\theta_S'$ is equal to:
\[G\lambda_Q\]
where $G$ is some function on $U \times V$.
This means that the $\omega$-dual of $\theta$
is:
\[\lambda = \lambda_F + W + G\lambda_Q.\]
Because $W=0$ near the horizontal boundary and because
the horizontal subspaces are tangent to the horizontal
boundary, we have that $\lambda$ is transverse to the
horizontal boundary.
In order to show that $\lambda$ is transverse to the
vertical boundary we need to ensure that we can
make $G$ very large compared to $\lambda_F + W$.
This can be done by making $K$ sufficiently large.

\qed

\section{A cofinal family compatible with a Lefschetz fibration} \label{section:lefschetzcofinal}

In this section we construct a family of Hamiltonians 
$H_k : \widehat{E} \rightarrow \R$ which
behave well with respect to the Lefschetz fibration, so that
\[SH_*^l(E) := \underset{k}{\varinjlim} \nsmbox{ } SH_*(E,H_k,J) = SH_*(E).\]
This would be obvious if $H_k$ belonged to the ``usual'' class
(i.e. linear of slope $k$ on the cylindrical end) but it
is not obvious that our Hamiltonians are
linear with respect to some cylindrical end. Throughout this section,
$(E,\pi)$ is a compact convex Lefschetz fibration.
We let $\Theta,\Omega,\theta,\omega$ be defined as in
section \ref{section:lefschetzfibrations}.

\begin{theorem} \label{thm:admissiblecomplexstr}
Let $H : \widehat{E} \rightarrow \R$
be Lefschetz admissible for $E$ with non-degenerate orbits. 
Then the space of regular almost complex structures \\
${\mathcal J_{\nsmbox{reg}}}(\widehat{E},H)$ is of second category in the space 
${\mathcal J}^h(\widehat{E})$ of admissible almost complex structures with
respect to the $C^\infty$ topology.
\end{theorem}
With a regular almost complex structure, it is possible
to define symplectic homology $SH_*(\widehat{E},H,J)$ with the pair of pants product
$SH_k(\widehat{E},H,J)\otimes SH_l(\widehat{E},H,J) \rightarrow SH_{l+k-n}(\widehat{E},2H,J)$.
Regular almost complex structures are almost complex structures
such that some natural section of a Banach bundle associated with the
almost complex structure is transverse to the zero section. 
This section is described by some linearized
version of the perturbed Cauchy-Riemann equations.
This theorem comes from using results in \cite{FHS:transversalitysymplectic}.
Viterbo in \cite[section 1.1]{Viterbo:functorsandcomputations}
shows us why the Theorem is true in the context of open symplectic
manifolds where the almost complex structure is fixed outside a large compact set
and where there is a maximum principle as in 
\cite[Lemma 1.8]{Viterbo:functorsandcomputations} or \cite[Lemma 1.5]{Oancea:survey}. This argument
applies if we use the maximum principle in \ref{lemma:lefschetzmaximumprinciple}
below.
This ensures that the moduli spaces of Floer trajectories
are manifolds.
For non-generic $(H,J)$, $SH_*(H,J)$ is defined via small
perturbations, and is independent of choice of small
perturbation via continuation map techniques.
We also need a maximum principle to ensure that the
Floer moduli spaces have compactifications.

Let $W$
be a connected component of $\partial S$ where $S$ is the base.
Now $\widehat{S}$ has a cylindrical end $W \times [0,\infty)$.
Let $r_S$ be the coordinate for $[1,\infty)$. 
Let $u : \D \rightarrow \widehat{E}$ satisfy
Floer's equations for some $J \in {\mathcal J}^h{\widehat{E}}$
and some admissible Hamiltonian $H$. Here $\D$ is the unit disk
parameterized by coordinates $(s,t)$.
We can write $H = \pi^*H_S + \pi_1^*H_F$ as in
Definition \ref{defn:lefschetzadmissiblehamiltonian}.
We assume that $H_F=0$ on $F$.
\begin{lemma} \label{lemma:lefschetzmaximumprinciple}
The function $f := r_S \circ \pi \circ u$ cannot have an interior maximum for $r_S$ large. 
\end{lemma}
\proof
Let $f$ have an interior maximum at $q \in \D$.
Let $U$ be a small neighbourhood of $u(q)$.
The symplectic form $\omega$ on $\widehat{E}$
splits the tangent space of $E$ into vertical planes
and horizontal planes. 
Let $V$ be the vertical plane field, and
let $P$ be the horizontal plane field
(the $\omega$-orthogonal of vertical tangent spaces of $\pi$).
Let $\omega_S$ be the symplectic form on the base $S$,
then $\omega_P := \pi^* \omega_S|_P$ is non-degenerate.
This means that there exists a function 
\[g : \pi^{-1}(W \times [0,\infty)) \rightarrow (0,\infty)\]
such that $g\omega_P = \omega|_P$.
We may assume that $J(P) \subset P$ because $J$ is
compatible with $\widehat{E}$ if $r_S$ is large.
Let $p$ be the natural projection $TE \rightarrow P$ induced
by the splitting $TE = V \bigoplus P$.

Floer's equation for $u$ splits up into a horizontal part
associated to $P$ and a vertical part associated to $V$.
The horizontal part can be expressed as:
\[p(\frac{\partial u}{\partial s}) + J p(\frac{\partial u}{\partial t}) = -J\frac{1}{g}G\]
where $G$ is a vector field on $P$ which is the
$\omega_P$-orthogonal to $d\pi^*H_S|_P$ in $P$.
Hence $u' := \pi \circ u$ satisfies the equation:
\[\frac{\partial u'}{\partial s} + j \frac{\partial u'}{\partial t} = -j\frac{1}{u^*(g)}X_{H_S}\]
where $j$ is the complex structure of $S$, and $X_{H_S}$ is the
Hamiltonian vector field of $H_S$ in $S$.
Rearranging the above equation gives:
\[u^*(g)\frac{\partial u'}{\partial s} + j u^*(g)\frac{\partial u'}{\partial t} = -jX_{H_S}.\]
Now locally around the point $q$, we can choose a reparameterization
of the coordinates $(s,t)$ to new coordinates $(s',t')$ so that $u'$
satisfies:
\[\frac{\partial u'}{\partial s'} + j \frac{\partial u'}{\partial t'} = -jX_{H_S}\]
(i.e. $\frac{\partial s'}{\partial s} = \frac{\partial t'}{\partial t} = \frac{1}{u^*(g)}$ and 
$\frac{\partial t'}{\partial s} = \frac{\partial s'}{\partial t} = 0$ ).
The above equation is Floer's equation which doesn't have a maximum
by \cite[Lemma 1.5]{Oancea:survey}.
This gives us a contradiction as we assumed $f$ had a maximum at $q$.
\qed

\bigskip
If $H_s$ is a non-decreasing sequence of Lefschetz admissible Hamiltonians, then
we can use the same methods as above to prove a maximum principle for
$u$ satisfying the parameterized Floer equations
\[ \partial_s u + J_{s,t}(u(s,t)) \partial_t u = \nabla^{g_t} H_{s,t}.\]
Note we also have a maximum principle in the vertical direction as well.  We have that the region $A$ as defined in 
Definition \ref{defn:lefschetzadmissiblehamiltonian} looks like
$\partial F \times [1,\infty) \times \widehat{S}$. 
Let $r_F$ be the coordinate for $[0,\infty)$ in this
product. Let $\pi_1 : A \twoheadrightarrow \partial F \times [1,\infty)$
be the natural projection. If a Floer trajectory $u$
has an interior maximum with respect to $r_F$ for $r_F$ large, then
 $\pi_1 \circ u$ satisfies Floer's equations on $F$ and
hence has no maximum by \cite[Lemma 1.5]{Oancea:survey}.
This gives us a contradiction.
Hence $r_F \circ u$ has no maximum for $r_F$ large.
The above maximum principles and the regularity
result from Theorem \ref{thm:admissiblecomplexstr}
ensures that $SH_*(\widehat{E},H)$ is well defined. 

\begin{defn} \label{defn:periodspectra}
Let $M$ be a manifold with contact form $\alpha$.
Let \\ $S : \{\nsmbox{Reeb orbits}\} \rightarrow \R$,
$S(o) := \int_o \alpha$.
Then the {\bf period spectrum} ${\mathcal S}(M)$ is
the set $\nsmbox{im}(S) \subset \R$. We say that the
period spectrum is discrete and injective if
the map $S$ is injective and the period
spectrum is discrete in $\R$.
\end{defn}
\begin{defn} \label{defn:actionspectrum}
Let $H$ be a Hamiltonian on a symplectic manifold $M$.
Then the {\bf action spectrum} ${\mathcal S}(H)$ of
$H$ is defined to be:
\[{\mathcal S}(H) := \left\{ A_H(o) : o \nsmbox{ is a 1-periodic orbit of } X_H \right\}.\]
$A_H$ is the action defined in section \ref{section:symplectichomology}.
\end{defn}

We let $F$ be a smooth
fibre of $(E,\pi)$ and $\Theta_F := \Theta|_F$. Also we let $S$
be the base of this fibration.
Let $r_S$ and $r_F$ be the ``cylindrical'' coordinates on
$\hat{S}$ and $\hat{F}$ respectively
(i.e. $\omega_S = d(r_S\theta_S)$ on the cylindrical end at infinity and similarly
with $r_F$).
Let $W$ be some connected component of the boundary of $S$.
Let
$C := \pi^{-1}(W) \times [1,\infty)$.
Note: we will sometimes
write $r_S$ instead of $\pi^*r_S$ so that calculations
are not so cluttered. We hope that this will
make things easier to understand for the reader.

The boundary of $E$ is a union of $2$ manifolds whose boundaries
meet at a codimension $2$ corner. We can smooth out this corner
so that $E$ becomes a compact convex symplectic manifold $M$
such that the completion $\widehat{M}$ is exact symplectomorphic
to $\widehat{E}$. This means we can view $M$ as an exact submanifold
of $\widehat{E}$. We will let $\partial M \times [1,\infty)$ be
the cylindrical end of $\widehat{E} = \widehat{M}$ and we will
let $r$ be the coordinate for the interval $[1,\infty)$.
We will assume that the period spectrum of $\partial M$ is discrete and injective.
Let $\varrho_p : \widehat{E} \rightarrow \R$
be an admissible Hamiltonian on $\widehat{M}=\widehat{E}$ with slope $p$ with respect to the cylindrical
end $\partial M \times [1,\infty)$ where $p$ is a positive integer. 
We will also assume
that $\varrho_p < 0$ inside $M$ and that $\varrho_p$ tends to $0$ 
in the $C^2$ norm inside $M$
as $p$ tends to infinity, and that $\varrho_p = h_p(r)$ in the cylindrical end.
We assume that $h_p'(r) \geq 0$ for all $r$ and $h_p'(r) = p$ for
$r \geq 2$. We also assume that $h_p''(r) \geq 0$ for all $r$.
We can perturb the boundary of $M$ to ensure
that no positive integer is in the period spectrum of $\partial M$
and hence $p$ is not in the action spectrum.
Hence the family $(\varrho_p)_{p \in \N_+}$ is a cofinal family
of admissible Hamiltonians.
\begin{theorem} \label{thm:admissiblecofinal}
There is a cofinal family of Lefschetz admissible Hamiltonians $K_p : \widehat{E} \rightarrow \R$ 
and a family of almost complex structures $J_p \in {\mathcal J_{\nsmbox{reg}}}(\widehat{E},K_p)$
such that for $p \gg 0$:

\begin{enumerate}
\item 
The periodic orbits of $K_p$ of positive action are in 1-1
correspondence with the periodic orbits of $\varrho_p$. This correspondence
preserves index. Also the moduli spaces of Floer trajectories
are canonically isomorphic between respective orbits.

\item $K_p < 0$ on $E \subset \widehat{E}$.

\item $K_p|_E$ tends to $0$ in the $C^2$ norm on $E$ as $p$ tends to
infinity.

\end{enumerate}

\end{theorem}

This theorem implies that:
\begin{equation} \label{equation:lefschetzisomorphism}
\underset{p}{\varinjlim} \nsmbox{ } SH_*^{[0,\infty)}(K_p) = \underset{p}{\varinjlim} \nsmbox{ } SH_*(\varrho_p)
\end{equation}
$SH_*^{[0,\infty)}(K_p) := SH_*(K_p) / SH_*^{(-\infty,0)}(K_p)$ where
$SH_*^{(-\infty,0)}$ is the symplectic homology group generated
by orbits of negative action.
We also have:
\begin{equation} \label{equation:positiveaction}
\underset{p}{\varinjlim} \nsmbox{ } SH_*(K_p) = \underset{p}{\varinjlim} \nsmbox{ } SH_*^{[0,\infty)}(K_p)
\end{equation}
This is because there exists a cofinal family of Lefschetz admissible Hamiltonians
$G_p$ such that:
\begin{enumerate}
\item $G_p < 0$ on $E \subset \widehat{E}$.

\item $G_p|_E$ tends to $0$ in the $C^2$ norm on $E$ as $p$ tends to
infinity.

\item \label{item:Gpositiveaction}
All the periodic orbits of $G_p$ have positive action.
\end{enumerate}
Property (\ref{item:Gpositiveaction}) of $G_p$ will follow from Lemma
\ref{lemma:lefschetzactionbounds}.
Using the fact that both $K_p$ and $G_p$ are cofinal,
tending to $0$ in the $C^2$ norm on $E$ and are non-positive on $E$,
there exist sequences $p_i$ and $q_i$ such that:
\[K_{p_i} \leq G_{q_i} \leq K_{p_{i+1}}\]
for all $i$. Hence:
\[\underset{p}{\varinjlim} \nsmbox{ } SH_*^{[0,\infty)}(G_p) = \underset{p}{\varinjlim} \nsmbox{ } SH_*^{[0,\infty)}(K_p).\]
Property (\ref{item:Gpositiveaction}) of $G_p$ implies:
\[\underset{p}{\varinjlim} \nsmbox{ } SH_*^{[0,\infty)}(G_p) = \underset{p}{\varinjlim} \nsmbox{ } SH_*(G_p).\]
This gives us equation (\ref{equation:positiveaction}).
Combining this with equation (\ref{equation:lefschetzisomorphism}) gives:
\[\underset{p}{\varinjlim} \nsmbox{ } SH_*(K_p) = \underset{p}{\varinjlim} \nsmbox{ } SH_*(\varrho_p).\]
This proves Theorem \ref{theorem:lefschetzcofinalfamily}.

Before we prove Theorem \ref{thm:admissiblecofinal},
we need two preliminary Lemmas.
We need a preliminary Lemma telling us something
about the flow of a Lefschetz admissible Hamiltonian.
We let $H = \pi^*H_S + \pi_1^*H_F$ be as in
Definition \ref{defn:lefschetzadmissiblehamiltonian}.
We assume that the slope of $H_S$ and $H_F$ is strictly less
than some constant $B>0$.
We set $H_F$ to be zero in $F$,
and $H_F$ to be equal to
$h_F(r_F)$ in the region $r_F \geq 1$
such that $h_F'(r_F) \geq 0$ and $h_F''(r_F) \geq 0$.
We also assume that for some very small $\epsilon > 0$,
$h_F'$ is constant for $r_F > \epsilon$
and not in the period spectrum of $\partial F$ so that
all the orbits lie in the region $r_F \leq \epsilon$.
We define $H_S$ in exactly the same way so that
it is zero in $S$ and equal to $h_S(r_S)$ on the cylindrical
end of $\widehat{S}$ where $h_S$ has the same properties as $h_F$.
The action of an orbit of $H_F$ in the cylinder $r_F \geq 1$
is $r_F h_F'(r_F) - h_F(r_F)$  and similarly the action of an orbit
of $H_S$ in $r_S \geq 1$ is $r_S h_S'(r_S) - h_S(r_S)$, so
we can choose $\epsilon$ small enough so that the actions of the orbits
lie in the interval $[0,B]$ because the slope of $H_S$ and $H_F$ is less than $B$.
We have from Section \ref{section:lefschetzfibrationproofs},
$\theta = \Theta + k\pi^*\theta_S$
where $\Theta$ is the $1$-form associated to
the Lefschetz fibration (it is a $1$-form such that
$\Theta|_F$ makes each fibre $F$ into a compact convex
symplectic manifold. Also $\theta_S$ is the $1$-form making
the base $S$ into a compact convex symplectic manifold.
The constant $k$ is some large constant.
\begin{lemma} \label{lemma:lefschetzactionbounds}
For $k$ large enough,
there exists a constant $\Xi$ depending only on $E$ and $\theta$
(not on $H$) such that the action of any orbit
of $H$ is contained in the interval
$[0,\Xi B]$.
\end{lemma}
\proof
Inside $E$, we have that the Hamiltonian is $0$
so all the orbits have action $0$ there. In the
region $A$  as defined in
Definition \ref{defn:lefschetzadmissiblehamiltonian}, we have that
the orbits come in pairs $(\gamma,\Gamma)$
where $\gamma$ is an orbit from $H_S$ and
$\Gamma$ is an orbit of positive action from $H_F$.
The action of $(\gamma,\Gamma)$ is the sum of the actions
of $\gamma$ and $\Gamma$. Both these actions are positive.
Also their actions are bounded above by $B$.

So we only need to consider orbits outside the region
$A \cup E$. The Hamiltonian $\pi_1^*H_F$ is zero in this
region so we only need to consider $\pi^*H_S$. We will
consider the orbits of $\pi^*H_S$ in the region $r_S \geq 1$.
In this region, there are no singular fibres of the Lefschetz fibration,
so we have a well defined plane field $P$ which is the $\omega$-orthogonal
plane field to the vertical plane field which is the plane field tangent to the
fibres of $\pi$. The Hamiltonian flow only
depends on $\omega|_{P}$ and not the vertical plane field because
$\pi^*H_S$ restricts to zero on the vertical plane field.
The symplecic form $\omega|_{P}$
is equal to $G k \pi^*d\theta_S|_P$ for some function $G>0$.
This means that the Hamiltonian vector field associated
to $\pi^*H_S$ is $\frac{1}{G}$ times the horizontal lift
of the Hamiltonian vector field associated to $H_S$ in $S$.
Let $V$ be this horizontal lift.
The construction of the completion of a Lefschetz fibration
before Definition \ref{defn:halfconvexlefschetzcompletion}
ensures that the region $r_S \gg 1$ is a product
$W \times [1,\infty)$ where $r_S$ parameterizes the second factor
of this product and $\Theta$ is a pullback of a $1$-form on $W$
via the natural projection $W \times [1,\infty) \rightarrow W$.
This means that $\Theta$ is invariant under translations in the
$r_S$ direction (i.e. under the flow of the vector field
$\frac{\partial}{\partial r_S}$ which is $\frac{1}{r_S}$ times
the horizontal lift of $\lambda_S$ where $\lambda_S$
is the Liouville flow in $\widehat{S}$). We also have that
$d\theta_S$ is invariant under translations in the $r_S$ direction
(i.e. under the flow of $\frac{1}{r_S}\lambda_S$).
Hence the symplectic structure $\omega$ is also invariant under
translations in the $r_S$ direction for $r_S \gg 0$. 
This means that the function $G$ is bounded above and below by positive constants
as the symplectic structure 
 is invariant under translations in the $r_S$ direction
 and if we travel to infinity in the fibrewise direction (i.e. if we travel
into the region $A$), then $G=1$.
We want bounds on
the function $V(\theta)$ because the function $G$ is bounded.
Let $Y$ be the Hamiltonian flow of $r_S$ in $\widehat{S}$ and
let $\tilde{Y}$ be its horizontal lift to $P$.
We have that $Y(\theta_S)=1$. This means that
$\tilde{Y}(\pi^*\theta_S) =1$. We also have
that $\tilde{Y}(\Theta)$ is bounded because $\Theta$
is invariant in the $r_S$ direction for $r_S$ large
and $\tilde{Y}(\Theta) = 0$ if we are near infinity in the fibrewise
direction. 
We choose the constant $k$ large enough so that
$\tilde{Y}(\theta) = \tilde{Y}(\Theta) + k\tilde{Y}(\pi^*\theta_S) > 0$.
This function is also bounded above because $\tilde{Y}(\Theta)$ is
bounded and $k\tilde{Y}(\pi^*\theta_S) = kY(\theta_S) = k$.
This choice of $k$ only depends on the Lefschetz fibration and not
on $H$.
Now, $V = h_S'(r_S)\tilde{Y}$. Because $h_S'$ bounded below by $0$, we have
that $V(\theta)$ is bounded below by $0$ and bounded above by some constant
multiplied by the slope of $H_S$. All the orbits of $H$ lie in some
compact set where $H$ is $C^0$ small, so the action of an orbit is
near $\int_o V(\theta) dx$ where the integral is taken over an orbit $o$
and $dx$ is the volume form on $o$ giving it a volume of $1$. 
This means that the action of these orbits is in the interval $[0,\Xi B]$ for
some constant $\Xi$. This completes our theorem.
\qed

\bigskip

The manifold $\widehat{M} = \widehat{E}$ has a cylindrical end
$\partial M \times [1,\infty)$. We let $r$ be the radial coordinate
of this cylindrical end. The we define set $\{r \leq R\}$
to be equal to $M \cup (\partial M \times [1,R])$. We define
the sets $\{r_F \leq  R\}$ and $\{r_S \leq R\}$ in a similar way.

\begin{lemma} \label{lemma:lefschetzcontactboundarysize}
There exists a constant $\varpi > 0$ such that
for all $R \geq 1$, we have that $\{r \leq R\} \subset \{r_S \leq \varpi R\}$
and $\{r \leq R\} \subset \{r_F \leq \varpi R\}$.
\end{lemma}
\proof
We will deal with $r_S$ first.
The level set $r = R$ is equal to the flow of $\partial M$
along the Liouville vector field $\lambda$ for a time $\log (R)$.
Hence, all we need to do is show that $dr_S(\lambda)$ is bounded above
by $e^\varpi r_S$. This means that if $p$ is a point in $\partial M$,
then the rate at which $r_S(p)$ increases as we flow $p$ along $\lambda$ is
bounded above by $e^\varpi r_S(p)$. Hence if we flow $p$ for a time
$\log(R)$ to a point $q$, then $r_S(q) \leq \varpi R$ which is our result.

We will now show $dr_S(\lambda)$ is bounded above
by $e^\varpi r_S$ to finish the first part of our proof. We let $\Theta$ be a $1$-form associated
to $E$ as constructed before Definition 
\ref{defn:halfconvexlefschetzcompletion}.
Then $\theta = \Theta + \pi^*\theta_S$ where $\theta_S$
is a convex symplectic structure for the base $\widehat{S}$.
We have that $\omega = d\Theta + \pi^*d\theta_S$.
The construction before Definition \ref{defn:halfconvexlefschetzcompletion}
ensures that the region $r_S \gg 1$ is a product
$W \times [1,\infty)$ where $r_S$ parameterizes the second factor
of this product and $\Theta$ is a pullback of a $1$-form on $W$
via the natural projection $W \times [1,\infty) \rightarrow W$.
This means that $\Theta$ is invariant under translations in the
$r_S$ direction. Hence $d\Theta$ is also invariant under these
translations. Also $\pi^* d\theta_S$ is invariant under translations
in the $r_S$ direction. All of this means that the vector
field $V$ defined as the $\omega$-dual of $\Theta$ is invariant
under these translations for $r_S$ large. This implies that
$dr_S(V)$ is bounded. 

Let $V'$ be the $\omega$-dual of $\pi^*\theta_S$.
Let $\lambda_S$ be the Liouville vector field in $\widehat{S}$.
Then $V' = G L$ where $L$ is the horizontal lift
of $\lambda_S$ and $G : \widehat{E} \rightarrow \R$ is
defined in the proof of Lemma \ref{lemma:lefschetzactionbounds}.
The proof of Lemma \ref{lemma:lefschetzactionbounds}, tells
us that $G$ is a bounded function.
Also, $dr_S(\lambda_S) = r_S$, hence
\[dr_S(V') = G dr_S(L) = G dr_S(\lambda_S) = Gr_S\]
Hence $dr_S(V')$ is bounded above by some constant
multiplied by $r_S$.
Finally, we have that $\lambda = V + V'$ which means
that there exists a $\varpi>0$ such that $dr_S(\lambda)$
is bounded above by $e^\varpi r_S$.

We will now deal with $r_F$. This is slightly
more straightforward because the Lefschetz fibration
is a product $\partial F \times [1,\infty) \times \widehat{S}$
and $\theta$ splits up in this product as $\Theta + \pi^*\theta_S$,
where we can view $\Theta$ as $1$-form on $\partial F \times [1,\infty)$.
We need to bound $dr_F(\lambda)$. In this case, because everything
splits in this product, we have that $dr_F(\lambda) = dr_F(\Lambda)$
where $\Lambda$ is the $\omega$-dual of $\Theta$.
This is equal to $r_F \leq e^\varpi r_F$ as $\varpi > 0$. 
Hence we have that $r \leq R$ implies that
$r_F \leq \varpi R$.
\qed

\proof of Theorem \ref{thm:admissiblecofinal}.
%
%
Let $\varrho_p$ be the Hamiltonian as above. We will write $\varrho=\varrho_p$ for simplicity.
The idea of the proof is to modify the Hamiltonian $\varrho$ outside some large compact set so that
it becomes Lefschetz admissible and in the process only create
orbits of negative action without changing the orbits of $\varrho$ or the Floer
trajectories connecting orbits of $\varrho$. We will do this in three sections.
In section (a), we will modify $\varrho$ to a Hamiltonian $\varsigma$ 
so that it becomes constant outside a large
compact set $\kappa$ while only adding orbits of negative action. This is 
exactly the same as the construction due to Hermann \cite{Hermann:holomorphic}.
In section (b) we will consider a Lefschetz admissible Hamiltonian $L$ which
is $0$ in the region $\kappa$, but has action bounded above so that
the orbits of $L + \varsigma$ outside $\kappa$ have negative action.
We define our cofinal family $K_p := L + \varsigma$.
(c) we ensure that the Floer trajectories and pairs
of pants satisfying Floer's equation connecting
orbits of positive action stay inside the region
$r \leq 2$.

(a)
We have that $p$ is the slope of the Hamiltonian $\varrho$ and this
is not in the period spectrum of $\partial M$. Hence, we define
$\mu := \mu(p) > 0$ to be smaller than the distance between $p$
and the action spectrum.
Define:
\[A = A(p) := 3p / \mu > 1.\]
\begin{figure}[b]
\end{figure}
\begin{figure}[b]
\end{figure}
\begin{figure}[b]
\end{figure}

We can assume that $A>4$ because we can choose
$\mu$ to be arbitrarily small.
%
Remember that $\widehat{E}=\widehat{M}$ where $M$ is a compact convex symplectic manifold,
and that $r$ is the radial coordinate for the cylindrical end of $\widehat{M}$.
We define $\varsigma$ to be equal to $\varrho$ on $r \leq A - 1$.
On the region $r \geq 1$, we have that $\varrho$ is equal to $h_p(r)$. 
We will
just write $h$ instead of $h_p$.
Set $\varsigma = k(r)$ for $r \geq 1$ with non negative derivative.
This means that in the region $1 \leq r \leq A - 1$ we have that
$h(r) = k(r)$.
Hence in $r \leq A - 1$ we have that $k''(r) \geq 0$
and $k'(r) \geq 0$, and in the region $2 \leq r \leq A - 1$
we have $k'(r) = p$. Also we have that $\varsigma$ is $C^2$ small and negative for $r$ near $1$.
Because $\varsigma$ is $C^2$ small, we can also assume that $p$ is large enough so that
for $r$ near $1$, $k' \ll p$. Because $\varrho_p$ is cofinal, we can assume that $p$
is large enough so that $h(2) = k(2) > 0$.
Both these previous facts mean that $p(A-2) < k(A-1) < p(A - 1)$.
Outside this region, we define $k$ to be a function with the following constraints: 
For $r \geq A$ set $k(r)$ to be constant
and equal to $C$ where $C = p (A - 1)$.
In the region $A - 1 \leq r$, $k'' \leq 0$.
We assume that $k' \geq 0$ for all $r \geq 1$.
Here is a picture:

\begin{fig} \label{fig:graphslope}

 \[
 \xy
 (5,0)*{}="A"; (5,-50)*{}="B";
 "A"; "B" **\dir{-};
 (0,-35)*{}="C"; (60,-35)*{}="D";
 "C"; "D" **\dir{-};
 (5,-40)*{}="E"; (60,-40)*{}="F";
 (5,-10)*{}="E1"; (60,-10)*{}="F1";
 "E1"; "F1" **\dir{.};
 (40,-35)*{}="VL1"; (40,-10)*{}="VL2";
 "VL1"; "VL2" **\dir{.};
 (45,-35)*{}="VR1"; (45,-10)*{}="VR2";
 "VR1"; "VR2" **\dir{.};
 (18,-35)*{}="G"; 
 "G"; "G" + (7,3); **\crv{"G"+(5,0) & "G"+(2,0)};
 "G"+(7,3); "VL2"+(0,-3) **\dir{-};
 "VL2"+(0,-3); "VR2"; **\crv{"VL2"+(3,0) & "VR2" + (-5,-3)};
 "VR2"; "F1" **\dir{-};
 "E1"+(-2,0) *{C};
 "VL1"+(-4,-2) *{\scriptstyle A-1};
 "VR1"+(2,-2) *{\scriptstyle A};
 "G" + (0,-2) *{\scriptstyle 1};
 "G" + (7,-2) *{\scriptstyle 2};
 "VL2" + (-10,2) *{{\bf \varsigma}};
 "VL2" + (-8,-12) *{\scriptstyle p}
 \endxy
 \]

\end{fig}

We want to show that the additional orbits
of $\varsigma$ only have
negative action. All these orbits lie in the region $r \geq 2$.
In fact because $p$ is not in the action spectrum, they lie
in the region $r \geq A-1$.
In the region $\{r : p - \mu < k'(r) \leq p\}$,
we have that $\varsigma$ has no periodic orbits. Also, the action
of a periodic orbit is $k'(r)r - k(r)$. Combining these two
facts implies that the action of a periodic orbit
in the region $2 \leq r$ is less than
\[(p - \mu)r - k(r) \leq (p - \mu)A - p(A - 2)\]
\[ = -\mu A + 2p = - \mu  \frac{3p}{\mu} + 2p = -p < 0\]
Hence we have a Hamiltonian $\varsigma$ equal to $\varrho$ in the region $r \leq 2$
and such that it is constant and equal to $C = p(A - 1)$
in the region $r \geq A-1$ and such that all the additional periodic orbits
created have negative action.

(b) Lemma \ref{lemma:lefschetzactionbounds} tells us that there
exists a cofinal family of Lefschetz admissible Hamiltonians
$\Lambda_p$ such that the
action spectrum of $\Lambda_p$ is bounded above by some
constant $\Xi$ multiplied by the slope of $\lambda_p$.
We can assume that both the slopes of $\lambda_p$ are equal to
$\sqrt(p)$ (if $\sqrt{p}$ is in the action spectrum of the
fibre or the base, then we perturb this value slightly to ensure
that $\Lambda_p$ has orbits in a compact set).
This means that the action of $\Lambda_p$ is bounded above
by $\Xi \sqrt{p}$.
The Hamiltonian $\Lambda_p$ is equal to zero in $E$.
We will now define a Hamiltonian
$L_p$ as follows:
We let $\varpi$ be defined as in Lemma 
\ref{lemma:lefschetzcontactboundarysize}.
Set $L_p = 0$ in the region
$\{r_S \leq \varpi A \} \cap \{r_F \leq \varpi A\}$.
In the region $\{r_S \geq 1\} \cup \{r_F \geq 1\}$,
we have that $\Lambda_p$ is a function of the form
$\pi_1^*h_F(r_F) + \pi^*h_S(r_S)$.
Here, $\pi_1$ is the natural projection:
$\partial F \times [1,\infty) \times \widehat{S} \rightarrow
\partial F \times [1,\infty)$ (this is the same as the
projection defined just before Definition
\ref{defn:lefschetzadmissiblehamiltonian}).
So, we set the function $\pi_1^*h_F(r_F)$ to be zero outside
the domain of definition of $\pi_1$. Also,
$\pi^*h_S$ is zero outside the region $r_S \geq 1$.
We define $L_p$ to be 
\[\pi_1^*h_F(r_F - \varpi A) + \pi^*h_S(r_S - \varpi A)\]
in the region  $\{r_S \geq \varpi A\} \cup \{r_F \geq \varpi A\}$.
Hence we have a well defined function $L_p$.
Because $L_p$ has scaled up, we have that the
action spectrum of $L_p$ is equal to $\varpi A$ multiplied
by the action spectrum of $\Lambda_p$. Hence,
we have that the action spectrum of $L_p$
is bounded above by $\varpi A \Xi \sqrt{p}$.

Because $\{r \leq A\} \subset \{r_S \leq A\} \cap \{ r_F \leq A\}$,
we can add $L_p$ to $\varsigma$ without changing the orbits of $\varsigma$
in the region $r \leq A$. Also, the action of the orbits
of $\varsigma + L_p$ in the region $r \geq A$ is
bounded above by $\varpi A \Xi \sqrt{p} - p(A - 1)$.
So for $p$ large enough we have that the additional
orbits added are of negative action.

(c) We choose an almost complex structure 
$J \in {\mathcal J}^h(\widehat{E})$ such that
on some neighbourhood of the hypersurface $r = 2$,
$J$ is admissible. Then \cite[Lemma 7.2]{SeidelAbouzaid:viterbo}
and the comment after this Lemma ensure that
no Floer trajectory or pair of pants satisfying
Floer's equation connecting orbits inside $r < 2$
can escape $r \leq 2$.
Hence our Hamiltonian $K_p := \varsigma + L_p$ has all the
required properties.


\subsection{A better cofinal family for the Lefschetz fibration} \label{section:bettercofinal}

In this section we will prove Theorem \ref{theorem:lefschetzcofinalsequence}.
We consider a compact convex Lefschetz fibration $(E,\pi)$ fibred over
the disc $\D$. 
Basically the cofinal family is such that
$H_F=0$. This means that the boundary of $F$ does not contribute
to symplectic homology of the Lefschetz fibration.  The key idea
is that near the boundary
of $F$ the Lefschetz fibration looks like a product
$\D \times \nsmbox{nhd}(\partial F)$ and because symplectic homology
of the disc is $0$ we should get that the boundary contributes nothing.
Statement of Theorem \ref{theorem:lefschetzcofinalsequence}:
\[SH_*(E) \cong SH^{\nsmbox{lef}}_*(E).\]
We will define $F,S (=\D),r_S,r_F,\pi_1$ as in the previous section.
This means that the compact convex sympectic manifold $F$ is a fibre of $E$ and
$S$ is the base which in this section is equal to $\D$.
We also have that $r_S$ is a radial coordinate for the cylindrical end
of $\widehat{S}$ which we also identify with $\pi^*r_S$.
The map $\pi_1$ is the natural projection
$(\partial F \times [1,\infty)) \times \widehat{S} \twoheadrightarrow
(\partial F \times [1,\infty))$ where
$(\partial F \times [1,\infty)) \times \widehat{F}$ is a subset of $\widehat{E}$.
The function  $r_F$ is a radial coordinate for the cylindrical end
of $\widehat{F}$ which we also identify with $\pi_1^*r_F$.
Before we prove Theorem \ref{theorem:lefschetzcofinalsequence},
we will write a short lemma on the $\Z$ grading of $SH_*(E)$.
\begin{lemma} \label{lemma:chernclassoffibre}
Let $\widehat{F} := \pi^{-1}(a) \subset \widehat{E}$ ($a \in \D$). 
Suppose we have trivialisations of ${\mathcal K}_{\widehat{E}}$ 
and ${\mathcal K}_{\widehat{S}}$ (these are the canonical
bundles for $\widehat{E}$ and $\widehat{S}$ respectively);
these naturally induce a trivialisation of ${\mathcal K}_{\widehat{F}}$
away from $F$. If we smoothly move $a$, then this
smoothly changes the trivialisation.
\end{lemma}
\proof of Lemma \ref{lemma:chernclassoffibre}.

We choose a $J \in {\mathcal J}^h(E)$.
The bundle $E$ away from $E^{\nsmbox{crit}}$ has a connection induced
by the symplectic structure. Let $A \subset \widehat{E}$ be defined as in
Definition \ref{defn:lefschetzadmissiblehamiltonian}. Let $U$
be a subset of $A$ where 
\begin{enumerate}
\item $\pi$ is $J$ holomorphic.
\item \label{item:extendtrivialisation}
$U$ is of the form $r \geq K$ where $r$ is
the coordinate for $[1,\infty)$ in $A$
(see Definition \ref{defn:lefschetzadmissiblehamiltonian}).
\end{enumerate}
This means that in $U$, we have that the horizontal plane
bundle ${\mathbb H}$ is $J$ holomorphic. Choose a global holomorphic section
of  ${\mathcal K}_{\widehat{S}}$ and lift this to a section $s$ of 
${\mathbb H}$.
Choose a global holomorphic section $t$ of  ${\mathcal K}_{\widehat{E}}$.
The tangent bundle of $\widehat{F}$ is isomorphic to the
$\omega$-orthogonal bundle $T$ of ${\mathbb H}$. 
This is also a holomorphic bundle.
Let $\Lambda^k T$ be the highest exterior power of $T$.
There exists a unique holomorphic section $w$ of $\Lambda^k T$ such that
$s \wedge w = t$. Hence, $w$ is our nontrivial holomorphic section
of $T$ in $U \cup \widehat{F}$.
This can be extended to $A \cup \widehat{F}$ by property (\ref{item:extendtrivialisation}).
\qed

In the following proof, whenever we talk about indices of
orbits of $\widehat{F}$ outside $F$, we do this with
respect to the trivialisation of Lemma \ref{lemma:chernclassoffibre}
above. We do not deal with orbits inside $F$ so this trivialisation
is sufficient.

\proof of Theorem \ref{theorem:lefschetzcofinalsequence}.
We start by defining a cofinal family of Lefschetz admissible Hamiltonians
$H^\lambda = \pi^*H_S^\lambda + \pi_1^*H_F^\lambda$.
To avoid cluttered notation, we suppress the $\lambda$ and
just write $H, H_S, H_F$ instead of $H^\lambda, H_S^\lambda, H_F^\lambda$
unless we need to explicitly deal with $\lambda$.
We assume that the period spectrum of $\partial F$ is discrete and injective
and also that the Reeb orbits of $\partial F$ are non-degenerate.
We assume that $H_F=0$ on $F$ and is equal to $h_F(r_F)$ outside
$F$ with $h_F'(r_F) > 0$ and $h_F''(r_F) \geq 0$ when $r_F > 1$.
The orbits of $H_F$ consist of constant orbits in $F$
and ${\mathbb S}^1$ families of orbits corresponding to periodic
Reeb orbits outside $F$. We can perturb
$H_F$ by a very small amount outside
$F$ so that each ${\mathbb S}^1$ family of Reeb orbits becomes
a pair of non-degenerate orbits (see \cite[Section 3.3]{Oancea:survey}).
Hence, we have a Hamiltonian $H_F$ which is
equal to $0$ inside $F$ and all its
orbits outside $F$ are non-degenerate.
We set the slope of $H_F$ at infinity to be equal to 
$\lambda \notin {\mathcal S}(\partial F)$.
The completion $\widehat{\D}$ of the disc $\D$ is symplectomorphic to $\C$
with the standard symplectic structure.
We have that $r_S(z) = |z|^2$ where $z \in \C$.
The function $r_S$ is defined only on the cylindrical end, but we will
extend it to the interior of $\D$ by the function $|z|^2$.
We set $H_S=\kappa(\lambda)r_S$ on $\C$ 
where $\kappa(\lambda)$ is some function of the slope $\lambda$ of $H_F$
such that it is never a multiple of $\pi$.
Here $H_S$ has exactly one periodic orbit at $0 \in \C$ of index $(2 a + 1)$
where $a$ is the integer satisfying $a \pi < \kappa(\lambda) < \pi(a + 1)$
(see \cite[Section 3.2]{Oancea:survey}).
For each $\lambda$, we choose $\kappa(\lambda)$ to be large enough so that
the index of the only orbit of $H_S$ is greater than $M + \lambda + 1$ where
$M = -\nsmbox{ind}(\Gamma)$ such that $\Gamma$ is an orbit of $H_F$ of lowest index.
We also assume that as $\lambda$ tends to infinity, $\kappa(\lambda)$ also tends to infinity.
The Hamiltonians $H^\lambda$ form a cofinal family of Lefschetz admissible Hamiltonians.
Hence for some $J \in {\mathcal J}^h(\widehat{E})$, we have $SH_*(E) = \underset{\lambda}{\varinjlim} \nsmbox{ } SH_*(E,H^\lambda,J)$.

We will also define a half Lefschetz admissible Hamiltonian $\bar{H}^\lambda := \pi^*H_S^\lambda + \pi_1^*\bar{H}_F^\lambda$.
When appropriate, we will write $\bar{H}, \bar{H}_F$ instead of $\bar{H}^\lambda, \bar{H}_F^\lambda$.
We define $\bar{H}_F = 0$ inside $F$ and $\bar{H}_F := \bar{h}_F(r_F)$ in the region $\{r_F \geq 1\}$
where the derivative $\bar{h}_F '$ is so small, that $\bar{H}_F$ has no periodic orbits
in the region $\{r_F > 1\}$.
The Hamiltonians $\bar{H}^\lambda$ form a cofinal family of half admissible Hamiltonians, hence
$SH_*^{\nsmbox{lef}}(E) = \underset{\lambda}{\varinjlim} \nsmbox{ } SH_*(E,\bar{H}^\lambda,J)$.

We will construct a natural continuation map $SH_*(\bar{H}^\lambda,J) \rightarrow SH_*(H^\lambda,J)$
such that it is an isomorphism in all degrees less than $\lambda$.
This will prove the theorem for the following reason:
Because continuation maps are natural, we can take direct limits with respect to $\lambda$, so that
we get a map $SH_*(E) \rightarrow SH_*^{\nsmbox{lef}}(E)$. This map must be an isomorphism
because if we choose an integer $b$, then for $\lambda > b$ we have
$SH_b(\bar{H}^\lambda,J) \rightarrow SH_b(H^\lambda,J)$ is an isomorphism
which implies that $SH_b(E) \rightarrow SH_b^{\nsmbox{lef}}(E)$ is an isomorphism
(because we are taking a direct limit as $\lambda$ tends to infinity).
Hence $SH_*(E) \rightarrow SH_*^{\nsmbox{lef}}(E)$ is an isomorphism.

We will now show that the continuation map $SH_*(\bar{H}^\lambda,J) \rightarrow SH_*(H^\lambda,J)$
is an isomorphism in all degrees less than $\lambda$.
From now on, we will write $H, H_S, H_F, \bar{H}, \bar{H}_F$
instead of $H^\lambda, H_S^\lambda, H_F^\lambda, \bar{H}^\lambda, \bar{H}_F^\lambda$.
The region $\{r_F \geq 1\}$ is a product $([1,\infty) \times \partial F) \times \C$
and the orbits of the Hamiltonian $H$ come in pairs $(\gamma, \Gamma)$,
where $\gamma$ is the orbit of $H_S$ and $\Gamma$ is a non-constant orbit of $H_F$.
The index of this orbit is the sum $\nsmbox{ind}(\gamma) + \nsmbox{ind}(\Gamma)$,
hence its index is greater than $\lambda + 1$.
This means that all the orbits of index $\leq \lambda + 1$ are disjoint from the region $\{r_F \geq 1\}$.
We have that $H_F$ is not quite a function of $r_F$ in the region $\{r_F \geq 1\}$,
as we perturbed it so that it had non-degenerate orbits.
Having said that we can assume that for some $\delta > 0$, we have 
$H_F = h_F(r_F)$ in the region $\{1 \leq r_F \leq 1 + \delta\}$.
We also assume that the Lefschetz admissible almost complex structure $J$
is of the form $j + J_F$ in the region $\{r_F \geq 1\}$
viewed as a subset of $\widehat{E}$, where $j$ is an complex structure on $\C$ and
$J_F$ is convex on the cylindrical end of $\widehat{F}$.
Any Floer trajectory satisfying Floer's equation with respect to $(H,J)$
connecting orbits outside the region $\{r_F > 1\}$, must stay outside
this region for the following reason:
If $u : S \rightarrow \widehat{E}$ is such a curve, then let
$\bar{S} := u^{-1}(\{r_F > 1\})$.
We can project $u|_{\bar{S}} : \bar{S} \rightarrow \{r_F > 1 \}$ down from $\{r_F > 1\}$ to
$(1,\infty) \times \partial F$. Lemma 7.2 from \cite{SeidelAbouzaid:viterbo} then
tells us that this curve cannot exist. Hence $\bar{S}$ is empty and the claim follows.
A similar argument shows that if we had a pair of pants satisfying Floer's equation
with respect to Hamiltonians similar to $H$ connecting orbits outside $\{r_F \geq 1\}$,
then the curve must also be disjoint from $\{r_F \geq 1\}$.
Similarly any Floer trajectory or pair of pants satisfying Floer type equations with respect to $(\bar{H},J)$
must be disjoint from $\{r_F \geq 1\}$ because all orbits of $\bar{H}$ are disjoint from this region.
Let $K_s : \widehat{F} \rightarrow \R$ be a monotone increasing sequence of admissible Hamiltonians
joining $\bar{H}_F$ and $H_F$ such that in the region $\{1 \leq r_F \leq 1 + \delta\}$, we have
that $K_s$ is a function of $r_F$ only.
Again if we have a trajectory satisfying the Floer continuation equations with respect to $(K_s,J)$
joining orbits of $\bar{H}$ and $H$ which are disjoint from $\{r_F \geq 1\}$,
then by \cite[Lemma 7.2]{SeidelAbouzaid:viterbo} we have that this trajectory
is also disjoint from $\{r_F \geq 1\}$.
Combining all these facts, we get that the continuation map induced by $(K_s, J)$
is an isomorphism in all degrees less than $\lambda$ because
$\bar{H} = H$ outside the region $\{r_F \geq 1\}$ and all orbits of $H$ and $\bar{H}$ of degree less than
$\lambda + 1$ along with all Floer trajectories and pairs of pants connecting them
are disjoint from this region.
This completes the proof of Theorem \ref{theorem:lefschetzcofinalsequence}.
\qed


\section{$SH^{\nsmbox{lef}}_*(\widehat{E})$ and the Kaliman modification} \label{section:transferisomorphism}

In this section we prove Theorem \ref{theorem:lefschetztransferisomorphism}.
Throughout this section we assume that $E'$ and $E''$ are Lefschetz
fibrations as described in section \ref{subsection:symhomlefschetz}.
We recall the situation:
\begin{enumerate}
   \item \label{item:subfibrationlef2}
$E''$ is a subfibration of $E'$ over the same base.
   \item \label{item:monodromysupportlef2}
The support of the parallel transport maps of $E'$ are
contained in the interior of $E''$.
   \item \label{item:complexcurvecriterionlef2}
There exists a complex structure $J_{F'}$ 
(coming from a Stein domain) on 
$F'$ such that
any $J_{F'}$-holomorphic curve in $F'$
with boundary in $F''$ must be contained in $F''$.
 \end{enumerate}
We wish to prove that 
$SH_*^{\nsmbox{lef}}(E') \cong SH_*^{\nsmbox{lef}}(E'')$ as rings.

\proof of Theorem \ref{theorem:lefschetztransferisomorphism}.
Fix $\lambda > 0$. The value $\lambda$ is going to be the slope of
some Hamiltonian, we can always perturb $\lambda$
slightly so that it isn't in the action spectrum of
the boundary.
By Theorem \ref{theorem:steintoalmoststein} we can choose an almost complex
structure $J_{F',1}$ on $\widehat{F'}$ after a convex deformation away
from $F'$ such that it is convex with respect to some cylindrical end
at infinity and such that any $J_{F',1}$-holomorphic curve in $F'$
with boundary in $F''$ must be contained in $F''$. The reason is
because we can ensure that $J_{F',1} = J_{F'}$ in $F' \subset \widehat{F'}$
and that any $J_{F',1}$-holomorphic curve with boundary in $F'' \subset F'$
is contained in $F'$ by Theorem \ref{theorem:steintoalmoststein}
hence is contained in $F''$ by property
(\ref{item:complexcurvecriterionlef2}) above. Supposing we have a Hamiltonian
$H_{F'}$ which is of the form $h_{F'}(r_{F'})$ on the cylindrical end
where $r_{F'}$ is the radial coordinate and $h' \geq 0$ and $H_{F'} = 0$
elsewhere. Then, any curve (Floer cylinder or pair of pants)
with boundary in $F''$ satisfying Floer's
equations with respect to $H_{F'}$ and $J_{F',1}$ must be contained in $F''$.
We choose $h'$ small enough so that $H_{F'}$ has no periodic orbits
in the region $r_{F'} > 1$.
The convex deformation mentioned in Theorem \ref{theorem:steintoalmoststein}
fixes $F' \subset \widehat{F'}$ hence it
induces a convex deformation on $\widehat{E}$. This is because
the region where we deform $\widehat{E}$ looks like a product
$\C \times (\widehat{F'} \setminus F')$.
From now on we assume that the fibres of $\widehat{E}$ have this
almost complex structure $J_{F',1}$ with this cylindrical end.

A neighbourhood of $\partial F''$ in $F'$ is symplectomorphic to
$L:=(-\epsilon, \epsilon) \times \partial F''$ with the symplectic
form $d(r\alpha'')$. Here, $r$ is a coordinate in $(-\epsilon, \epsilon)$
and $\alpha''$ is the contact form for $\partial F''$.
We also choose $\epsilon$ small enough so that $L$ is
disjoint from the support of the parallel transport maps in $F'$.
Let $\bar{F}'' := F'' \setminus ((-\epsilon/3,0] \times \partial F'')$.
We can choose an almost complex structure $J' \in {\mathcal J}(\widehat{F'})$
with the following properties: 
\begin{enumerate}
\item There exists a $\delta >0$ 
such that any holomorphic curve meeting both boundaries of
$[-\epsilon ,-\epsilon /2] \times \partial F''$ has area greater
than $\delta$. This is true by the monotonicity lemma \cite[Lemma 1]{Oancea:kunneth}. 
\item $J'=J_{F',1}$ on $\widehat{F'} \setminus \bar{F}''$.
This means that any curve (cylinder or pair of pants)
satisfying Floer's equations with respect to $H_{F'}$ and $J'$ with
boundary in $F''$ is contained entirely in $F''$.
\end{enumerate}
Construct an almost complex structure $J$ on $\widehat{E'}$ as follows:
The parallel transport maps on $\widehat{F'} \setminus \bar{F''}$
are trivial, hence there is a region $W$ of $\widehat{E'}$
symplectomorphic to $\C \times (\widehat{F'} \setminus \bar{F''})$.
We set $J|_W$ to be the product almost complex structure $J_{\C} \times J'$
where $J_{\C}$ is the standard complex structure
on $\C$. We then extend $J|_W$ to some
$J$ compatible with the symplectic form $\omega'$
such that $\pi'$ is $J$-holomorphic outside some large compact set.
Let $H$ be a Hamiltonian of the form $\pi^*K + \pi_1^*H_{F'}$ where $K$ is
admissible of slope $\lambda$ on the base $\C$ and
$\pi_1 : W \twoheadrightarrow {\widehat{F'} \setminus \bar{F}''}$ is the natural projection map.
$J$ has the following properties:
\begin{enumerate}
\item \label{item:areaholomorphic}
Any curve $u$ satisfying Floer's equations with respect to $H$ 
meeting both boundaries of $([-\epsilon ,-\epsilon /2] \times \partial F'') \times \C \subset \widehat{E'}$
must have energy $\geq$ $\delta$.
($u$ can be a cylinder or a pair of pants).
\item \label{item:trajectoryconstraint}
Any such $u$ connecting orbits
inside $E''$ must be entirely contained in $E''$.
\end{enumerate}
Property (\ref{item:trajectoryconstraint}) is true because:
Let $u$ be a curve satisfying Floer's equations connecting
orbits in $E''$, then composing $u|_{u^{-1}(W)}$ with the natural projection
$W \twoheadrightarrow \widehat{F'} \setminus \bar{F}''$ 
gives us a curve $w$ with boundary in $F''$.
This means that $w$ is contained in $F''$,
and hence $u$ is contained in $E''$.
Also if $u$ meets both boundaries of $([-\epsilon ,-\epsilon /2] \times \partial F'') \times \C \subset \widehat{E'}$, then the projected curve $w$
 has energy $\geq \delta$ which means
that $u$ has energy $\geq \delta$.
Hence Property (\ref{item:areaholomorphic}) is true.

We perturb $H_{F'} : \widehat{F'} \rightarrow \R$ slightly so that:
\begin{enumerate}
\item It is equal to $0$ in $F''$.
\item The only periodic orbits of $H_{F'}$ are constant orbits.
\item \label{item:hfactionbound}
The action spectrum of $H|_{\widehat{F'} \setminus F''}$ is discrete and injective
and contained in $(-\delta/4,\delta/4)$.
\item We leave $H_{F'}$ alone on the cylindrical end.
\item \label{item:negativeaction}
All the orbits in $\widehat{F'} \setminus F''$ are of negative action and non-degenerate.
\item \label{item:smallslopeinmiddle}
we can ensure that $H_{F'}$ has very small positive slope
with respect to the cylindrical end of $F''$ on the region
$(\epsilon/2,\epsilon) \times \partial F'' \subset F'$.
\end{enumerate}
Let $\delta_1>0$ be the smallest distance between $0$ and
the action value of an orbit of $H_{F'}$ of negative action.
Here we fix some integer $m>0$.
We can assume that the critical points of our Lefschetz
fibration in $\C$ form a regular polygon with centre the origin.
Draw a straight line from the origin to each critical point
and let $G$ be the union of these lines.
Let $X := \frac{r}{2} \frac{\partial}{\partial r}$ be
an outward pointing Liouville flow. We choose a loop $l$
around $G$ so that the disc $V$ with $\partial V = l$
has volume  $v$ where $v$ can be chosen
arbitrarily small, and such that $X$ is
transverse to this loop.
This forms
a new cylindrical end $\varrho$ for $\C$.
Now let $H_\lambda^V$ be a Hamiltonian on $V$ with slope $\lambda$.
We assume that $H_\lambda^V$ has the following properties:
\begin{enumerate}
\item All the orbits are non-degenerate.
\item \label{item:actionbound}
The action of any orbit is in the region $[0,2v\lambda]$.
\item All orbits of index $\leq m + n$ are exact.
($2n$ is the dimension of our symplectic manifold)
\item \label{item:smallaction}
$2v\lambda < \nsmbox{min}(\delta/2,\delta_1)$
\item $H_\lambda^V$ is constructed in the same
way as the Hamiltonian $H_S$ mentioned before Lemma
\ref{lemma:lefschetzactionbounds}, but we do this with
respect to the new cylindrical end $\varrho$ of $\C$.
\end{enumerate}
%
We let $K_\lambda = \pi^*(H_\lambda^V) + \pi_1^* H_{F'}$.
Let $B := \C \times  (\widehat{F'} \setminus F'') \subset \widehat{E'}$.
The Hamiltonian $K_\lambda$
is of the form $\pi^*(H_\lambda^V) + \pi_1^* H_{F''}$ on $B$.
Hence the orbits on $B$ come in pairs
$(\gamma, \Gamma)$, where $\gamma$ corresponds to a periodic
orbit of $H_\lambda^V$ and $\Gamma$ is a constant periodic orbit of $H_{F''}$. 
The action difference between two orbits is 
$\leq \delta$ due to property (\ref{item:hfactionbound})
for $H_{F'}$ and properties (\ref{item:actionbound}) and (\ref{item:smallaction})
for $H_\lambda^V$.
This implies that any Floer trajectory connecting orbits $(\gamma_1,\Gamma_1)$
and $(\gamma_2,\Gamma_2)$ inside $B$ must stay
inside $B$, due to property (\ref{item:areaholomorphic})
from the properties of $J$. Also the definition of $\delta_1$
combined with properties (\ref{item:actionbound}) and (\ref{item:smallaction})
for $H_\lambda^V$ imply that all orbits inside $B$ have negative action. 

Hence we have a subcomplex $C^\lambda_{B}$ generated by orbits
in $B$. The Hamiltonian $K_\lambda$ in the region $\widehat{E} \setminus B$
is equal to $\pi^*(H_\lambda^V)$.
Hence, using Lemma \ref{lemma:lefschetzactionbounds},
we have that all the orbits of $\pi^*(H_\lambda^V)$
have non-negative action.
This means that we have a quotient complex
$C^\lambda_{E''} := C^\lambda_{E'} / C^\lambda_{B}$ 
where $C^\lambda_{E'}$ is
the complex generated by all orbits.
By property (\ref{item:smallslopeinmiddle}) of $H_{F'}$
and property (\ref{item:trajectoryconstraint}) of $J$
we have that the limit of $H_*(C^\lambda_{E''})$
as $\lambda$ tends to infinity is isomorphic to $SH^{\nsmbox{lef}}_*(E'')$.
So in order to show that $SH^{\nsmbox{lef}}_*(E'')$ is isomorphic
to $SH^{\nsmbox{lef}}_*(E')$ we need to show that the homology of $C^\lambda_B$
is zero in all degrees $\leq m$. This would imply that
$SH^{\nsmbox{lef}}_*(E'') = SH^{\nsmbox{lef}}_*(E')$ in all degrees $< m$. 
As $\lambda$ increases, we can make
$m$ increase with $\lambda$, and this means that
$SH^{\nsmbox{lef}}_*(E'') = SH^{\nsmbox{lef}}_*(E')$ in all degrees.
The point is that all the orbits in $C^\lambda_B$
have negative action, so any Floer trajectory
starting at one of these orbits must also
finish at one of these orbits. Also,
any Floer trajectory connecting two of these orbits
must be contained in $B$. We have that $B$ is a product.
We can assume that it has the product almost complex structure.
This means that Floer trajectories between $(\gamma_1,\Gamma_1)$
and $(\gamma_2,\Gamma_2)$ come in pairs $(u,U)$ where
$u$ is a Floer trajectory in $\C$ connecting $\gamma_1$
and $\gamma_2$ and $U$ is a Floer trajectory connecting
$\Gamma_1$ and $\Gamma_2$. Hence by a K\"unneth formula,
we have that the homology of $C^\lambda_B$
is isomorphic to $SH_*(H^V_\lambda)\otimes SH_*^{(-\infty,0)}(H_{F'})$,
where $SH_*^{(-\infty,0)}(H_{F'})$ is the subcomplex generated by
orbits negative action. We have that in all degress $\leq m$,
that $SH_*(H^V_\lambda) = 0$, and that $SH_*^{(-\infty,0)}(H_{F'})$
is $0$ in negative degrees which means that
$SH_*(H^V_\lambda)\otimes SH_*^{(-\infty,0)}(H_{F'})$ is
zero in all degrees $\leq m$. This in turn implies
that the homology of $C^\lambda_B$ is zero in all degrees
$\leq m$. This gives us our isomorphism:
\[SH^{\nsmbox{lef}}_*(E'') \cong SH^{\nsmbox{lef}}_*(E').\]
\qed

\section{Brieskorn Spheres} \label{section:brieskornspheres}

In this section we will mainly be studying the variety
$V$ as constructed in section \ref{section:exoticconstruction}
and also the variety $M'' := \C^4 \setminus V$.
The variety $V$ is equal to $\{z_0^7+z_1^2+z_2^2+z_3^2 = 0\} \subset \C^4$.

\subsection{Parallel transport} \label{subsection:brieskornparalleltransport}

There is a natural symplectic form on $\C^4$ (induced
from an ample line bundle on its compactification $\P^4$).
We have a holomorphic map $P := z_0^7+z_1^2+z_2^2+z_3^2$ with
one singular point at $0$. We can view $P$ as a fibration
which is compatible with this symplectic form
as in Definition \ref{defn:symplecticfibration}.
These fibrations have a natural connection which
is induced from the symplectically orthogonal plane
fields to the fibres.
We prove:
\begin{theorem} \label{theorem:paralleltransportP}
Parallel transport maps are well defined for $P$.
\end{theorem}
\proof
We first of all compactify $\C^4$ to $\P^4$.
We let $P'$ be a holomorphic section of $E:=\O_{\P^4}(7)$:
\[P'([z_0:\cdots:z_4]) := z_0^7 + z_4^5 z_1^2 + z_4^5 z_2^2 + z_4^5 z_3^2.\]
This is equal to $P$ on the trivialisation $z_4=1$.
We also have another section $Q$ defined by
\[Q([z_0:\cdots:z_4]) := z_4^7.\]
The map $P$ can be extended to a rational map
$P'' : \P^4 \dashrightarrow \P^1$, where
$P'' = \frac{P'}{Q}$.
Fix an identification $\C^4 = \P^4 \setminus Q^{-1}(0)$.
We now have that
\[P = \frac{P'}{Q}.\]
Let $\|.\|_E$ be a positive curvature metric on the ample bundle $E$.
We have a symplectic structure and K\"ahler form
defined in terms of the plurisubharmonic function
\[\phi = -\log\|Q\|_E^2.\]
In order to show that $P$ has well defined parallel transport maps 
we need to construct bounds on derivatives similar to the
main theorem in
\cite[section 2]{FukayaSeidelSmith:exactlagragiancotangent}.
We take the vector field $\partial_z$ on the base $\C$. It has
a unique lift with respect to the K\"ahler metric which is:
\[\xi := \frac{\nabla P}{\|\nabla P\|^2}.\]
Here $\|.\|$ is the K\"ahler metric and
$\nabla P$ is the gradient of $P$
with respect to this metric.
Take a point $p$ on $D := \{z_4=0\}$. We can assume without loss
of generality that this lies in the chart $\{z_1=1\}$.
In this chart we have that the metric
$\|.\|_E = e^\sigma |.|$ where $\sigma$ is a smooth function
and $|.|$ is the standard Euclidean metric with respect to this chart.
Then:
\[\phi = -\log\|Q\|^2 = -\log|Q|^2-\sigma.\]
The notation $\lesssim$ means that one term is less than
or equal to some constant times the other term.
Hence we get:
\begin{equation}
\begin{aligned}
B := \big| \xi.\phi \big|  \leq 
\frac{|\langle \nabla P,\nabla \sigma\rangle |}{||\nabla P||^2} + 
\frac{2 |Q| \cdot |\langle \nabla Q, \nabla P \rangle|}
{||\nabla P||^2 \cdot |Q|^2} \\
 \lesssim 
\frac{1}{||\nabla P||} + \frac{||\nabla Q||}{||\nabla P|| \cdot |Q|} 
\label{equation:parallelbound}
\end{aligned}
\end{equation}
We get similar equations to \ref{equation:parallelbound}
in the other charts $\{z_i=1\}$.
If we can show that for any compact set $T \subset \C$
the function $B$ is bounded above by
a constant $K$ in the region $T_1 := P^{-1}(T) \setminus A$ where
$A = \{|z_0|,|z_1|,|z_2|,|z_3| \leq 1\}$, then we have well 
defined parallel transport maps. 
This is because we get similar bounds if we lift other vectors
of unit length (i.e. $c\partial_z$ where $c \in U(1)$). Hence if we
have a path, then 
$\big| \xi.\phi \big|$
is bounded above by a constant on this path. This ensures that the
transport maps do not escape to infinity.
In the chart $\{z_3=1\}$,
we have that for $1 \leq i \leq 2$, 
\[\partial_i P = 2z_i / z_4^2,\]
\[\partial_0 P = 7z_0^6 / z_4^7.\]
We have the following bounds on derivatives: 
\[|\nabla Q| \lesssim \frac{|Q|}{|z_4|}.\]

Combining this with equation \ref{equation:parallelbound} gives:
\[B 
\lesssim \frac{1 + |z_4|}{|z_4| \left( \sum_{j=0,j \neq 3}^4 |\partial_i P| \right)}\]
\[\lesssim C := (1+|z_4|)/ \left[ \frac{7|z_0|^6}{|z_4|^6} + \frac{2|z_1|}{|z_4|} + \frac{2|z_2|}{|z_4|}
+ |\partial_4 P||z_4| \right]\]
\[ \lesssim (1+|z_4|) / \left( |z_0|^6/|z_4|^6 + |z_1|/|z_4| + |z_2|/|z_4| \right).\]
Hence on the chart $\{z_4=1\}$,
\begin{equation}
\begin{aligned}
B \lesssim (1+|z_3|^{-1}) / \left( |z_0|^6 + |z_1| +|z_2| \right)
\end{aligned} \label{equation:bound1}
\end{equation}
By symmetry we also have
\begin{equation}
\begin{aligned}
B \lesssim (1+|z_2|^{-1}) / \left( |z_0|^6 + |z_1| + |z_3| \right)
\end{aligned} \label{equation:bound2}
\end{equation}
\begin{equation}
\begin{aligned}
B \lesssim (1+|z_1|^{-1}) / \left( |z_0|^6 + |z_2| + |z_3| \right)
\end{aligned} \label{equation:bound3}
\end{equation}
In the chart $\{z_0=1\}$ we have:
\[B 
\lesssim \frac{1 + |z_4|}{|z_4| \left( \sum_{j=1}^4 |\partial_i P| \right)}\]
\[\lesssim (1+|z_4|)/ \left[ \frac{2|z_1|}{|z_4|} + \frac{2|z_1|}{|z_4|} + \frac{2|z_2|}{|z_4|}
+ |\partial_4 P||z_4| \right].\]
So:
\[B \lesssim (1+|z_4|)|z_4| / \left( |z_1| + |z_2|+ |z_3| \right)\]
so in the chart $\{z_4=1\}$, we get a bound:
\begin{equation}
\begin{aligned}
B \lesssim (1+|z_0|^{-1}) /  \left( |z_1| + |z_2|+ |z_3| \right)
\end{aligned} \label{equation:bound4}
\end{equation}
Suppose for a contradiction that there is a sequence of
vectors $(z_0^i,z_1^i,z_2^i,z_3^i)$ lying in $T_1$
such that $B$ tends to infinity as $i$ tends to infinity.
If (after passing to a subsequence) $z_0^i$ tends to infinity, 
then equation \ref{equation:bound4}
tells us that $z_1^i,z_2^i,z_3^i$ are all bounded. But 
this is impossible as $(z_0^i,z_1^i,z_2^i,z_3^i)$ lies in $T_1$
which means that $z_0^i$ is bounded. Similarly, using equations
\ref{equation:bound1},\ref{equation:bound2},\ref{equation:bound3}
we get that $z_j^i$ is bounded. Hence $B$ is bounded
away from the compact set $\{|z_0|,|z_1|,|z_2|,|z_3| \leq 1\}$.
This means that $B$ is bounded when restricted to $T_1$, so we have
well defined parallel transport maps.
\qed

\bigskip

Let $(\C^4,\theta)$ be the convex symplectic manifold induced by
the compactification $\C^4 \hookrightarrow \P^4$.
Because parallel transport maps for $P$ are well defined
we can use ideas from \cite[section 19b]{Seidel:fukayalefschetz}
to deform the $1$-form $\theta$ on $\C^4$ through
a series of $1$-forms $\theta_t$ such that: 
\begin{enumerate}
\item each $\omega_t := d\theta_t$ is compatible with $P$ 
as in Definition \ref{defn:symplecticfibration} and $\theta_t$ 
is a convex symplectic deformation on $\C^4$.
\item  \label{item:productboundary}
The parallel transport maps of $P$ with respect to the connection induced by $\omega_1$
are trivial at infinity.
This means that near infinity, $P$ looks like the natural projection
$C \times \C \twoheadrightarrow \C$ where $C$ is the complement of some compact set in $V$.
\item For a smooth fibre $F$ of $P$,
$(F,\theta_1)$ is exact symplectomorphic to $(F,\theta_0)$.
\end{enumerate}

We have that $M'' = \C^4 \setminus P^{-1}(0)$, so we can restrict $P$
to a fibration $P'' = P|_{M''} : M'' \rightarrow \C^*$. Let $\theta_S$
be a convex symplectic structure on $\C^*$ with the property that
$\theta_{M'',t} := \theta_t|_{M''} + {P''}^*\theta_S$ is a convex symplectic
structure for $M''$.
Let $\theta''$ be a convex
symplectic structure on $M''$ constructed as in Example
\ref{example:algebraicsteinmanifold}. It is convex deformation
equivalent to $(M'',\theta_{M'',0})$ as follows:
Let $F$ be a fibre of $P''$, then 
$(F,\theta_{M'',0}|_F)$ is convex
deformation equivalent to $(F,\theta''|_F)$
by Lemma \ref{lemma:algebraicuniqueness}
as both convex structures come from Stein structures
constructed algebraically as in Example \ref{example:algebraicsteinmanifold}. 
This deformation is $(1-t)\theta_{M'',0}|_F + t\theta''|_F$.
The following family of $1$-forms $\Theta_t := (1-t)\theta_{M'',0} + t\theta''$
induces a convex symplectic deformation (we might have to add $\pi^*\theta_S'$
to $\theta_{M'',0}$ and $\theta''$ where $\theta_S$ is a convex symplectic
structure on $\C^*$ and $d\theta_S$ is sufficiently large).
The reason why it is a convex symplectic deformation is as follows:
We can ensure that $\theta_S$ comes from a Stein function $\phi_S$ on $\C^*$.
Also, $\theta_0$ comes from some Stein function $\phi : \C^4 \rightarrow \R$,
hence $\theta_{M'',0}$ comes from a Stein function $\phi_0 := \phi|_{M''} + {P''}^*\phi_S$.
The $1$-form $\theta''$ comes from a Stein function $\phi_1$.
Hence $\Theta_t$ comes from a Stein function of the form
$\phi_t := (1-t)\phi_0 + t\phi_1$. The set of singular points
of $\phi_t|_F$ for all $t$ lie inside a compact set $K_F$
(independent of $t$) for each fibre $F$. 
Let $K$ be the union of all the compact sets $K_F$
for each fibre $F$ in $M''$. We can choose
$\theta_S$ large enough so that outside some annulus $A$ in $\C^*$,
$\phi_t$ has no singularities in $K \cap {P''}^{-1}(\C^* \setminus A)$.
Also, there are no singularities of $\phi_t$ outside $K$.
Hence, all the singularities of $\phi_t$ stay inside some
compact set independent of $t$. This means that $\phi_t$
is a Stein deformation. This means that $(M'',\theta_{M'',1})$ is
convex deformation equivalent to $(M'',\theta_{M'',0})$ which
is convex deformation equivalent to $(M'',\theta'')$.

Hence on $(M'',\theta'')$, we have that the parallel transport maps
of $P''$ are trivial at infinity after a convex symplectic deformation
to $(M'',\theta_{M'',1})$.

\subsection{Indices} \label{subsection:indices}

Let $P'' : M'' \rightarrow \C^*$, $P''(z) = P(z)$.
Let $F$ be a smooth fibre of $P''$. This fibre has a natural exhausting
plurisubharmonic function $\phi$
as in Example \ref{example:algebraicsteinmanifold}.
We can modify $\phi$ to an exhausting plurisubharmonic function $\phi'$
which is complete by \cite[Lemma 6]{SS:rama}. We denote this
new Stein manifold by $\widehat{F}$.
The following theorem
is about indices of a cofinal family of Hamiltonians on $\widehat{F}$.
\begin{theorem} \label{theorem:nontrivialboundary}
There is a cofinal family of Hamiltonians $H_\lambda$ on $\widehat{F}$ with
the following properties:
\begin{enumerate}
\item There exists some convex symplectic submanifold $T$ of $F$ such that
$\widehat{T}$ (the symplectic completion of $T$) is exact symplectomorphic
to $\widehat{F}$.
\item $H_\lambda = 0$ on $T$.
\item if $y$ is a periodic orbit of $H_\lambda$ not in $T$ then
 $\nsmbox{ind}(y) \geq 2$.

\item \label{item:finitenumberoforbits}
For each $k \in \Z$ there exists an $N > 0 $ (independent of $\lambda$) such that
the number of periodic orbits of $H_\lambda$ of index $k$ is bounded above by $N$.

\item \label{item:orbitsofindex2and3}
If we don't count critical points from the interior, then
there is exactly one orbit of index $2$ and one orbit of index $3$
such that the action difference between these two orbits tends to $0$ as $\lambda$ tends
to infinity.
Also the number of Floer cylinders connecting these orbits is even.
\end{enumerate}
\end{theorem}
This theorem is proved by analysing the Conley-Zehnder
indices of a Reeb foliation on the Brieskorn sphere
$V \cap S$, where $S$ is the unit sphere in $\C^4$.
This result needs the following two lemmas:
\begin{lemma} \label{lemma:brieskornboundary}
$\widehat{F}$ is the completion of some convex symplectic
submanifold $T$ with boundary the Brieskorn sphere
$V \cap S$.
\end{lemma}
\proof of Lemma \ref{lemma:brieskornboundary}.
By Theorem \ref{theorem:paralleltransportP}, we have that
$V \setminus 0$ is symplectomorphic to $F \setminus K$
where $K$ is a compact set.
Hence there exists a cylindrical end of $\widehat{F}$ which is symplectomorphic
to the cylindrical end of $V$ induced by flowing $V \cap S$ by parallel transport.
\qed

\begin{lemma} \label{lemma:brieskornconleycalculation}
There is a contact form on the Brieskorn sphere $V \cap S$ such that
all the Reeb orbits are non-degenerate and they have
Conley-Zehnder indices $\geq 2$. Also, there is exactly
one orbit of index $2$ and no orbits of index $3$ and finitely
many orbits of degree $k$ for each $k$.
\end{lemma}
\proof
In \cite{Ustilovsky:infinitecontact} Ustilovsky constructs
a contact form such that all the Reeb orbits are non-degenerate
and such that their reduced Conley-Zehnder index is
$\geq 2(n-2)$ where $n=3$ in our case.
Ustilovsky defines the reduced Conley-Zehnder
index to be equal to the Conley-Zehnder index $+ (n-3)$.
This means that the Reeb orbits have Conley-Zehnder index
$\geq n-1 = 2$. He also shows for each $k \in \Z$,
there are finitely many orbits of Conley-Zehnder index $k$.
He shows that there are no orbits of odd index and
the orbit of lowest index has index $2$.
\qed
\proof of Theorem \ref{theorem:nontrivialboundary}.

By Lemma \ref{lemma:brieskornboundary}
$\widehat{F}$ has a convex cylindrical end which is symplectomorphic
to $[1,\infty) \times \Sigma$ where $\Sigma$ is the Brieskorn
sphere $V \cap S$.
We choose a Hamiltonian which
is constant on the interior of $F$ and equal to
$h(r)$ on the cylindrical end, where $r$ parameterizes
$[1,\infty)$. We also assume that $h'(r)$ is
constant and not in the period spectrum of $B$ at infinity.
Also, near each orbit in the cylindrical end, we assume that $h'' > 0$.
The flow of the Hamiltonian at the
level $r=k$ is the same as the flow of
$X_H := -h'(k)R$, where $R$ is the Reeb flow.
The Conley-Zehnder indices from Lemma
\ref{lemma:brieskornconleycalculation} are computed by trivialising
the contact plane bundle. We can trivialise the symplectic
bundle by first trivialising the contact plane bundle
and then trivialising its orthogonal bundle.
We trivialise the orthogonal bundle by giving it a basis
$(\frac{\partial}{\partial r}, R)$. The symplectic form
restricted to this basis is the standard form:
\[\left( \begin{array}{cc}
0 & 1 \\
-1 & 0 \end{array} \right).\]
The Hamiltonian flow in this trivialisation is the matrix:
\[\left( \begin{array}{cc}
1 & 0 \\
h'' t & 1 \end{array} \right)\]
along the orthogonal bundle. This is because $R$ is invariant under
this flow and the Lie bracket of $X_H = -h'R$ with $\frac{\partial}{\partial r}$
is $h''R$.
The Robbin-Salamon index of this family of matrices is $\frac{1}{2}$.
We calculate this index by perturbing this family of matrices by a function
$\xi : [0,1] \rightarrow \R$ where $\xi(0) = \xi(1) = 0$ as follows:
\[\left( \begin{array}{cc}
1 & \xi(t) \\
h'' t & 1 \end{array} \right).\]
Choosing $\xi$ so that its derivative is non-zero whenever $\xi = 0$
ensures that the path is generic enough to
enable us to compute its Robbin-Salamon index.

Remember that the Robbin-Salamon index of an orbit is equal to the Conley-Zehnder
index taken with negative sign.
Lemma \ref{lemma:brieskornconleycalculation} tells us the Conley-Zehnder
indices of all the Reeb orbits.
The flow $X_H = -h'R$ of the Hamiltonian
has orbits in the opposite direction to Reeb orbits. Hence the
Robbin-Salamon index (restricted to the contact plane field)
of an orbit of $X_H = -h'R$ is the
same as the Conley-Zehnder index of the corresponding Reeb orbit.
Hence the Robbin-Salamon index of some orbit of the
Hamiltonian on the level set $r=k$ is equal
to $C+\frac{1}{2}$ where $C$ is the Conley-Zehnder index
of the associated Reeb orbit as calculated in
Lemma \ref{lemma:brieskornconleycalculation}.
Hence the indices of these orbits are
$\geq 2+\frac{1}{2}$.

The problem is that these orbits are degenerate. This is why
their index is not an integer.
As in \cite[section 3]{Oancea:survey} we can
perturb each circle of orbits to a pair of non-degenerate orbits.
Let $C'$ be a circle of orbits.
We choose a Morse function $f$ on $C'$. If we
flow $f$ along $X_H$ (the Hamiltonian flow of $H$)
we get a time dependent Morse function $f_t = f \circ \phi_{-t}$
($\phi_t$ is the Hamiltonian flow).
Extend $f_t$ so that it is defined as a function on a
neighbourhood of $C'$.
Let $H+f_t$ be our new Hamiltonian.
The orbits near $C'$ now correspond to critical
points $p$ of $f$.
The Robbin-Salamon index of such an orbit is:
\[i(C') + \frac{1}{2}\nsmbox{sign}(\nabla_p^2 f)\]
where $i(C')$ is the Robbin-Salamon index of the manifold
of orbits. The symbol `$\nsmbox{sign}$' means the number
of positive eigenvalues minus the number of negative eigenvalues.
In our case we can choose $f$ so that it
has $2$ critical points $p_1,p_2$ such that
\[\nsmbox{sign}(\nabla^2_{p_1} f)=1,\nsmbox{ sign}(\nabla^2_{p_2} f)=-1.\]
Hence, if the Conley-Zehnder index of a Reeb orbit $C$ is
$k$, then we can perturb $H$ so that the associated
Hamiltonian orbits have Robbin-Salamon index
(or equivalently Conley-Zehnder index taken with negative sign)
$k+0$ and $k+1$.
This means all the non-constant orbits of $H$ have Robbin-Salamon
index $\geq 2$.

We now need to show that there are a finite number of orbits in each degree.
This follows directly from Lemma \ref{lemma:brieskornconleycalculation} which says
that there are finitely many Reeb orbits in each degree.
Finally this same lemma says that there is only one Reeb orbit with Conley-Zehnder index $2$
and no Reeb orbits with Conley-Zehnder index $3$.
So the Hamiltonian $H$ has one orbit of Robbin-Salamon index $2$
and one orbit of index $3$. We can also ensure that the actions
of these orbits are arbitrarily close by letting the associated
Morse function $f$ be $C^2$ small.
There are an even number of Floer cylinders connecting
the orbit of index $3$ with the orbit of index $2$
by \cite[Proposition 2.2]{CieliebakFloerHoferWysocki:SymhomIIApplications}.
\qed

\begin{lemma} \label{lemma:homologyofbrieskorn}
We have $H^i(M'') = 0$ for $i \geq 2$.
\end{lemma}
\proof
$M'' = \C^4 \setminus V$. Theorem \ref{theorem:brieskorntopology}
tells us that $V$ is homeomorphic to $\R^6$.
This means that there is a neighbourhood $B$
of $V$ which retracts onto $V$ whose boundary $\partial B$
satisfies $H^i(\partial B) = 0$ for $i \geq 2$.
The Mayor-Vietoris sequence involving $B$, $M''$ and $B \cup M'' = \C^4$
ensures that $H^i(M'') = 0$ for $i \geq 2$.
\qed

\subsection{Symplectic homology of these varieties} \label{section:symplectichomologypremodified}

We wish to show that the symplectic homology of the variety
$M'' := \C^4 \setminus V$ has only finitely many idempotents using
the results of the previous two sections. We will then show
that it has at least two idempotents: $0$ and $1$.
First of all we
need the following lemma:
We let $R$ be an algebra over $\Z / 2 \Z$ which is graded by a finitely generated abelian
group $G$. This means that as a vector space, $R = \bigoplus_{g \in G} R_g$ with the property that
if $a \in R_{g_1}$ and $b \in R_{g_2}$ then the product
$ab$ is contained in $R_{g_1.g_2}$.

\begin{lemma} \label{lemma:finitelymanyidempotents}
If $a$ is an idempotent in $R$ then
$a \in \bigoplus_{g \in G_n} R_g$
where $G_n$ is the subgroup of torsion elements of $G$.
\end{lemma}
\proof
We have $a = a_{g_1} + \dots + a_{g_n}$ where $g_i \in G$ and $a_{g_i} \in R_{g_i}$.
Suppose for a contradiction we have that $a = a^2$ and $g_1$ is not torsion.
Then $a^2 = a_{g_1}^2 + \dots + a_{g_n}^2$.
The group $G/G_n$ is a free $\Z$ algebra, hence there is a group
homomorphism $p : G \rightarrow G/G_n \rightarrow \Z$ such that $p(g_1) \neq 0$.
The map $p$ gives $R$ a $\Z$ grading. Let $b$ be an element of $R$.
It can be written uniquely as $b = b_1 + \dots + b_k$ where $b_i$
are non-zero elements of $R$ with grading $q_i \in \Z$. We can define $f(b)$ as
$\nsmbox{min}\{|q_j| \neq 0 \}$.
Note that $f(b)$ is well
defined only if at least one of the $q_i$'s are non-zero.
Because $p(g_1) \neq 0$, we have that $f(a)$ is well defined and positive.
We also have that $f(a^2) \geq 2f(a)$ which means that $a^2 \neq a$.
This contradicts the fact that $a$ is an idempotent.
\qed

\bigskip

The vector space $SH_{4+*}(M'')$
 is a ring bi-graded by the Robbin-Salamon index 
and the first homology group. 
We write $4+*$ here because the unit has Robbin-Salamon index $4$.
The previous lemma shows us that
any idempotent must have grading $4$ in $SH_*(M'')$ and be in a torsion
homology class.
%

We have a map $P'' : M'' \rightarrow \C^*$.
At the end of section \ref{subsection:brieskornparalleltransport} we had a
convex symplectic structure
$(M'',\theta_{M'',1})$.
Let $A$ be a large annulus in the base $\C^*$ which is a compact convex
symplectic manifold.
Let $(F'',\theta_{M'',1})$ be a fibre of $P''$. Choose a
compact convex symplectic manifold (with corners)
$\bar{M}''$ such that $(\bar{M}'',\bar{P}'' := P''|_{\bar{M}''},\theta_{M'',1})$
is a compact convex Lefschetz fibration with fibres $\bar{F}''$
and base $A \subset \C^*$. We can also ensure that $\partial \bar{F}''$ is
transverse to $\lambda_1$ (the associated Liouville vector field of $\F''$)
and there are no singularities of $\lambda_1$
outside $\bar{F}''$ in $F''$. Hence the completion of $\bar{M}''$
is $\widehat{M}''$.

Let $(\widehat{E}'',\pi'')$ be the completion of $(\bar{M}'',\bar{P}'',\theta_{M'',1})$
(so that $\widehat{\bar{M}''} = \widehat{E''}$).
We wish to use the results of section \ref{section:lefschetzcofinal} to
show that $SH_*(E'')$ has finitely many idempotents, and hence $SH_*(M'')$
has finitely many idempotents.
Let $H$ be a Lefschetz admissible Hamiltonian for $\widehat{E''}$.
Let $C$ be the cylindrical end of $\widehat{F''}$.
We may assume that this cylindrical end is of
the form $(S_V \times [1,\infty),r_F \alpha_F)$ where
$(S_V,\alpha_F)$ is the Brieskorn sphere described in
\ref{subsection:indices} and $r_F$ is the coordinate for $[1,\infty)$.
The Hamiltonian
$H$ is of the form ${\pi''}^* H_{S''} + {\pi_1''}^* H_{F''}$ as in
Definition \ref{defn:lefschetzadmissiblehamiltonian}.
By Lemma \ref{lemma:finitelymanyidempotents}, we have that
any idempotent must come from a linear combination of
orbits of $H$ in torsion homology classes
as long as $H$ is large enough (i.e. it is large enough in
some cofinal sequence of Lefschetz admissible Hamiltonians).
Away from $C \times \C^* \subset \widehat{E''}$ we have that the Hamiltonian
flow of $H$ is the same as the flow of $L := {\pi''}^*H_{S''}$. Let $X$ be the
Hamiltonian vector field associated to $L$, and let $X_{S''}$ be the
Hamiltonian vector field in $\C^*$ associated to $H_{S''}$.
Then the value of $X$ at a point $p$ is some positive multiple of the horizontal lift
of $X_{S''}$ to the point $p$.
We can assume that $H_{S''}$ has exactly two contractible
periodic orbits of index $0$ and $1$ corresponding to
Morse critical points of $H_{S''}$ (as any Reeb orbit of $\C^*$
is not contractible). We can also make $H_{S''}$ $C^2$ small
away from the cylindrical ends of $\C^*$ so that the only
Floer cylinders connecting contractible orbits correspond
to Morse flow lines.
Hence, any contractible orbit of $X$ must project down to a constant
orbit of $X_{S''}$.
We let $H_{F''}$ be a Hamiltonian as in Theorem \ref{theorem:nontrivialboundary} above in Section
\ref{subsection:indices}. We let our almost complex structure $J$
when restricted to  $C \times \C^* \subset \widehat{E''}$
be equal to the product almost complex structure $J_F \times J_{\C^*}$
where $J_F$ is an admissible almost complex structure on $\widehat{F}$
and $J_{\C^*}$ is the standard complex structure on $\C^*$.
The contractible orbits in this cylindrical end come in pairs $(\Gamma,\gamma)$
where $\Gamma$ is an orbit in $\widehat{F}$ and $\gamma$
is a contractible orbit in $\C^*$.
Because there are only $2$ contractible orbits in $\C^*$ and there are finitely
many orbits in each degree in $\widehat{F}$, we have
finitely many contractible orbits of index $4$ for $H$.
Hence:
\begin{theorem} \label{theorem:finitelymanyidempotentsofMprimeprime}
The ring $SH_{4+*}(M'')$ has only finitely many idempotents.
\end{theorem}

\bigskip

We now wish to show that $SH_*(M'')$ has at least $2$ idempotents.
To do this we show that $SH_*(M'') \neq 0$, and hence has a unit
by \cite[Section 8]{Seidel:biasedview}. This means that
$SH_*(M'')$ has $0$ and $1$ as idempotents.
The Hamiltonian $H$ has non-degenerate orbits in
$C \times \C^* \subset \widehat{E''}$, so we perturb
$H$ away from this set to make all its orbits non-degenerate.
In $E'' \subset \widehat{E''}$ we can ensure that $H$ is $C^2$
small and $J$ is independent of $t$,
hence the only orbits in this region are critical points of $H$ and
the only Floer cylinders correspond to Morse flow lines.
The orbits corresponding to critical points of $H$
have Robbin-Salamon index $\geq 3$ because $H^i(M'') = 0$ for $i > 1$
by Lemma \ref{lemma:homologyofbrieskorn}.
Hence all orbits have index $\geq 2$.
There is only one orbit of index $2$. 
This orbit is
in the cylindrical end $C \times \C^* \subset \widehat{E''}$.
Hence the orbit is of the form $(\Gamma_m,\gamma_m)$
where $\gamma_m$ has index $0$ and $\Gamma_m$ has index $2$.
This orbit is closed because there are no orbits of lower index.
Suppose for a contradiction this orbit is exact, then there exists a Floer
cylinder connecting an orbit $\beta$ of index $3$ with
 $(\Gamma_m,\gamma_m)$. This orbit $\beta$ must
be contractible, so it is either a critical point,
or it is of the form $(\Gamma_1,\gamma_1)$ in  $C \times \C^* \subset \widehat{E''}$.
The action of $(\Gamma_m,\gamma_m)$ is larger than the action
of a critical point and hence $\beta$ cannot be a critical point.
Hence $\beta$ is of the form $(\Gamma_1,\gamma_1)$.
Suppose that the index of $\Gamma_1$ is $3$.
We have $\gamma_m = \gamma_1$ and by Theorem \ref{theorem:nontrivialboundary}
we can ensure that the action difference between
$\Gamma_1$ is arbitrarily close to $\Gamma_m$.
Similarly if $\gamma_1$ has index $1$ then we can ensure
that $\Gamma_m = \Gamma_1$ and the action difference
between $\gamma_m$ and $\gamma_1$ is arbitrarily small.
This means that the action difference between
$(\Gamma_m,\gamma_m)$
and $(\Gamma_1,\gamma_1)$ is arbitrarily small.
This means that if we have a Floer cylinder connecting $(\Gamma_m,\gamma_m)$
and $(\Gamma_1,\gamma_1)$ then Gromov compactness ensures that
it must stay in the region  $C \times \C^* \subset \widehat{E''}$
(Because the action of $(\Gamma_m,\gamma_m)$ tends to
the action of $(\Gamma_1,\gamma_1)$,
we get a sequence of Floer cylinders converging to a Floer
cylinder of energy $0$ which cannot exit
 $C \times \C^* \subset \widehat{E''}$).
Because all the Floer cylinders stay inside
 $C \times \C^* \subset \widehat{E''}$,
the number of Floer cylinders connecting
 $(\Gamma_m,\gamma_m)$ and  $(\Gamma_1,\gamma_1)$ is
equal to the number of Floer cylinders
connecting $\Gamma_m$ and $\Gamma_1$ multiplied by the
number of Floer cylinders connecting
$\gamma_m$ and $\gamma_1$.
We need to show that the number of Floer cylinders connecting
$(\Gamma_1,\gamma_1)$ and  $(\Gamma_m,\gamma_m)$ is even and
by the previous comment, this means we only need to show that
the number of Floer cylinders connecting $\Gamma_1$ and $\Gamma_m$
is even or the number of Floer cylinders connecting
$\gamma_1$ and $\gamma_m$ is even.
But the number of Floer cylinders connecting $\Gamma_1$
and $\Gamma_m$ is even if $\Gamma_1$ has index $3$ 
(by part (\ref{item:orbitsofindex2and3}) of Theorem \ref{theorem:nontrivialboundary})
and similarly $\gamma_1$ is closed if it has index $1$
(so there are an even number of Floer cylinders connecting $\gamma_1$
and $\gamma_m$).
Hence the number of Floer cylinders connecting
these two orbits is even and so
$(\Gamma_m,\gamma_m)$ is not exact.
Hence $SH_*(M'') \neq 0$.

This completes the proof of the main theorem
\ref{thm:infinitelymanysteinmanifolds} subject to checking ring
addition under end connect sums.

\section{Appendix A: Lefschetz fibrations and the Kaliman modification} \label{section:appendix}

Let $X$,$D$,$M$ be as in Example \ref{example:algebraicsteinmanifold}.
This means that $X$ is a projective variety with $D$ an effective ample divisor
and $M = X \setminus D$ an affine variety.
Let $Z$ be an irreducible
divisor in $X$ and $q \in (Z \cap M)$ a point in
the smooth part of $Z$. We assume there is a rational function $m$
on $X$ which is holomorphic on $M$ such that 
$\overline{m^{-1}(0)}$ is reduced and irreducible and
$Z = \overline{m^{-1}(0)}$.
Let $M' := \nsmbox{Kalmod}(M,(Z \cap M),\{q\})$,
and let $M'' := M \setminus Z$.
Suppose also that $\nsmbox{dim}_{\C}X \geq 3$.
{\it Here is the statement of Theorem
\ref{theorem:KalimanLefschetz}:
There exist Lefschetz fibrations $E'' \subset E'$
respectively satisfying the conditions of Theorem
\ref{theorem:fibrationaddition} such that
$E'$ (resp. $E''$) is convex deformation equivalent to $M'$ (resp. $M''$).
}

The rest of this section is used to prove this theorem.
We will start with several preliminary lemmas.
\begin{lemma}  \label{lemma:symplecticsubmanifold}
There are Stein functions $\phi'$ (resp. $\phi''$)
on $M'$ (resp. $M''$) such that $M''$ becomes a symplectic
submanifold of $M'$.
\end{lemma}
\proof
Let $m'$ be the pullback of $m$ to $\nsmbox{Bl}_q X$.
Let $Z'$ be the divisor defined by the zero set of $m'$.
Let $Z''$ be the divisor defined by the zero set of $\frac{1}{m'}$,
so that $Z'$ is linearly equivalent to $Z''$.
By abuse of notation, we write $D$ as the total transform of $D$
in $\nsmbox{Bl}_q X$.

Let $\tilde{Z}$ be the proper transform of $Z$.
We can choose an effective ample divisor
$D'$ with support equal to $\tilde{Z} \cup D$ (as a set)
so that $D' -  Z''$ is effective.
We have that $Y_1 := D'$ and $Y_2 := Y_1 - Z'' + Z'$ are linearly
equivalent effective ample divisors. Let $E$ be a line bundle associated
to $Y_1$ and let $s_1,s_2$ be sections so that $s_i^{-1}(0) = Y_i$.
There is a metric $\|.\|$ of positive curvature on $E$.
We define $\phi' := -dd^c \log(s_1)$ and $\phi'' :=  -dd^c \log(s_2)$.
\qed

\bigskip

Moving the point $q$ within the smooth part of $Z \cap M$
induces a Stein deformation of $M'$ and $M''$ by a slight modification
of the above lemma.

We now need a technical lemma involving convex symplectic manifolds of
finite type.
Let
$(M,\theta_1)$, $(M,\theta_2)$ be convex symplectic manifolds.
Suppose that
$\theta_1=\theta_2$ inside some codimension $0$ submanifold $C$
such that $(C,\theta_1)$ is a compact convex symplectic manifold.

\begin{lemma} \label{lemma:convexsymplecticdeformations}
If all the singular points of $\theta_1$ and $\theta_2$
are contained in $C$, then $(M,\theta_1)$ is convex
deformation equivalent to $(M,\theta_2)$.
\end{lemma}
\proof
The interior $C^o$ of $C$ has the structure of a 
finite type non-complete convex symplectic manifold constructed
as follows:
The boundary of $C$ has a collar neighbourhood in $C$
of the form $N := (-\epsilon,1] \times \partial C$, with
$\theta_1 = r\alpha$. Here $r$ is the coordinate on $(-\epsilon,1]$,
and $\alpha$ is a contact form on $\partial C$. We let 
$\psi : C^o \rightarrow \R$ be an exhausting function,
which is of the form $h(r)$ on $N$ and such that $h(r) \rightarrow \infty$
as $r \rightarrow 1$. For some $N \gg 0$, we have
that $\psi^{-1}(l)$ is transverse to the associated
Liouville field $\lambda_1$ for all $l \geq N$.
Let $\phi_1$ be the function associated to the convex symplectic
structure $(M,\theta_1)$. We may assume that $\phi_1^{-1}(l)$
is transverse to $\lambda_1$ for all $l \geq N$ as well.
We can smoothly deform the function $\phi_1$ into the function
$\psi$ through a series of exhausting functions $\phi_t$ 
(the domain of $\phi_t$ smoothly changes within $M$ as $t$ varies)
such that
$\phi_t^{-1}(N +k)$ is transverse to $\lambda_1$ for each $k \in \N$.
This induces a convex symplectic deformation from $(M,\theta_1)$
to $(C^o,\theta_1|_{C^o})$. Similarly we have a convex symplectic
deformation from $(M,\theta_2)$ to  $(C^o,\theta_1|_{C^o})$.
Hence,  $(M,\theta_1)$ is convex
deformation equivalent to $(M,\theta_2)$.
\qed

We need another similar lemma about deformation equivalence.

\begin{lemma} \label{lemma:deformationequivalenceexactness}
Suppose $(M,\theta_1)$ and $(M,\theta_2)$ are convex
symplectic manifolds such that $\theta_1 = \theta_2 + dR$
for some function $R$, then $(M,\theta_1)$ is convex
deformation equivalent to $(M,\theta_2)$.
\end{lemma}
\proof
Let $\phi_1$ (resp. $\phi_2$) be the function associated
with the convex symplectic structure $(M,\theta_1)$ (resp. $(M,\theta_2)$).
Choose constants $c_1 < c_2 < \cdots$ and
$d_1 < d_2 < \cdots$ tending to infinity
such that $M_i^1 := \phi_1^{-1}(-\infty,c_i]$
(resp. $M_i^2 := \phi_2^{-1}(-\infty,c_i]$) are compact convex
symplectic manifolds. Also we assume that:
\[M_i^1 \subset M_i^2 \subset M_{i+1}^1 \subset M_{i+1}^2\]
for all $i$.
Let $R' : M \rightarrow \R$ be a function such that
$R' = 0$ on a neighbourhood of $\partial M_i^1$ and
$R' = R$ on a neighbourhood of $\partial M_i^2$ for all $i$.
Let $\theta_3 := \theta_1 + dR'$. We will show that
both $(M,\theta_1)$ and $(M,\theta_2)$ are convex deformation
equivalent to $(M,\theta_3)$.
Let $R_t : M \rightarrow \R$ be a family of functions such that
$R_t = 0$ on a neighbourhood of $\partial M_i^1$ for all $i$
and such that $R_0 = 0$ and $R_1 = R'$.
Then $(M,\theta_1 + dR_t)$ is a convex deformation from
$(M,\theta_1)$ to $(M,\theta_3)$ because
$(M_i^1,\theta_1 + dR_t)$ is a compact convex symplectic manifold for all $i$.
Also let $R_t' : M \rightarrow \R$ be a family of functions such
that $R_t' = R$ on a neighbourhood of $\partial M_i^2$ and such that
$R_0' = R$ and $R_1' = R'$. Then $(M,\theta_1 + dR_t')$ is a
convex deformation from $(M,\theta_2)$ to $(M,\theta_3)$.
Hence $(M,\theta_1)$ is convex deformation equivalent to $(M,\theta_2)$.
\qed

\bigskip

We let $E$ be an ample line bundle on $X$, and $s,t$ sections
of $E$. We assume that $s$ is non-zero on $M$.
Let $t$ be a holomorphic section of $E$, and let $p := t/s$ be a map
from $M$ to $\C$.

\begin{defn} \label{defintion:algebraiclefschetzfibration}

We call $(M,p)$ an {\bf algebraic
Lefschetz fibration} if:

\begin{enumerate}

\item $\overline{t^{-1}(0)}$ is smooth, reduced
and intersects each
stratum of $D$ transversally.

\item $p$ has only nondegenerate critical points and there is at most
one of these points on each fibre.

\end{enumerate}

\end{defn}

An algebraic Lefschetz fibration $(M,p)$ has
a symplectic form $\omega$ constructed as
in Example \ref{example:algebraicsteinmanifold}.
This means that $\omega$ is compatible
with $p$. These are not exact Lefschetz fibrations
since the horizontal boundary is not trivial, but they are
very useful since our examples arise in this way.
\begin{theorem} \label{theorem:algebraicparalleltransport}
Parallel transport maps for an algebraic Lefschetz fibration are well defined.
\end{theorem}
This is basically proved in
\cite[section 2]{FukayaSeidelSmith:exactlagragiancotangent},
but, there is a subtle distinction between the above theorem
and theirs.
In \cite[section 2]{FukayaSeidelSmith:exactlagragiancotangent},
the Stein structure and the Lefschetz fibration are constructed
from the same compactification $(X,D)$ of $M$. In our case
they come from different compactifications. The proof
can be easily adjusted to this case. This is due to the fact
that if we have two different metrics on $M$ induced from
compactifications $(X_1,D_1)$ and $(X_2,D_2)$, then
the $C^2$ distance between them is bounded.

We need the following technical lemma so that we can relate
algebraic Lefschetz fibrations with ordinary Lefschetz fibrations.
We let $(E',\pi')$, $(E'',\pi'')$ be algebraic Lefschetz fibrations
such that $\pi''|_{E''} = \pi'$. Let $\theta'$ (resp. $\theta''$) be a convex
symplectic structure 
on $E'$ (resp. $E''$) constructed as in Example \ref{example:algebraicsteinmanifold}
such that $d\theta'' = d\theta'|_{E''}$.
We assume that the real dimension of $E'$ and $E''$ is $4$ or higher.
\begin{lemma} \label{lm:exacttoalgebraiclf}
Suppose that all the singular points of $\pi'$ are contained
in $E''$. Then there exists a convex symplectic structure $\theta'_1$ (resp. $\theta''_1$)
on $E'$ (resp. $E'')$ such that:
\begin{enumerate}
\item $(E',\pi',\theta'_1)$ (resp. $(E'',\pi'',\theta''_1)$) are Lefschetz
fibrations without  \\ boundary.
\item ${d\theta'_1}|_{E''} = d\theta''_1$.
\item  All the parallel transport maps are trivial on a neighbourhood $N$
of $E' \setminus E''$, and
$E' \setminus N$ is relatively compact
when restricted to each fibre.
\item For each smooth fibre $F'$ of $\pi'$,
$\theta'|_{F'} = \theta'_1|_{F'} + dR$ for some
compactly supported function $R$.
We have a similar statement for $(E'',\pi'')$.
\item $(E',\theta'_1)$ (resp. $(E'',\theta''_1)$) is convex symplectic
deformation equivalent to  $(E',\theta')$ (resp. $(E'',\theta'')$).
\end{enumerate}
\end{lemma}
%
\begin{figure}[b]
\end{figure}
\begin{figure}[b]
\end{figure}
\begin{figure}[b]
\end{figure}
\proof of Lemma \ref{lm:exacttoalgebraiclf}
We divide this proof into 3 sections.
In the first section we construct the Lefschetz fibration without boundary
 $(E',\pi',\theta'_1)$.
In the second section we construct $(E'',\pi'',\theta''_1)$.
In section 3 we show that $(E',\theta'_1)$ (resp.  $(E'',\theta''_1)$)
 is convex deformation equivalent to $(E',\theta')$ (resp. $(E'',\theta'')$). 

{\bf Step 1}
We will use ideas from \cite[section 19b]{Seidel:fukayalefschetz}.
The map $\pi'$ has well defined parallel transport maps by
Theorem \ref{theorem:algebraicparalleltransport}.
We have the same for $(E'', \pi'')$.
Suppose without loss of generality that $0 \in \C$ is a regular point of these
fibrations. Let $Q' := {\pi'}^{-1}(0)$, 
$Q'' := Q \cap E''$.
Consider the family of radial lines in $\C$ coming out of $0$.
Let $L$ be one of these radial lines which passes through
a critical value $l$ of $\pi'$. We can write
$L = L_1 \cup L_2$ where $L_1$ is the line joining
$0$ and $l$, and $L_1 \cap L_2 = \{l\}$.
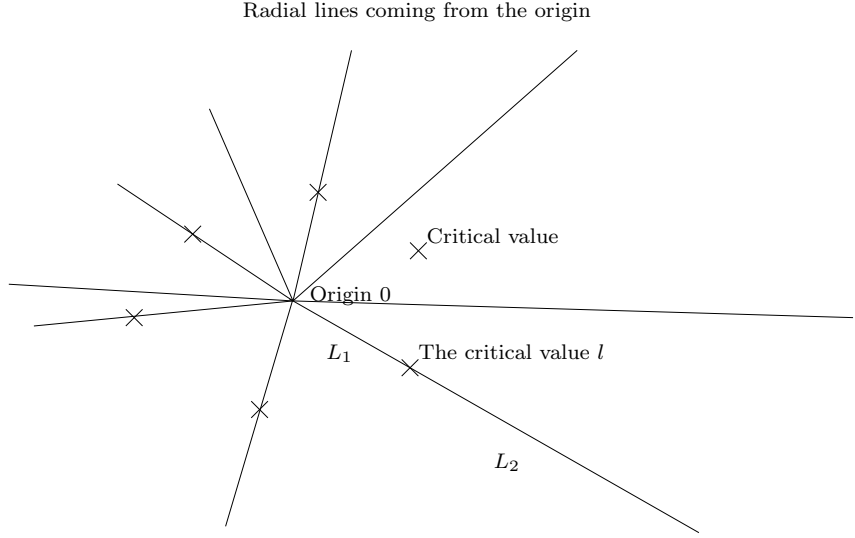
\begin{figure}[H]
 \centerline{
   \scalebox{1.0}{
     \input{radiallines.pstex_t}
   }
 } 
\label{figure:radiallines}
\end{figure}
We now have vanishing thimbles $V_1$ and $V_2$
of $l$ covering $L_1$ and $L_2$.
A vanishing thimble covering a line $L_i$
is just the set of points in ${\pi'}^{-1}(L_i)$
which parallel transport along $L_i$ into the critical point associated
to the critical value $l$.
Let $V$ be the union of all such thimbles for
all radial lines passing through critical values
of $\pi'$.
We can use this to construct
a map $\rho : E' \setminus V \rightarrow \C \times Q'$. The map is
constructed as follows:
Let $x$ be a point in $E' \setminus V$. Then we can parallel
transport $x$ along a radial line
to a point $a$ in $Q'$. Then $\rho(x) := (\pi'(x),a)$.
Let $X = V \cap Q'$ and $W := Q' \setminus X$. Then
$\rho|_{E' \setminus V} : E' \setminus V \rightarrow \C \times W$
is a diffeomorphism. Let $\varpi := (\rho|_{E' \setminus V})^{-1}$.
Let $\theta_{Q'} := \theta'|_{Q'}$.
From now on, if we have a differential form $q$ on $\C \times W$,
then we will just write $q$ instead of $\rho^*q$ to clean up notation.

Because parallel transport maps are exact, we have:
$\theta'|_{E' \setminus V} = \theta_{Q'} + \kappa' + dR'$ where
$\kappa'$ is a 1-form satisfying $i^*\kappa' = 0$ for all maps $i$ where
$i$ is the inclusion map of any fibre of $\pi'$ into $E'$,
and $R$ is some function on $\C \times W$.
Let $\bar{f} : W \rightarrow \R$
be a function which is equal to $1$ near $X$ and is $0$
outside some relatively compact neighbourhood of $X$. We extend $\bar{f}$
by parallel transport along these radial lines to a map
$g : E' \setminus V \rightarrow \R$. Then we extend
$g$ to a map $f : E' \rightarrow \R$ as $g$ is
constant near $V$.
We will also assume that $g$ is only non-zero inside $E''$
because parallel transport maps are well defined for
$(E'',\pi'',\theta'')$, hence $V \subset E''$.
We define
\[\theta'_f := \theta_{Q'} + g\kappa' + d(gR').\]
This form extends over $V$ because $\theta'_f = \theta'$ near $V$
(where $g=1$).
The $1$-form $\theta'_f$ makes $\pi'$ into a Lefschetz fibration without
boundary where each of the fibres have a convex symplectic structure.
We define $\theta'_1 := \theta'_f$.

{\bf Step 2}
Let $Q'' := {\pi''}^{-1}(0) \subset E''$.
We also have that:
\[\theta'' = \theta_{Q''} + \kappa'' + dR''.\]
Here, $\kappa''$ is  a 1-form on $E''$ satisfying $i^*\kappa'' = 0$
for all maps $i$ where $i$ 
is the inclusion map of a fibre of $\pi''$ into $E''$,
and $R''$ is some function on $E'' \cap (\C \times W)$.
Because $d\theta'' = d\theta'$, we have that $d\kappa' = d\kappa''$.
This means that $\beta := \kappa' - \kappa''$ is a closed $1$-form in $E''$.
We can also show that $\beta$ is exact as follows:
Let $l : {\mathbb S}^1 \rightarrow E''$ be a loop. 
Because we are in dimension $4$
or higher, 
we can perturb the loop so that it doesn't intersect
the radial vanishing thimbles described above.
We can then deform $l$ using parallel transport
to a loop $l'$ contained in a smooth fibre $F$. We have $\beta|_F=0$
which means that $\int_l \beta = \int_{l'} \beta = 0$. Hence $\beta = dL$
for some $L : E'' \rightarrow \R$.
We define:
\[\theta''_f := \theta_{Q''} + g\kappa' + d(gL) + d(gR').\]
We have $d\theta''_f = d\theta'_f$, hence this makes $(E'',\pi'')$
into a well defined symplectic subfibration of $E'$.
We define $\theta''_1 := \theta''_f$.

{\bf Step 3}
We can deform $\bar{f}$ through functions which are trivial at infinity
to some $\bar{f}'$ where 
$\bar{f}' = 0$ outside some large compact set, and
$(\bar{f}')^{-1}(1)$ contains
an arbitrarily large compact set $K \subset F$. We can construct
$f' : E' \rightarrow \R$ using $\bar{f}'$ in the same way that
we constructed $f$ from $\bar{f}$ and the deformation from
$\bar{f}$ to $\bar{f}'$ induces a deformation from $f$ to $f'$.
We can choose a convex symplectic structure $\theta_S$
on the base so that $(E',\theta'_f + {\pi'}^*\theta_S)$
and $(E',\theta'_{f'} + {\pi'}^*\theta_S)$ are convex
symplectic manifolds.
Hence $\theta'_f + {\pi'}^*\theta_S$ is convex deformation equivalent to
$\theta'_{f'} + {\pi'}^*\theta_S$.
If we choose $K$ large enough we get that
$\theta'_{f'} + {\pi'}^*\theta_S$ is convex deformation equivalent to
$(E',\theta'+{\pi'}^*\theta_S)$ by Lemma
\ref{lemma:convexsymplecticdeformations} and
Lemma \ref{lemma:deformationequivalenceexactness},
and hence is convex deformation equivalent
to $(E',\theta')$.

Because $\theta''_f$ is described in a very similar way to $\theta'_f$,
we can use exactly the same argument as above to show that
$(E'',\theta'')$ is convex deformation equivalent to
$(E'',\theta''_f + \pi^*\theta_{S,1})$. The $1$-form
$\theta_{S,1}$ is a convex symplectic structure on the base
making $\theta''_f + \pi^*\theta_{S,1}$ into a convex symplectic structure.

\qed

Let $X$,$D$,$M$ be as in Theorem \ref{theorem:KalimanLefschetz}.
This means that $Z$ is an irreducible
divisor in $X$ and $q \in (Z \cap M)$ is a point in
the smooth part of $Z$. There is a rational function $m$
on $X$ which is holomorphic on $M$ such that
$\overline{m^{-1}(0)}$ is reduced and irreducible and
$Z = \overline{m^{-1}(0)}$.
We have $M' := \nsmbox{Kalmod}(M,(Z \cap M),\{q\})$,
and $M'' := M \setminus Z$.
We also have $\nsmbox{dim}_{\C}X \geq 3$.

\begin{lemma} \label{lemma:algebraicsubpencil}
There exist algebraic Lefschetz fibrations
\[p' : M' \rightarrow \C \nsmbox{, } p'' : M'' \rightarrow \C\]
such that $p''$ is a subfibration of $p'$ (i.e. $p' \circ (\nsmbox{inclusion}) = p''$).
Also, if $F'$ (resp. $F''$) is a page of $p'$ (resp. $p''$),
then $F'$ is the proper transform of $F''$ in $\nsmbox{Bl}_q X$.
The singularities of $p'$ are contained in $M''$.
\end{lemma}
\proof
Let $Q$ be an effective ample line bundle on $X$ with
support equal to $D$ and such that $Q'' := \overline{m^{-1}(0)} + Q$
is ample.
Let $s'',t''$ be sections of $Q''$
such that ${s''}^{-1}(0) = \overline{m^{-1}(0)} + Q$.
We choose $t''$ such that
\[p'' = \frac{t''}{s''} :M'' \rightarrow \C\] 
is some algebraic Lefschetz fibration on $M''$.
Let $\bar{F}''$ be the closure of one of the smooth fibres of $p''$ in $M$.
We can move the point $p$ to somewhere in the smooth part of
$\bar{F}'' \cap Z$ as the
smooth part of $Z$ is connected
(as $Z$ is irreducible);
$M'$ is unchanged up to Stein deformation. NB here we use
$\nsmbox{dim}_\C X \geq 2$.

Remember $b$ is the blowdown map $b : \nsmbox{Bl}_q X \rightarrow X$.
Let $s' := b^*s''$ and $t' := b^*t''$.
Let $\Delta$ be the exceptional divisor $b^{-1}(p)$.
The divisor ${s'}^{-1}(0)$ is equal to $\Delta + \nsmbox{other divisors}$.
We can choose an effective divisor $K'$ with support equal
to the boundary divisor $D'$ of $M'$ in $\nsmbox{Bl}_q X$ such
that $K'' := K' - \Delta$ is ample.
Hence, we can choose a meromorphic section $h$ of $K''$
whose zero set is contained in $D'$, and such that $h$ has
a pole of order $1$ along the exceptional divisor and such that
$h$ is holomorphic away from $D' \cup \Delta$. Let $L$ be
the line bundle associated to $K''$.
This means that $s' \otimes h \in H^0(\O(L \otimes b^*Q''))$ is non-zero
away from $D'$.
We let $p' := \frac{t' \otimes h}{s' \otimes h}$.
This means that $p'|_{M''} = p''$.
Because $q$ is in the smooth locus of $\bar{F}'' \cap Z$ and
$\bar{F}''$ is transverse to $Z$, we have that the closure of any
smooth fibre of $p'$ intersects each stratum of $D'$
transversally.

We can choose holomorphic coordinates $z_1 \cdots z_n$
on an open set $U$ of $p$ and a holomorphic
trivialisation of $Q''$ such that $s'' = z_1$
and $t'' = z_2$. We then blow up at the point $p$.
Locally around $p$, we have a subvariety of $U \times \P^n$
defined by $Z_iz_j = Z_jz_i$ where $Z_1 \cdots Z_n$ are projective
coordinates for $\P^n$. We choose the chart $Z_1 = 1$.
This has local holomorphic coordinates $z_1,Z_2,Z_3,\cdots,Z_n$.
We can choose a trivialisation of $K''$ so that
the section $h$ is equal to $1 / z_1$. This means that locally
$b^*s'' = Z_1$ and $b^*t'' = Z_2z_1$. Hence locally,
$s' = 1$ and $t' = Z_2$ which means that $p' = Z_2$.
This means that $p'$ has no singular points near $\Delta$.
Hence $p'$ is also an algebraic Lefschetz fibration which coincides
with $p''$ away from $\Delta$ and such all the singular points of $p'$
are the same as the singular points of $p''$.
\qed

\begin{lemma} \label{lemma:holomorphiccurvesinF}
Let $F'$ (resp. $F''$) be a fibre of $p'$ (resp. $p''$).
Let $K$ be a compact set in $F''$. There is a Stein
structure $J$ on $F'$ (depending on $K$) such that
any $J$-holomorphic
$u : T \rightarrow F'$, where $T$ is a compact
Riemann surface with boundary, has the property that
$u(T) \subset F''$ if $u(\partial T) \subset F''$.
\end{lemma}
\proof
Let $G$ be the closure of $F''$ in $M$.
Then $F'$ is biholomorphic to 
$\nsmbox{Kalmod}(G,G \cap Z,\{q\})$.
Let $\phi_G$ be a Stein function for $G$. The compact set
$K$ is contained in $\phi_G^{-1}(C)$ for some large $C$.
Let $q'$ be a point outside
$\phi_G^{-1}(C)$ which is contained in the smooth part of $G \cap Z$.
Because $\dim_\C G \cap Z \geq 2$, we can assume that $q$ and $q'$
are in the same irreducible component $U$ of $G \cap Z$. This is
where we use the assumption that $\nsmbox{dim}_{\C}X \geq 3$.
The manifold $G':=\nsmbox{Kalmod}(G,G \cap Z,\{q'\})$ is naturally a
Stein manifold by Example \ref{example:algebraicsteinmanifold}.
By the comment after Lemma \ref{lemma:symplecticsubmanifold},
we have that $G'$ is Stein deformation equivalent
to $F'$ such that it also induces a Stein deformation on $M''$.
The Stein deformation is induced from moving $q'$
smoothly down to $q$ inside the smooth part of $U \subset G \cap Z$
(Note: $U$ is irreducible, hence the smooth part of $U$ is connected). 
This induces a Stein deformation of $M'$ and $M''$ which in turn induces
a Stein deformation of $F'$ and $F''$.

From now on we assume that the symplectic structures on $F'$
and $G'$ are complete by \cite[Lemma 6]{SS:rama}. We can also
ensure that the above Stein deformation between $F'$ and $G'$
is complete and finite type by the same lemma, hence by
\cite[Lemma 5]{SS:rama} we have a symplectomorphism
$h: F' \rightarrow G'$ induced by this Stein deformation.
Let $J_{G'}$ be the natural
complex structure on $G'$. Let $J := h^*J_{G'}$.
Then if $T$ is a $J$-holomorphic curve in $F'$ with
boundary inside $K$ then $h(T)$ is a holomorphic
curve in $G'$ with boundary in $h(K)$. 
We can blow down this curve to give a holomorphic curve $T'$
in $G$ with boundary in $b(h(K))$. 
We can ensure that $b(h(K))$ is contained in $\phi_G^{-1}(C)$.
If $T$ passes through the blowup of $q$, then $T'$ passes
through $q'$. This means that $\phi_G \circ T'$ has an interior
maximum outside $\phi_G^{-1}(C)$, but this is impossible.
Hence the curve $T$ must be contained in $F''$.
\qed

\bigskip

We can now apply the above lemmas to prove Theorem \ref{theorem:KalimanTransfer}.
We can apply 
Lemma \ref{lm:exacttoalgebraiclf} to $p'$ and $p''$
to get symplectic fibrations $(M',p',\theta'_1)$ and
$(M'',p'',\theta''_1)$. 
These fibrations are Lefschetz without boundary. 
We can
cut down the fibres to Stein domains $\bar{F'}$ and $\bar{F''}$
where $\bar{F''}$, $\bar{F'}$ are large enough so that the support of
all the monodromy maps of $(M',p',\theta'_1)$ are contained
in $\bar{F''}$ and $\bar{F''} \subset \bar{F'}$.
We can also remove the cylindrical end from the base. This
will make $p'$ and $p''$ into Lefschetz fibrations $(E',\pi')$
and $(E'',\pi'')$ respectively. Note that if we have a holomorphic
curve $T$ in $\bar{F'}$ with boundary in $\bar{F''}$, 
Lemma \ref{lemma:holomorphiccurvesinF} implies that it is contained in $F'' \cap \bar{F'}$.
The Stein maximum principle \cite[Lemma 1.5]{Oancea:survey}
ensures that $T$ is contained in $\bar{F''}$.
Hence we get that Theorem \ref{theorem:KalimanTransfer} is a consequence 
of Lemma \ref{lm:exacttoalgebraiclf} and Lemma
\ref{lemma:holomorphiccurvesinF}.

\section{Appendix B: Stein structures and cylindrical ends}

The problem with Stein structures is that the complex
structure associated with them does not behave well with
respect to cylindrical ends.
Cylindrical ends here means that near infinity, the convex
symplectic manifold is exact symplectomorphic to
$(\Delta \times [1,\infty), r\alpha)$ where $r \in [1,\infty)$
and $\alpha$ is a contact form on $\Delta$.
The almost complex structure is convex with respect to this
cylindrical end if $dr \circ J = -\alpha$. We will
deal with this problem in this section.

Let $(M,J,\phi)$ be a complete finite-type Stein manifold with
$\theta = -d^c\phi$ and $\omega = d\theta$. Let
$c \gg 0$ be greater than the highest critical value of $\phi$.
\begin{theorem} \label{theorem:steintoalmoststein}
There exists a complete finite type convex symplectic structure
$(M,\theta_1)$ with the following properties:
\begin{enumerate}
\item
It has a cylindrical end with an almost
complex structure $J_1$ which is convex at infinity.
\item $J_1 = J$ and $\theta_1 = \theta$ in the region $\{\phi \leq c\}$.
\item \label{item:holomorphiccurveconstraint}
Any $J_1$ holomorphic curve with boundary in $\{\phi = c\}$ is contained
in  $\{\phi \leq c\}$.
\item
It is convex deformation equivalent to $(M,\theta)$ via
a convex deformation $(M,\theta_t)$ where
$\theta_t|_{\{\phi \leq c\}} = \theta|_{\{\phi \leq c\}}$ for $t \in [0,1]$.

\end{enumerate}
\end{theorem}

\proof
Let $\lambda := \nabla \phi$ and $\Delta := \phi^{-1}(c+1)$.
We define $G : M \rightarrow \R$, $G =  \frac{1}{\|\nabla \phi\|^2}$
where
$\|.\|$ is the norm defined using the metric $\omega(\cdot,J \cdot)$.
Let
$\lambda' := G\lambda$.
Let $F_t : M \rightarrow M$ be the flow of $\lambda'$.
This exists for all time because $\phi$ is unbounded and
${\mathcal L}_{\lambda'} \phi = 1$
which implies that $\phi$ increases linearly with $t$
(${\mathcal L}$ here means Lie derivative).
We have an embedding $\Phi : \Delta \times [1,\infty) \rightarrow M$
defined by $\Phi(a,r) = F_{\log{r}}(a)$ where $a \in \Delta \subset M$
and $r \in [1,\infty)$.  Also,
${\mathcal L}_{\lambda'} \theta = G\theta$.
Hence,
$\Phi^*(\theta) = f\alpha$ where $f : \Delta \times [1,\infty) \rightarrow \R$,
 $f(a,r) := 1+\int_0^r (G \circ \Phi)(a,t) dt$ and
$\alpha$ is the contact form $\theta|_{\Delta}$ on $\Delta$.

We will now deform the $1$-form $f\alpha$ to a $1$-form
$f'\alpha$ such that $f' = f$ near $r= 1$ and $f' = r$ near infinity.
We define $\theta_1$ to be equal to $f'\alpha$ in this cylindrical
end and equal to $\theta$ away from this end.
This means that for $r$ large, we have a cylindrical end
with $1$-form $f'\alpha = r\alpha$.
If we have a function
$g : \Delta \times [1,\infty) \rightarrow \R$, then
$d(g\alpha)$ is non-degenerate if and only if
$\frac{\partial g}{\partial r} > 0$. Also, the
Liouville vector field associated to $g\alpha$ is
$(g / \frac{\partial g}{\partial r}) \frac{\partial}{\partial r}$,
and hence we have that this Liouville vector field is transverse
to every level set $\{r = \nsmbox{const}\}$ and pointing outwards.
If $(g / \frac{\partial g}{\partial r})$ is bounded above by any polynomial,
then the respective Liouville vector field is complete.
We define $f' : \Delta \times [1,\infty) \rightarrow \R$ such that
$f' = f$ near $r = 1$, $f' = r$ near infinity and $\frac{\partial f'}{\partial r} > 0$.
This gives a complete finite type convex symplectic structure $\theta_1$
on $M$ as we can extend $f'\alpha$ outside $M$ as $f' = f$ near $r=1$.
We can join $f$ to $f'$ via a smooth family of functions
with $f_t$ ($t \in [0,1]$) where
$\frac{\partial f_t}{\partial r} > 0$ and such that $f_t = f$ near $r=1$.
This gives us a convex deformation from $\theta$ to $\theta'$.

We now need to construct our almost complex structure $J_1$.
We have ${\mathcal L}_{\lambda'} \phi = G d\phi(\nabla \phi) = 1$.
This means that $\Phi^*(\phi) = \log{r}$ so the level sets of $\phi$
coincide with the level sets of $\log{r}$.
By abuse of notation we will just write $J$ for the pullback
$\Phi^*J$ and we will write $\phi$ for $\log{r}$.
We have two orthogonal symplectic vector subbundles of
the tangent bundle $\Phi^*(TM) = T(\Delta \times [1,\infty))$ whose
direct sum is the entire tangent bundle
(the symplectic structure we are dealing with here is $\theta_1$). These are:
$V_1 := \nsmbox{Ker}(\theta_1) \cap \nsmbox{Ker}(dr)$ and
$V_2 := \nsmbox{Span}(\frac{\partial}{\partial r}, X_r)$ where $X_r$
is the Hamiltonian flow of $r$. The problem is that $J$
is not necessarily compatible with $d\theta_1$, so we need
to deform it so that it is. However, near $r=1$, $J$ is in fact
compatible with $d\theta_1$ because $\theta= \theta_1$
in some region $\Xi := \{r \leq 1 + \epsilon\}$.
Inside $\Xi$, we have that:
$J|_{V_1}$ and $J|_{V_2}$ are holomorphic subbundles of $\Phi^*(TM)$.
There exists a complex structure $J_{V_1}$ (resp. $J_{V_2}$)
on the vector bundle $V_1$ (resp. $V_2$)
compatible with $d\theta_1|_{V_1}$ (resp.  $d\theta_1|_{V_2}$)
such that, $J_{V_1} = J|_{V_1}$
(resp. $J_{V_2} = J|_{V_2}$) when restricted to $\Xi$.
Because $V_1 \bigoplus V_2 = \Phi^*(TM)$, this gives us an almost
complex structure $J_1$ on $\Phi^*(TM)$ compatible with $d\theta_1$
which is equal to $J$ in the region $\Xi$.
We can choose $J_{V_1}$ and $J_{V_2}$ so that
$J_{V_2}(\frac{\partial}{\partial r}) = -\frac{1}{r}X_r$ for $r \gg 0$
and $J_{V_1}$ is invariant under the flow of
$\frac{\partial}{\partial r}$ for $r \gg 0$. This ensures that $J_1$
is convex at infinity. Also, we have that $r$ is plurisubharmonic
with respect to $J_1$ hence any $J_1$ holomorphic curve with
boundary in $\{r=1\}$ is contained in $\{r \leq 1\}$. Hence
property (\ref{item:holomorphiccurveconstraint}) is satisfied.
\qed

\section{Appendix C: Transfer maps and handle attaching} \label{section:transfermapsandhandleattaching}

The purpose of this section is to show that symplectic homology is additive
as a ring under end connect sums. This was already done by Cieliebak
in \cite{Cieliebak:handleattach} but without taking into account the ring structure.
Throughout this section, let $(M,\theta)$, $(M',\theta')$ be compact convex 
symplectic manifolds such that $M'$ is an exact submanifold of $M$ of codimension $0$.
We let $C:= N \times [1,\infty)$ be a cylindrical end of $M$ where $\theta = r\alpha$,
$\alpha$ is a contact form on $N$, and $r$ is the coordinate for $[1,\infty)$. 
Similarly we have a cylindrical end $C'$ of $M'$.
Let $H : M \rightarrow \R$ be an admissible Hamiltonian with
an almost complex structure $J$, convex at infinity.
Let $SH_*^{(-\infty,a)}(M,H,J)$ be the group generated by
orbits of action $< a$. For $b \geq a$, we define
\[SH_*^{[a,b)}(M,H,J):= SH_*^{(-\infty,b)}(M,H,J)/ SH_*^{(-\infty,a)}(M,H,J).\]

\subsection{Weak cofinal families} \label{subsection:weakcofinal}

\begin{defn} \label{defn:weakadmissible}
We say that the pair
$(H,J)$ is {\bf weakly admissible} if there exists an
$f : N \rightarrow \R$ and a constant $b$ such that
for $r \gg 0$, $H = re^{-f} + b$ and
$d(re^{-f}) \circ J = -\theta$.
\end{defn}

Every admissible pair $(H,J)$ is weakly admissible with $f=\nsmbox{const}$.
Symplectic homology $SH_*(M)$ is defined as a direct limit of $SH_*(M,H,J)$
with respect to admissible pairs ordered by $\leq$. We wish to replace
``admissible'' with ``weakly admissible''.
The
reason why we wish to do this is because in section
\ref{subsection:handleattaching} we carefully construct a cofinal
family of weakly admissible pairs to show that
symplectic homology behaves well under end connect sums.
We construct a partial order $\leq$
on weakly admissible pairs as follows:
$(H_0,J_0) \leq (H_1,J_1)$ if and only if $H_0 \leq H_1$. 
We will show that:
\[SH_*(M) :=  \underset{(H,J)}{\varinjlim} \nsmbox{ } SH_*(M,H,J)\]
where the direct limit is taken over weakly admissible
pairs $(H,J)$ ordered by $\leq$. Note that a family
of weakly admissible Hamiltonians $(H_s,J_s)$ is cofinal
with respect to $\leq$
if the corresponding functions $f_s : N \rightarrow \R$
tend uniformly to $-\infty$ as $s$ tends to $\infty$.
In order to ensure that this direct limit exists,
we will show that if $(H_0,J_0) \leq (H_1,J_1)$,
then there is a natural map of rings
$SH_*(M,H_0,J_0) \rightarrow SH_*(M,H_1,J_1)$.

This map will be a continuation map.
In order for a continuation map to be well defined, 
we need a family of Hamiltonians $T_s$
joining $H_0$ and $H_1$ such that solutions
of the parameterized Floer equation 
$\partial_s u + J_t \partial_t u = \nabla^{g_t} T_s$
joining orbits of $H_0$ and $H_1$
stay inside some compact set. To ensure this,
we flatten $H_i$ so that it is constant outside
some large compact set, and so that all the additional
orbits created have very negative action. We do this as
follows:

Let $D$ be a constant such that any orbit of $H_0$ or $H_1$
has action greater than $D$. Then 
$SH_*(M,H_i,J_i) \cong SH_*^{[D,\infty)}(M,H_i,J_i)$.
Near infinity, we have that $H_i = R_i + b_i$ where
$R_i = re^{-f_i}$. We wish to create a new Hamiltonian $K_i$
such that $K_i = H_i$ on $R_i \leq B$ where $B \gg 0$,
and such that $K_i$ is constant in $\{R_i > B+1\}$
where all the additional orbits have action less than $D$.
We assume that all the orbits of $H_0$ and $H_1$ lie in a compact set.
We have a cylindrical end
$C_i := N_i \times [K,\infty)$ where $N_i$ is the contact manifold
$\{re^{-f_i}=1\}$ with contact form $\theta|_{N_i}$ and $R_i$
is the coordinate for $[K,\infty)$. So, $H_i$ is linear
with slope $1$ on this cylindrical end. Because all the orbits
of $H_i$ lie in a compact set, there are no Reeb orbits of length $1$
in the contact manifold $N_i$. Choose $\epsilon > 0$ such that
the length of any Reeb orbit of $N_0$ or $N_1$ is of distance more than $\epsilon$
from $1$.
We assume that $B$ is large enough so that $H_i$
is linear with respect to the cylindrical end $C_i$ in $R_i \geq B$
and such that $B\epsilon > -b_i - D$ for $i=0,1$.
Finally we let $K_i$ be equal to $H_i$ in the region
$\{R_i \leq B\}$, and $K_i = k_i(R_i)$ where $k_i$ is
constant for $R_i \geq B+1$, and $k_i = R_i + b_i$ near $B$
and $k_i' \leq 1$.
This means that the orbits of $K_i$ in $\{R_i \leq B\}$ are the
same as the orbits of $H_i$, and the orbits
of $K_i$ in $\{R_i > B\}$ have action 
\[R_ik_i' - k_i < R_i(1-\epsilon) - R_i - b_i < -B\epsilon - b_i < D.\]
Hence all the orbits of $K_i$ of action greater than $D$ are the
same as the orbits of $H_i$.

We now wish to create an almost complex structure $J_i^0$ as follows:
we let $J_i^0 = J_i$ for $R_i \leq B$ and for $R_i \geq B+1$, we let $J_i^0$
be convex with respect to the cylindrical end $C$ (i.e. $dr \circ J = -\theta$).
Let $u$ be a cylinder or pair of pants satisfying
Floer's equation (\cite[Formula 8.1]{Seidel:biasedview}) with respect
to $(K_i,J_i^0)$
such that each cylindrical end of $u$ limits to a periodic orbit 
(or multiple of a periodic orbit in the pair of pants case) inside
$\{R_i \leq B\}$. 
By \cite[Lemma 7.2]{SeidelAbouzaid:viterbo} we have that $u$ is contained in  $\{R_i \leq B\}$.
Note that we really perturb these Hamiltonians so that all the orbits are non-degenerate,
and lie in a compact set. Lemma 7.2 from \cite{SeidelAbouzaid:viterbo} still
works in this case, as we can ensure the Hamiltonian stays the same in the region
$B-1 < R_i < B$. From now on if we deal with Hamiltonians which are constant at
infinity, we are really perturbing them in such a way that this convexity argument
from \cite{SeidelAbouzaid:viterbo} still holds and such that $SH_*$ is well defined
for this Hamiltonian.
Hence
\[SH_*(M,H_i,J_i) \cong SH_*^{[D,\infty)}(M,H_i,J_i^0) \cong  SH_*^{[D,\infty)}(M,K_i,J_i^0).\]
We wish to create a continuation map:
\[SH_*^{[D,\infty)}(M,K_0,J_0^0) \rightarrow SH_*^{[D,\infty)}(M,K_1,J_1^0).\]

There exists a family of Hamiltonians $A_s$ connecting $K_0$ and $K_1$,
and such that $A_s$ is monotonically increasing and $A_s$ is constant
at infinity. We also join $J_0^0$ and $J_1^0$ with a family of almost
complex structures which are convex with respect to the cylindrical end $C$ (i.e. $dr \circ J = -\theta$).
If $f_0$ and $f_1$ are constant, then we have a monotone increasing family of admissible
Hamiltonians $H_s^1$ joining $H_0$ and $H_1$, and almost complex structures $J_s$
joining $J_0$ and $J_1$.
The standard continuation map
\[SH_*(M,H_0,J_0) \rightarrow SH_*(M,H_1,J_1)\] involves counting
solutions of a parameterized Floer equation with respect to $(H_s,J_s)$.
We wish to show that this map is the same as the above continuation
map from $K_0$ to $K_1$.
In order to do this we construct an explicit family $(A_S,J_s)$ of Hamiltonians
and almost complex structures so that the continuation map
\[SH_*^{[D,\infty)}(M,K_0,J_0^0) \rightarrow SH_*^{[D,\infty)}(M,K_1,J_1^0)\]
coincides with the standard continuation map \[SH_*(M,H_0,J_0) \rightarrow SH_*(M,H_1,J_1)\]
under the isomorphism $SH_*(M,H_i,J_i) \cong SH_*^{[D,\infty)}(M,K_i,J_i^0)$.
The Hamiltonians $H_s^1$ are of the form $h_s(r)$ for $r \geq P$ where $h_s'$ is constant.
We assume that the almost complex structures $J_s$ are convex with respect to the cylindrical end $C$
for $r \geq P$.
Hence \cite[Lemma 7.2]{SeidelAbouzaid:viterbo} ensures all the Floer trajectories
with respect to $(H_s,J_s)$ stay inside the compact set $r \leq P$.
Also there is a constant $P'$ such that $K_i$ is a function of $r$ for $r \geq P' \geq P$.
The definition of $K_i$ depends on a parameter $B$ which can be arbitrarily large.
We can choose $B$ large enough so that $K_i = H_i$ in the region $r \leq P'$.
We choose the functions $A_s$ joining $K_0$ and $K_1$ so that $A_s$ is a function
of $r$ for $r \geq P'$. We can also assume that $J_s^0 = J_s$.
In order to show that the maps are the same, we need to show that
any Floer trajectory associated to $(A_s,J_s^0)$ connecting
orbits inside $\{r \leq P\}$ is contained in $\{r \leq P\}$.
This follows from  \cite[Lemma 7.2]{SeidelAbouzaid:viterbo}.

\subsection{Transfer maps} \label{section:transfermaps}

In this section we will construct a natural ring homomorphism:
$SH_*(M) \rightarrow SH_*(M')$.
We say that $H$ is called {\bf transfer admissible}
if $H \leq 0$ on $M'$.
We have:
$SH_*(M) = \underset{(H,J)}{\varinjlim} \nsmbox{ } SH_*(M,H,J)$
where the direct limit is taken over transfer admissible
Hamiltonians ordered by $\leq$.

\begin{lemma} \label{lemma:transfermap}
We have an isomorphism of rings,
\[\underset{(H,J)}{\varinjlim} \nsmbox{ } SH_*^{[0,\infty)}(M,H,J) \cong SH_*(M')\]
where the direct limit is taken over transfer admissible Hamiltonians.
\end{lemma}
\proof
We construct a particular cofinal family of transfer admissible
Hamiltonians $H_i$ and show the above isomorphism of rings.
We can embed $\widehat{M'}$ into $\widehat{M}$ by Lemma
\ref{lemma:completionsubmanifold}. Our cylindrical
end $C'$ is then a subset of $M$. We assume that the
action spectrum ${\mathcal S} := {\mathcal S}(\partial M')$ is discrete
and injective.
Let
$k : \N \rightarrow {\R \setminus {\mathcal S}}$ be a function such that 
$k(i)$ tends to infinity as $i$ tends to infinity. 
Let
$\mu : \N \rightarrow \R$ be defined by $\nsmbox{dist}(k(i),{\mathcal S})$
($\nsmbox{dist}(a,B)$ is the shortest distance between $a$ and $B$).
From
now on we just write $k$ instead of $k(i)$,
and similarly for $\mu$.

Define:
\[A = A(i) := 6k / \mu > k > 1.\]
We can assume that $A>k>1$ because we can choose $k(i)$ to make
$\mu(i)$ arbitrarily small whilst $k(i)$ is large. We also let
$\epsilon:=\epsilon(i)$ tend to $0$ as $i$ tends to infinity.
We assume that $H_i|_{M'} \leq 0$, and has slope
$k(i)$ on $1 + \epsilon/k \leq r' \leq A - \epsilon/k$.
We also assume that on $1 \leq r' \leq A$, $H_i = h(r')$ for
some function $h$ where $h$ has non-negative derivative $\leq k$.
For $A \leq r' \leq A+1$, we assume that $H_i$ is constant.
Let $B$ be this constant. $B$ is arbitrarily
close to $k(A-1)$. We can assume that $B \notin {\mathcal S}$.
We now describe $H_i$ on the cylindrical end $C$.
We keep $H_i$ constant until we reach $r = A+1+P$ where
$P$ is some constant large enough so that
$\{r' \leq 1\} \subset \{r \leq P\}$. This means
that $\{r' \leq A+1\} \subset \{r \leq A+1 + P\}$ as long
as we embed $C'$ in the same way as Lemma
\ref{lemma:completionsubmanifold}.
We then let $H_i$ be of the form $f(r)$ for $r \geq A+1+P$
where $f' < \frac{1}{2}k$ and has slope $\frac{1}{2}k$ for 
$r > A+1 + P + \epsilon/k$.

Here is a picture of what we have:

\begin{fig} \label{fig:transferhamiltonian}

 \[
 \xy
 (5,0)*{}="A"; (5,-50)*{}="B";
 "A"; "B" **\dir{-};
 (0,-35)*{}="C"; (60,-35)*{}="D";
 {\ar@{->} "C";"D"};
 (5,-40)*{}="E"; (60,-40)*{}="F";
 "E"; "F" **\dir{.};
 (5,-10)*{}="E1"; (60,-10)*{}="F1";
 "E1"; "F1" **\dir{.};
 (40,-35)*{}="VL1"; (40,-10)*{}="VL2";
 "VL1"; "VL2" **\dir{.};
 (45,-35)*{}="VR1"; (45,-10)*{}="VR2";
 "VR1"; "VR2" **\dir{.};
 (18,-35)*{}="G";
 "C"+(5,-3); "G"+(0,-3) **\crv{"C"+(12,2) & "G"+(-13,-5)};
 "G"+(0,-3); "G"+(7,0) **\crv{"G"+(2,-3) & "G"+(3,-4)};
 "G"+(7,0); "VL2"+(0,-3) **\dir{-};
 "VL2"+(0,-3); "VR2"; **\crv{"VL2"+(3,0) & "VR2" + (-5,-3)};
 "VR2"; "F1" **\dir{-};
 (65,-10)*{}="HDRH";
 "F1"; "HDRH" **\dir{.};
 (65,-35)*{}="HDRL";
 "D"; "HDRL" **\dir{.};
 "E"+(-2,0) *{\scriptstyle -\epsilon};
 "E1"+(-2,0) *{\scriptstyle B};
 "VL1"+(-4,2) *{\scriptstyle A-\frac{\epsilon}{k}};
 "VR1"+(2,2) *{\scriptstyle A};
 "G" + (0,2) *{\scriptstyle 1};
 "G" + (40,-2) *{\scriptstyle r'};
 "G" + (7,2) *{\scriptstyle 1 + \frac{\epsilon}{k}};
 "VL2" + (-10,2) *{{\bf H_i}};
 "VL2" + (-8,-12) *{\scriptstyle k};
 (65,0)*{}="SHIFT";
 "SHIFT"+(0,-35)*{}="C1"; "SHIFT"+(40,-35)*{}="D1";
 {\ar@{->} "C1";"D1"};
 "SHIFT" + (0,-10); "SHIFT"+(20,-10) **\dir{-};
 "SHIFT"+(20,-7)*{}="G1";
 "G1"+(0,-3); "G1"+(7,0) **\crv{"G1"+(2,-3) & "G1"+(3,-4)};
 "G1"+(7,0); "G1"+(10,3) **\dir{-};
 "SHIFT" + (20,-10); "SHIFT" + (20,-35) **\dir{.};
 "SHIFT" + (27,-7); "SHIFT" + (27,-35) **\dir{.};
 "SHIFT" + (40,-37) *{\scriptstyle r};
 "SHIFT" + (25,-5) *{\scriptstyle \frac{1}{2}k};
 "SHIFT" + (20,-37) *{\scriptstyle A+1 + P};
 "SHIFT" + (30,-33) *{\scriptstyle A+1 + P + \epsilon/k};
 \endxy
 \]

\end{fig}

The action of an orbit on a level set $r'=a$ 
is $h'(a)a-h(a)$.
The orbits near $r'=1$ have positive action less than or equal
to $k$. Let $p$ be a point on an orbit $o$ lying in the region
$A-\epsilon/k \leq r' \leq A$. The slope $h'(r')$ of $H_i$ at $p$
is $\leq k - \mu$.
Hence, the orbits near $r'=A$ have action
$\leq (k - \mu)A - B = -\mu A + k \rightarrow -\infty$ as 
$i \rightarrow \infty$. So, we can assume that these orbits
have negative action. Also all the orbits in $\{r' \geq A,r \leq A+1 + P\}$
are fixed points, so have action $-B < 0$. Finally,
the orbits in $r > A+1 + P$ have action:
$\leq \frac{1}{2}k.(A+1+P) - B = -\frac{1}{2}k A +  \frac{3}{2}k + \frac{1}{2}Pk \rightarrow -\infty$
as $i \rightarrow \infty$.
Hence, we can assume that
all the orbits of $H_i$ of non-negative action lie in $r' < 1 + \epsilon/k$.

We now need to show that any differential connecting two orbits
of non-negative action is contained entirely in $r' < 1 + \epsilon/k$.
By \cite[Lemma 7.2]{SeidelAbouzaid:viterbo},
any differential connecting orbits of non-negative action
must be contained in $r' < 1 + \epsilon/k$ as all the orbits of non-negative
action are contained in this region. A similar application
of this lemma ensures that a pair of pants satisfying Floer type
equations with similar Hamiltonians where the ends converge
to orbits of non-negative action must be contained in $r' < 1 + \epsilon/k$.
This means that we have maps
$SH_*^{[0,\infty)}(M,H_i,J) \cong SH_*(M',H'_i,J')$,
where $H'_i : \widehat{M'} \rightarrow \R$ has slope $k$.
Taking direct limits gives us a ring isomorphism
\[\underset{(H,J)}{\varinjlim} \nsmbox{ } SH_*^{[0,\infty)}(M,H,J) \cong SH_*(M').\]
\qed

This lemma enables us to define a transfer map:
\[SH_*(M) \cong \underset{(H,J)}{\varinjlim} \nsmbox{ } SH_*(M,H,J) 
\rightarrow \underset{(H,J)}{\varinjlim} \nsmbox{ } SH_*^{[0,\infty)}(M,H,J) 
\cong SH_*(M').\]

A Hamiltonian is called {\bf weakly transfer admissible} if it is weakly admissible
and is negative when restricted to $M'$.
We can combine the above results with the results of section \ref{subsection:weakcofinal}
to construct the above transfer map using a cofinal family
of weakly transfer admissible Hamiltonians.
We will need to construct a cofinal family of weakly transfer admissible
Hamiltonians in section \ref{subsection:handleattaching} to show that
a particular transfer map is an isomorphism of rings.

\bigskip

Here is an application of the transfer map:
\begin{lemma} \label{lemma:nontrivialsymplectichomology}
If $SH_*(M)=0$, then $SH_*(M')=0$.
\end{lemma}
\proof
We have a commutative diagram:
\[
\xy
(0,0)*{}="A"; (40,0)*{}="B";
(0,-20)*{}="C"; (40,-20)*{}="D";
"A" *{H^{n-*}(M)};
"B" *{H^{n-*}(M')};
"C" *{SH_*(M)};
"D" *{SH_*(M')};
{\ar@{->} "A"+(10,0)*{};"B"-(10,0)*{}};
{\ar@{->} "C"+(10,0)*{};"D"-(10,0)*{}};
{\ar@{->} "A"+(0,-4)*{};"C"+(0,4)*{}};
{\ar@{->} "B"+(0,-4)*{};"D"+(0,4)*{}};
"A"+(20,3) *{a};
"C"+(20,3) *{c};
"A"+(2,-10) *{b};
"B"+(2,-10) *{d};
\endxy
\]
Suppose for a contradiction $SH_*(M') \neq 0$. Then
\cite[Section 8]{Seidel:biasedview} says that the
map $d$ is non-zero in degree $n$. Also the map
$a$ is an isomorphism in degree $n$. Hence 
$d \circ a$ is non-zero, and so $c \circ b = d \circ a$
is non-zero. This means that $SH_*(M) \neq 0$ and
we get a contradiction.
\qed
\begin{corollary} \label{corollary:embeddinginsubcritical}
If $M$ is subcritical and $SH_*(M') \neq 0$,
then $M'$ cannot be embedded in $M$ as an exact
codimension $0$ submanifold.
In particular, if $H_1(M')=0$ then $M'$ cannot be symplectically
embedded into $M$.
\end{corollary}
\proof
By one of the applications of \cite{Oancea:kunneth},
we have that $SH_*(M)=0$ because $M$ is subcritical.
The result follows from the above lemma.
\qed

\subsection{Handle attaching} \label{subsection:handleattaching}

In this section we will prove Theorem \ref{thm:handleattaching}.
Here is the statement of Theorem \ref{thm:handleattaching}:
Let $M,M'$ be finite type Stein manifolds of real dimension greater than $2$, then \\
$SH_*(M \#_e M') \cong SH_*(M) \times SH_*(M')$ as rings. 
Also the transfer map $SH_*(M \#_e M') \rightarrow SH_*(M)$
is just the natural projection
\[SH_*(M) \times SH_*(M') \twoheadrightarrow SH_*(M).\] 
To prove this, we will show that attaching a symplectic $1$-handle to a compact convex symplectic
manifold of dimension $\geq 2$ does not change symplectic homology.
We will describe in more detail what it means to attach a symplectic $1$-handle later.
In fact we show that if $A$ is a compact convex symplectic manifold, and
$A'$ is equal to $A$ with a symplectic $1$-handle attached, then
the natural transfer map $SH_*(A') \rightarrow SH_*(A)$ is an isomorphism.
This proves Theorem \ref{thm:handleattaching} for the following reason:
The Stein manifold $M$ (resp. $M'$) is convex deformation equivalent to
$\widehat{N}$ (resp. $\widehat{N'}$) where $N$ (resp. $N'$)
is a compact convex symplectic manifold.
The end connect sum $M \#_e M'$ is convex deformation equivalent to $\widehat{N''}$ where
$N''$ is the disjoint union $N \coprod N'$ with a symplectic $1$-handle joining each
connected component.
The symplectic homology of $N \coprod N'$ is the direct
product $SH_*(N) \times SH_*(N')$. Also, the transfer map $SH_*(N \coprod N') \rightarrow SH_*(N)$
is the natural projection $SH_*(N) \times SH_*(N') \twoheadrightarrow SH_*(N)$.
Hence assuming that adding a $1$-handle does not change symplectic homology, we have that
$SH_*(N'')$ is isomorphic to $SH_*(N) \times SH_*(N')$, hence $SH_*(M \#_e M')$ is isomorphic
to this product. This implies that $SH_*(M \#_e M') \cong SH_*(M) \times SH_*(M')$
because $M$ (resp. $M'$) is convex deformation equivalent to $N$ (resp. $N'$).
Also the natural transfer map from $SH_*(M \#_e M')$ to $SH_*(M)$ is the natural
projection $SH_*(M) \times SH_*(M') \twoheadrightarrow SH_*(M)$ because this corresponds
to the composition of maps 
\[SH_*(N'') \overset \cong \rightarrow SH_*(N \coprod N') \overset \cong \rightarrow SH_*(N) \times SH_*(N')\twoheadrightarrow SH_*(N)\]

We will now describe handle attaching in detail as in  
\cite[Section 2.2]{Cieliebak:handleattach}. The paper
\cite{Eliashberg:steintopology} or 
\cite[Theorem 9.4]{CieliebakEliashberg:symplecticgeomofsteinmflds}
ensures that this construction corresponds to attaching
a Stein $1$-handle.
We will define $\phi$, $p_i$, $q_i$, $X$, $\omega$, $\psi(x,y)$ 
as in \cite[Section 2.2]{Cieliebak:handleattach}. We will now remind
the reader what these variables are:
We set $k=1$, so we are describing $1$-handles only.
We let $\R^{2n}$ have coordinates $(p_1,q_1,\dots,p_n,q_n)$.
\[\omega := \sum_i dp_i \wedge dq_i,\]
\[\phi := \frac{1}{4} \sum_{i=1}^{n-1} \left( q_i^2 + p_i^2\right) 
 + q_n^2-\frac{1}{2}p_n^2,\]
\[X := \nabla \phi =  \frac{1}{2}\sum_{i=1}^{n-1} \left(
\frac{\partial}{\partial q_i} + \frac{\partial}{\partial p_i} \right)
+ 2\frac{\partial}{\partial q_n} -\frac{\partial}{\partial p_n}.\]
$\psi$ is a function of $x$ and $y$ where:
\[x := \sum_{i=1}^{n-1}\left( A_i q_i^2 + B_i p_i^2 \right),\]
\[y := B_n p_n^2,\]
and $A_i,B_i > 0$ are constants. It satisfies $X.\psi > 0$
provided that:
\[\frac{\partial \psi}{\partial x} \geq 0, \nsmbox{ }
\frac{\partial \psi}{\partial y} \leq 0, \nsmbox{ }
\frac{\partial \psi}{\partial x}(x,0) > 0, \nsmbox{ }
\frac{\partial \psi}{\partial y}(y,0) < 0, \nsmbox{ }
\]
and the partial derivatives are not simultaneously $0$.
We can choose $\psi$ so that the level sets
$\{\phi = -1\}$ and $\{\psi=1\}$ agree outside some compact set.
This ensures that when we glue the handle onto our
convex symplectic manifold, it still has a smooth boundary
so we don't have to smooth the handle once we have attached it.

The handle $H = H_1^{2n} := \{\phi \geq -1\} \cap \{\psi \leq 1\}$.
We define $\partial^-H$ to be the boundary $\{\phi=-1\}$.
We can ensure that the only $1$-periodic orbit of $\psi$ is
the critical point at the origin by \cite[Section 2.2]{Cieliebak:handleattach}.
We wish to construct a family of $1$-handles (constructed in the same
way as $H$) $(H_l)_{l \in \N}$ with the following properties:
\begin{enumerate}
\item $H_{l+1} \subset H_l$
\item The attaching region $\partial^-H_{l+1}$ is a subset of  $\partial^-H_l$
\item As $l$ tends to infinity, $H_l$ converges uniformly to the core of the handle.
\end{enumerate}
This can be done by shrinking $\psi$.
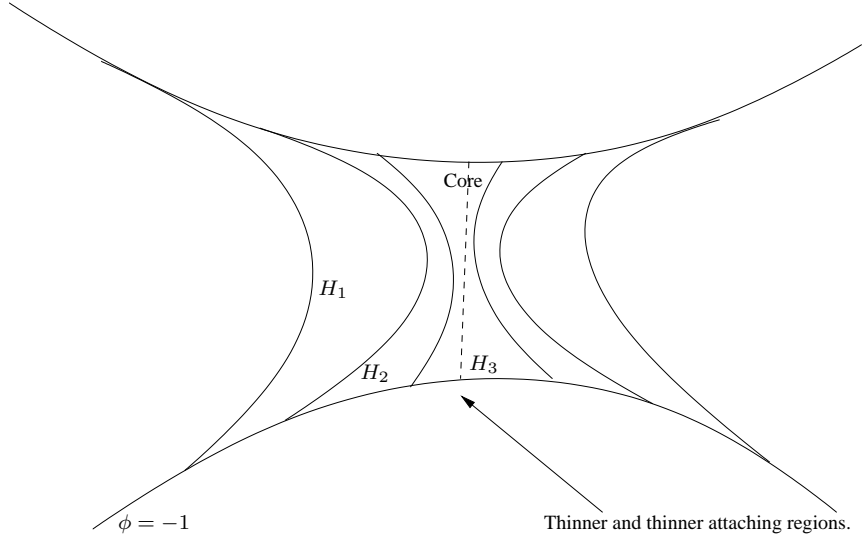
\begin{figure}[H]
 \centerline{
   \scalebox{1.0}{
     \input{handles.pstex_t}
   }
 }
\label{figure:handles}
\end{figure}
Let $M$ be a compact convex symplectic manifold. After a deformation,
we can assume that the boundary of $M$ has a region $A$ which is
contactomorphic to the attaching region $\partial^- H_1$. We
can also ensure that the period
spectrum of $\partial M$ is discrete and injective
(we might have to deform the region $A$ and $\phi$ 
slightly and hence all the handles).
We can use the region $A$ to attach the handle $H_l$ to $M$
to create a new compact convex symplectic manifold
$M_l := M \cup_{\partial^- H_l} H_l$.
We have that $M_{l+1} \subset M_l$ and the boundary of each
$M_l$ is transverse to the Liouville vector field on $M_1$.
Let $K$ be an admissible Hamiltonian on $\widehat M$. We assume that
$K$ has slope $S$ in a neighbourhood of $\partial M$. 
We choose $l$ large enough so that the
attaching region $P:=\partial^- H_l$ has the property that
a Reeb flowline outside $P$ intersecting $P$ twice
has length greater than $S$. We can now extend the Hamiltonian
$K$ to a Hamiltonian $K' : M_l \rightarrow \R$
using the function $B\psi$ where $B$
is some constant. Hence $\partial M_l$ is a level set of $K'$ and
$K'$ is linearly increasing on a neighbourhood of $\partial M_l$.
Hence we can extend $K'$ to an admissible Hamiltonian on $\widehat M_l$.
The periodic orbits of $K'$ are the same as the periodic orbits of $K$
with an extra fixed point at the origin of the $1$-handle.
We can ensure that the index of the extra fixed point at the origin
of the $1$-handle has index strictly increasing as $S$
increases (see the last part of the proof of Theorem 1.11 
in \cite[Section 3.4]{Cieliebak:handleattach}).
Because $\partial M_l$ is transverse to the Liouville field
of $\partial M_1$, we have that $\widehat M_1 = \widehat M_l$
and $K'$ is weakly admissible.
Hence we have a cofinal family of weakly transfer admissible Hamiltonians
$K'$. The only orbit outside $M \subset \widehat{M_1}$ has arbitrarily
large index, hence these $K'$'s induce a transfer isomorphism of rings
$SH_*(M) \rightarrow SH_*(M_1)$.
This proves Theorem \ref{thm:handleattaching}.

\bibliography{references}

\end{document}

%% file: radiallines.pstex_t
\begin{picture}(0,0)%
\includegraphics{radiallines.pstex}%
\end{picture}%
\setlength{\unitlength}{2763sp}%
\begingroup\makeatletter\ifx\SetFigFont\undefined%
\gdef\SetFigFont#1#2#3#4#5{%
  \reset@font\fontsize{#1}{#2pt}%
  \fontfamily{#3}\fontseries{#4}\fontshape{#5}%
  \selectfont}%
\fi\endgroup%
\begin{picture}(7599,4787)(2089,-6204)
\put(4801,-4111){\makebox(0,0)[lb]{\smash{{\SetFigFont{8}{9.6}{\rmdefault}{\mddefault}{\updefault}{\color[rgb]{0,0,0}Origin $0$}%
}}}}
\put(4201,-1561){\makebox(0,0)[lb]{\smash{{\SetFigFont{8}{9.6}{\rmdefault}{\mddefault}{\updefault}{\color[rgb]{0,0,0}Radial lines coming from the origin}%
}}}}
\put(4951,-4636){\makebox(0,0)[lb]{\smash{{\SetFigFont{8}{9.6}{\rmdefault}{\mddefault}{\updefault}{\color[rgb]{0,0,0}$L_1$}%
}}}}
\put(6451,-5611){\makebox(0,0)[lb]{\smash{{\SetFigFont{8}{9.6}{\rmdefault}{\mddefault}{\updefault}{\color[rgb]{0,0,0}$L_2$}%
}}}}
\put(5851,-3586){\makebox(0,0)[lb]{\smash{{\SetFigFont{8}{9.6}{\rmdefault}{\mddefault}{\updefault}{\color[rgb]{0,0,0}Critical value}%
}}}}
\put(5776,-4636){\makebox(0,0)[lb]{\smash{{\SetFigFont{8}{9.6}{\rmdefault}{\mddefault}{\updefault}{\color[rgb]{0,0,0}The critical value $l$}%
}}}}
\end{picture}%

%% file: handles.pstex_t
\begin{picture}(0,0)%
\includegraphics{handles.pstex}%
\end{picture}%
\setlength{\unitlength}{2763sp}%
\begingroup\makeatletter\ifx\SetFigFont\undefined%
\gdef\SetFigFont#1#2#3#4#5{%
  \reset@font\fontsize{#1}{#2pt}%
  \fontfamily{#3}\fontseries{#4}\fontshape{#5}%
  \selectfont}%
\fi\endgroup%
\begin{picture}(7734,4795)(2089,-5744)
\put(6226,-4261){\makebox(0,0)[lb]{\smash{{\SetFigFont{8}{9.6}{\rmdefault}{\mddefault}{\updefault}{\color[rgb]{0,0,0}$H_3$}%
}}}}
\put(5251,-4336){\makebox(0,0)[lb]{\smash{{\SetFigFont{8}{9.6}{\rmdefault}{\mddefault}{\updefault}{\color[rgb]{0,0,0}$H_2$}%
}}}}
\put(4876,-3586){\makebox(0,0)[lb]{\smash{{\SetFigFont{8}{9.6}{\rmdefault}{\mddefault}{\updefault}{\color[rgb]{0,0,0}$H_1$}%
}}}}
\put(3076,-5686){\makebox(0,0)[lb]{\smash{{\SetFigFont{8}{9.6}{\rmdefault}{\mddefault}{\updefault}{\color[rgb]{0,0,0}$\phi = -1$}%
}}}}
\end{picture}%

%% file: lefschetzsym.bbl
\providecommand{\bysame}{\leavevmode\hbox to3em{\hrulefill}\thinspace}
\providecommand{\href}[2]{#2}
\begin{thebibliography}{10}

\bibitem{AbbondandoloSchwarz:noteonfloerloop}
A.~Abbondandolo and M.~Schwarz, \emph{Note on {F}loer homology and loop space
  homology.}, Morse theoretic methods in nonlinear analysis and in symplectic
  topology, NATO Sci. Ser. II Math. Phys. Chem. \textbf{217} (2006), 75--108.

\bibitem{abbondandoloschwarz:cotangentloopproduct}
\bysame, \emph{{Floer homology of cotangent bundles and the loop product}},
  (2008), 1--137, \mbox{arXiv:0810.1995}.

\bibitem{SeidelAbouzaid:viterbo}
M.~Abouzaid and P.~Seidel, \emph{An open string analogue of {V}iterbo
  functoriality.},  (2007), 1--74, \mbox{arXiv:0712.3177}.

\bibitem{BEHWZ:compactnessfieldtheory}
F.~Bourgeois, Y.~Eliashberg, H.~Hofer, K.~Wysocki, and E.~Zehnder,
  \emph{Compactness results in symplectic field theory}, Geom.Topol. \textbf{7}
  (2003), 799--888, \mbox{arXiv:SG/0308183}.

\bibitem{Brieskorn:sphere}
E.~Brieskorn, \emph{{Beispiele zur Differentialtopologie von
  Singularit\"aten.}}, Invent. Math. \textbf{2} (1966), 1--14.

\bibitem{DimcaChoudary:complexhypersurfaces}
A.~D.~R. Choudary and A.~Dimca, \emph{Complex hypersurfaces diffeomorphic to
  affine spaces.}, Kodai Math. J. \textbf{17, no. 2} (1994), 171--178.

\bibitem{Cieliebak:handleattach}
K.~Cieliebak, \emph{Handle attaching in symplectic homology and the chord
  conjecture}, J. Eur. Math. Soc. (JEMS) \textbf{4} (2002), 115--142.

\bibitem{CieliebakEliashberg:symplecticgeomofsteinmflds}
K.~Cieliebak and Y.~Eliashberg, \emph{Symplectic geometry of {S}tein
  manifolds}, In preparation.

\bibitem{CieliebakFloerHoferWysocki:SymhomIIApplications}
K.~Cieliebak, A.~Floer, H.~Hofer, and K.~Wysocki, \emph{Applications of
  symplectic homology {II}:stability of the action spectrum}, Math. Z
  \textbf{223} (1996), 27--45.

\bibitem{Eliashberg:steintopology}
Y.~Eliashberg, \emph{Topological characterisation of {S}tein manifolds of
  dimension $> 2$}, Internat. J. Math. \textbf{1} (1990), 29--46.

\bibitem{Eliashberg:symplecticgeometryofplushfns}
\bysame, \emph{Symplectic geometry of plurisubharmonic functions, notes by {M}.
  {A}breu, in: Gauge theory and symplectic geometry ({M}ontreal 1995)}, NATO
  Adv. Sci. Inst. Ser. C Math. Phys. Sci. \textbf{488} (1997), 49--67.

\bibitem{FHS:transversalitysymplectic}
A.~Floer, H.~Hofer, and D.~Salamon, \emph{Transversality in elliptic {M}orse
  theory for the symplectic action}, Duke Math.J. \textbf{80} (1995), 251--292.

\bibitem{FukayaSeidelSmith:exactlagragiancotangent}
K.~Fukaya, P.~Seidel, and I.~Smith, \emph{Exact {L}agrangian submanifolds in
  simply-connected cotangent bundles},  (2007), \mbox{arXiv:SG/0701783}.

\bibitem{Gompf:handlebody}
R.~Gompf, \emph{Handlebody construction of {S}tein surfaces.}, Ann. of Math.
  (2) \textbf{148} (1998), 619--693.

\bibitem{Hermann:holomorphic}
D.~Hermann, \emph{Holomorphic curves and {H}amiltonian systems in an open set
  with restricted contact-type boundary,}, Duke. Math. J \textbf{103} (2000),
  335--374.

\bibitem{Kaliman:eisenman}
S.~Kaliman, \emph{Exotic analytic structures and {E}isenman intrinsic
  measures.}, Israel J. Math. \textbf{88} (1994), 411--423.

\bibitem{McLean:thesis}
M.~McLean, \emph{Ph{D} thesis}, Cambridge (2008), 1--106.

\bibitem{Oancea:survey}
A.~Oancea, \emph{A survey of {F}loer homology for manifolds with contact type
  boundary or symplectic homology}, Ensaios Mat. \textbf{7} (2004),
  \mbox{arXiv:SG/0403377}.

\bibitem{Oancea:kunneth}
\bysame, \emph{The {K}{\"u}nneth formula in {F}loer homology for manifolds with
  restricted contact type boundary}, Math. Ann. \textbf{334} (2006), 65--89,
  \mbox{arXiv:SG/0403376}.

\bibitem{Ramanujam:affineplane}
C.~Ramanujam, \emph{A topological characterisation of the affine plane as an
  algebraic variety}, Ann. of Math. (2) \textbf{94} (1971), 69--88.

\bibitem{Seidel:hochschildhomology}
P.~Seidel, \emph{Symplectic homology as {H}ochschild homology},
  \mbox{arXiv:SG/0609037}.

\bibitem{Seidel:longexactsequence}
\bysame, \emph{A long exact sequence for symplectic {F}loer cohomology},
  Topology \textbf{42} (2003), 1003--1063, \mbox{arXiv:SG/0105186}.

\bibitem{Seidel:fukayalefschetz}
\bysame, \emph{Fukaya categories and {P}icard-{L}efschetz theory}, to appear in
  the ETH Lecture Notes series of the European Math. Soc. (2006).

\bibitem{Seidel:biasedview}
\bysame, \emph{A biased view of symplectic cohomology}, Current Developments in
  Mathematics \textbf{2006} (2008), 211--253.

\bibitem{SS:rama}
P.~Seidel and I.~Smith, \emph{The symplectic topology of {R}amanujam's
  surface}, Comment. Math. Helv. \textbf{80} (2005), 859--881,
  \mbox{arXiv:AG/0411601}.

\bibitem{tDP:1989moh}
T.~tom Dieck and T.~Petrie, \emph{The {A}bhyankar-{M}oh problem in dimension
  3}, Lecture Notes in Math. \textbf{1375} (1989), 48--59.

\bibitem{tDP:1990conaff}
\bysame, \emph{Contractible affine surfaces of {K}odaira dimension one}, Japan.
  J. Math. (N.S.) \textbf{16} (1990), 147--169.

\bibitem{Ustilovsky:infinitecontact}
I.~Ustilovsky, \emph{Infinitely many contact structures on ${S}^{4m+1}$},
  Internat. Math. Res. Notices \textbf{14} (1999), 781--791.

\bibitem{Viterbo:functorsandcomputations}
C.~Viterbo, \emph{Functors and computations in {F}loer homology with
  applications, part {I}.}, Geom. Funct. Anal. \textbf{9} (1999), 985--1033.

\bibitem{Zaidenberg:1998exot}
M.~Zaidenberg, \emph{Lectures on exotic algebraic structures on affine spaces},
   (1998), 1--60, \mbox{arXiv:AG/9801075}.

\end{thebibliography}
